\documentclass[a4paper,12pt]{article}
\usepackage[a4paper,centering,scale=.7]{geometry}
\newif\ifShowMarginPar

\ShowMarginParfalse % or

\usepackage{amsmath,amssymb}
\usepackage{amsthm}
\usepackage{graphicx}

 \usepackage{pst-all}
\newcommand{\Z}{\mathbb Z}
\newcommand{\N}{\mathbb N}
\newcommand{\R}{\mathbb R}

\renewcommand{\P}{\mathbb P}
\newcommand{\E}{\mathbb E}
\newcommand{\bu}{\mathbf{u}}
\newcommand{\bv}{\mathbf{v}}
\newcommand{\bz}{\mathbf{z}}
\newcommand{\bw}{\mathbf{w}}
\newcommand{\bx}{\mathbf{x}}
\newcommand{\by}{\mathbf{y}}

\newcommand{\prob}{\xrightarrow{\P}}

\newtheorem{theorem}{Theorem}[section]
\newtheorem{lemma}[theorem]{Lemma}%[section]
\newtheorem{remark}[theorem]{Remark}%[section]

\newtheorem{prop}[theorem]{Proposition}%[section]
\newtheorem{cor}{Corollary}[theorem]%[section]

\title{Random directed forest and the Brownian web}
\author{
{Rahul Roy, Kumarjit Saha and Anish Sarkar }
\footnote{E-Mail: {\tt rahul, kumarjit10r, anish@isid.ac.in}}\\
{\it Indian Statistical Institute, New Delhi}}

\date{}

\begin{document}

\maketitle

\begin{abstract}
Consider the $d$ dimensional lattice $\Z^d$ where each vertex is \textit{open} or
\textit{closed} with probability $p$ or $1-p$ respectively. An open vertex $\bu := (\bu(1),
\bu(2),\dotsc,\bu(d))$ is connected by an edge to another open vertex which has the minimum
$L_1$ distance among all the open vertices $\bx$ with $\bx(d)>\bu(d)$. It is shown that this random
graph is a tree almost surely for $d=2$ and $3$ and it is an infinite collection of disjoint
trees for $d\geq 4$. In addition, for $d=2$, we show that when properly scaled, the family of
its paths converge in distribution to the Brownian web.
\end{abstract}

\vspace{0.1in}

\vspace{0.1in}

%\input{Intro.tex}
%\input{1ASRR2007.tex}
% !TEX program = XeLaTeX
% !TEX root = ./DSFrf.tex

\section{Introduction}
Let $  {\cal P} $ % = \{{\bf x}_1, {\bf x}_2, \dotsc\}$
be the points of a Poisson point process
on $\R^d$ of intensity $1$. For each  $ {\bf x} \in  {\cal P}$
let $ h({\bf x}) \in {\cal P}  $ be the Poisson point in the half-space $ \{
{\bf u} : {\bf u}(d) > {\bf x} (d) \} $
which has the minimum Euclidean distance from $ {\bf x} $, where $ {\bf v}(j) $ denotes the $j$ th  co-ordinate
of $  {\bf v} \in \R^d$.
%  let $h({\bf x}) \in \Xi$ be such that (i) ${\bf x}(d) < h({\bf x})(d)$ and
% (ii) $|| {\bf x} - h({\bf x})||_2 < ||{\bf x} - {\bf y}||_2$ for any ${\bf y}
% \in \Xi$ with ${\bf x}(d) < {\bf y}(d)$, where $|| \cdot ||_p$ denotes the $L_p$
% metric on $\mathbb R^d$.
%This point $h({\bf x})$ is almost surely unique.
The
directed spanning forest (DSF) is the random graph with vertex set ${\cal P}$
and edge
set $ \{ \langle{\bf x}, h({\bf x}) \rangle : {\bf x} \in {\cal P}\}$.
\ifShowMarginPar
\marginpar{Notation ${\cal P}$}
\fi
The study of the directed spanning forest (DSF) was initiated by Baccelli
{\it et al.}\/ \cite{BB07}. Coupier {\it et al.}\/ \cite{CT11}
proved that for $d=2$
the DSF is a tree almost surely. Ferrari {\it et al.}\/ \cite{FLT04} also
studied a directed random graph on a Poisson point process, however, the mechanism
used to construct edges in that model incorporates more independence than  is available in the
DSF. They proved that their random graph is a connected tree in dimensions $2$ and
$3$, and a forest in dimensions $4$ and more.

A similar construction,  like the DSF arising from a Poisson point process, can
be made from vertices of the integer lattice.
Let $\{U_{\bv}:\bv \in \Z^{d}\}$ be a collection of i.i.d. uniform $(0,1)$
random variables. Fix $0<p<1$ and let $V:=\{\bv \in \Z^{d} : U_{\bv} < p\}$ be the set
of {\it open}\/ vertices of $\Z^d$.
Given $\bu \in \Z^{d}$, let $\bv \in V$ be such that
\ifShowMarginPar
\marginpar{Notation $U_{\bv}$}
\fi
\begin{enumerate}
\item $\bu(d)< \bv(d)$,
\item there does not exist any $\bw \in V$ with $\bw(d) > \bu(d)$ such that
$||\bu - \bw||_{1} < ||\bu - \bv||_{1}$, and
\item for all $\bw\in V$ with $\bw(d)>\bu(d)$ and $||\bu - \bw||_{1} = ||\bu -
\bv||_{1}$ we have $U_{\bv} \leq U_{\bw}$.
\end{enumerate}
Here and henceforth  $  ||\bu ||_{1} $ denotes  the $L_1$ norm of $\bu$ on $ \R^d$.
Such a $\bv$ is almost surely unique and clearly, is a function of $ \bu $ and
$ {\bf W} := \{ U_{\bw}: \bw \in \Z^{d},   \bw(d) > \bu(d) \} $. We denote it by 
$h(\bu, {\bf W})$.
We will drop the second argument in $ h$ for the time being.
Let $  \langle \bu,h(\bu)  \rangle$ be the edge joining
$\bu$ and $h(\bu)$ and let $E$ denote the edge set given by,
\begin{equation*}
%\label{TreeStep}
E := \{ \langle \bu,h(\bu) \rangle:\bu \in V\}.
\end{equation*}
\ifShowMarginPar
\marginpar{Notation $h ( \bu, {\bf X})$}
\fi
In this paper, we study the undirected random graph $ G:= (V, E)$, which we will refer to as
the {\em discrete DSF} henceforth.
% \documentclass[a4paper,12pt]{article}
% \usepackage{pst-all}
% \begin{document}

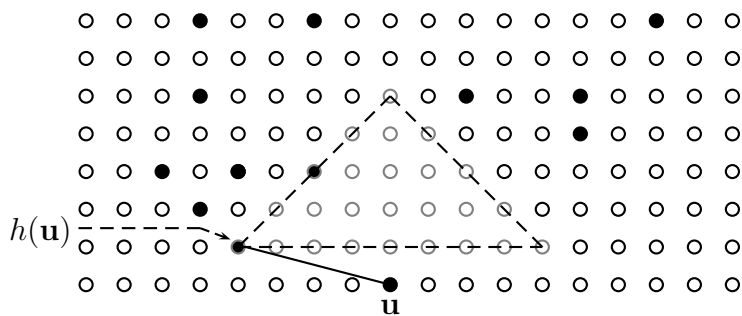
\begin{figure}[!htb]
\centering
%\begin{center}\leavevmode
%includegraphics[width=.5\textwidth]{hist1}

\begin{pspicture}(0,-0.45)(9,3.8)
\pscircle(0,0){.1}
\pscircle(0.5,0){.1}
\pscircle(1,0){.1}
\pscircle(1.5,0){.1}
\pscircle(2,0){.1}
\pscircle(2.5,0){.1}
\pscircle(3,0){.1}
\pscircle(3.5,0){.1}
\pscircle[fillcolor=black,fillstyle=solid](4,0){.1}
\pscircle(4.5,0){.1}
\pscircle(5,0){.1}
\pscircle(5.5,0){.1}
\pscircle(6,0){.1}
\pscircle(6.5,0){.1}
\pscircle(7,0){.1}
\pscircle(7.5,0){.1}
\pscircle(8,0){.1}
\pscircle(8.5,0){.1}

\pscircle(0,0.5){.1}
\pscircle(0.5,0.5){.1}
\pscircle(1,0.5){.1}
\pscircle(1.5,0.5){.1}
\pscircle[linecolor=gray,fillcolor=black,fillstyle=solid](2,0.5){.1}
\pscircle[linecolor=gray](2.5,0.5){.1}
\pscircle[linecolor=gray](3,0.5){.1}
\pscircle[linecolor=gray](3.5,0.5){.1}
\pscircle[linecolor=gray](4,0.5){.1}
\pscircle[linecolor=gray](4.5,0.5){.1}
\pscircle[linecolor=gray](5,0.5){.1}
\pscircle[linecolor=gray](5.5,0.5){.1}
\pscircle[linecolor=gray](6,0.5){.1}
\pscircle(6.5,0.5){.1}
\pscircle(7,0.5){.1}
\pscircle(7.5,0.5){.1}
\pscircle(8,0.5){.1}
\pscircle(8.5,0.5){.1}

\pscircle(0,1){.1}
\pscircle(0.5,1){.1}
\pscircle(1,1){.1}
\pscircle[fillcolor=black,fillstyle=solid](1.5,1){.1}
\pscircle(2,1){.1}
% \pscircle[linecolor=red](2.5,1){.1}
% \pscircle[linecolor=red](3,1){.1}
% \pscircle[linecolor=red](3.5,1){.1}
% \pscircle[linecolor=red](4,1){.1}
% \pscircle[linecolor=red](4.5,1){.1}
% \pscircle[linecolor=red](5,1){.1}
% \pscircle[linecolor=red](5.5,1){.1}
\pscircle[linecolor=gray](2.5,1){.1}
\pscircle[linecolor=gray](3,1){.1}
\pscircle[linecolor=gray](3.5,1){.1}
\pscircle[linecolor=gray](4,1){.1}
\pscircle[linecolor=gray](4.5,1){.1}
\pscircle[linecolor=gray](5,1){.1}
\pscircle[linecolor=gray](5.5,1){.1}
\pscircle(6,1){.1}
\pscircle(6.5,1){.1}
\pscircle(7,1){.1}
\pscircle(7.5,1){.1}
\pscircle(8,1){.1}
\pscircle(8.5,1){.1}

\pscircle(0,1.5){.1}
\pscircle(0.5,1.5){.1}
\pscircle[fillcolor=black,fillstyle=solid](1,1.5){.1}
\pscircle(1.5,1.5){.1}
\pscircle(2,1.5){.1}
\pscircle(2.5,1.5){.1}
% \pscircle[linecolor=red](3,1.5){.1}
% \pscircle[linecolor=red](3.5,1.5){.1}
% \pscircle[linecolor=red](4,1.5){.1}
% \pscircle[linecolor=red](4.5,1.5){.1}
% \pscircle[linecolor=red](5,1.5){.1}
\pscircle[linecolor=gray,fillcolor=black,fillstyle=solid](3,1.5){.1}
\pscircle[linecolor=gray](3.5,1.5){.1}
\pscircle[linecolor=gray](4,1.5){.1}
\pscircle[linecolor=gray](4.5,1.5){.1}
\pscircle[linecolor=gray](5,1.5){.1}
\pscircle(5.5,1.5){.1}
\pscircle(6,1.5){.1}
\pscircle(6.5,1.5){.1}
\pscircle(7,1.5){.1}
\pscircle(7.5,1.5){.1}
\pscircle(8,1.5){.1}
\pscircle(8.5,1.5){.1}

\pscircle(0,2){.1}
\pscircle(0.5,2){.1}
\pscircle(1,2){.1}
\pscircle(1.5,2){.1}
\pscircle(2,2){.1}
\pscircle(2.5,2){.1}
\pscircle(3,2){.1}
% \pscircle[linecolor=red](3.5,2){.1}
% \pscircle[linecolor=red](4,2){.1}
% \pscircle[linecolor=red](4.5,2){.1}
\pscircle[linecolor=gray](3.5,2){.1}
\pscircle[linecolor=gray](4,2){.1}
\pscircle[linecolor=gray](4.5,2){.1}
\pscircle(5,2){.1}
\pscircle(5.5,2){.1}
\pscircle(6,2){.1}
\pscircle[fillcolor=black,fillstyle=solid](6.5,2){.1}
\pscircle(7,2){.1}
\pscircle(7.5,2){.1}
\pscircle(8,2){.1}
\pscircle(8.5,2){.1}

\pscircle(0,2.5){.1}
\pscircle(0.5,2.5){.1}
\pscircle(1,2.5){.1}
\pscircle[fillcolor=black,fillstyle=solid](1.5,2.5){.1}
\pscircle(2,2.5){.1}
\pscircle(2.5,2.5){.1}
\pscircle(3,2.5){.1}
\pscircle(3.5,2.5){.1}
\pscircle[linecolor=gray](4,2.5){.1}
\pscircle(4.5,2.5){.1}
\pscircle[fillcolor=black,fillstyle=solid](5,2.5){.1}
\pscircle(5.5,2.5){.1}
\pscircle(6,2.5){.1}
\pscircle[fillcolor=black,fillstyle=solid](6.5,2.5){.1}
\pscircle(7,2.5){.1}
\pscircle(7.5,2.5){.1}
\pscircle(8,2.5){.1}
\pscircle(8.5,2.5){.1}

\pscircle(0,3){.1}
\pscircle(0.5,3){.1}
\pscircle(1,3){.1}
\pscircle(1.5,3){.1}
\pscircle(2,3){.1}
\pscircle(2.5,3){.1}
\pscircle(3,3){.1}
\pscircle(3.5,3){.1}
\pscircle(4,3){.1}
\pscircle(4.5,3){.1}
\pscircle(5,3){.1}
\pscircle(5.5,3){.1}
\pscircle(6,3){.1}
\pscircle(6.5,3){.1}
\pscircle(7,3){.1}
\pscircle(7.5,3){.1}
\pscircle(8,3){.1}
\pscircle(8.5,3){.1}

\pscircle(0,3.5){.1}
\pscircle(0.5,3.5){.1}
\pscircle(1,3.5){.1}
\pscircle[fillcolor=black,fillstyle=solid](1.5,3.5){.1}
\pscircle(2,3.5){.1}
\pscircle(2.5,3.5){.1}
\pscircle[fillcolor=black,fillstyle=solid](3,3.5){.1}
\pscircle(3.5,3.5){.1}
\pscircle(4,3.5){.1}
\pscircle(4.5,3.5){.1}
\pscircle(5,3.5){.1}
\pscircle(5.5,3.5){.1}
\pscircle(6,3.5){.1}
\pscircle(6.5,3.5){.1}
\pscircle(7,3.5){.1}
\pscircle[fillcolor=black,fillstyle=solid](7.5,3.5){.1}
\pscircle(8,3.5){.1}
\pscircle(8.5,3.5){.1}

\pscircle[fillcolor=black,fillstyle=solid](2,1.5){.1}
\psline(4,0)(2.1,0.5)
\rput(4,-0.3){$\mathbf{u}$}
\rput(-0.6,0.7){$h(\mathbf{u})$}
\psline[linestyle=dashed](2,0.5)(6,0.5)
\psline[linestyle=dashed](2,0.5)(4,2.5)
\psline[linestyle=dashed](6,0.5)(4,2.5)
\psline[linestyle=dashed](-0.1,0.75)(1.5,0.75)
\psline[linestyle=dashed]{>}(1.9,0.6)(1.5,0.75)
\end{pspicture}
\caption{The construction of ${h(\mathbf{u})}$ from ${\bf u}$ on $\Z^2$. The shaded points are open,
while the others are closed. Note that in order to get
${h(\mathbf{u})}$ from ${\bf u}$, we require information on the values
of the uniform random variables of the gray
vertices.}
\label{TreeStep1}
%\end{center}
\end{figure}

%\end{document}

%\end{document}

Similar models of random graphs are known in the physics literature as drainage
networks (see Scheidegger \cite{S67})
and have been studied  extensively
(see Rodr\'{i}guez-Iturbe {\it et al.}\/ \cite{RR97}).
Mathematically, for similar discrete processes but with a condition for constructing
edges which allows more independence, the dichotomy in dimensions of having
a single connected tree vis-a-vis a forest has been studied (see Gangopadhyay
{\it et al.}\/ \cite{GRS04},
Athreya {\it et al.}\/ \cite{ARS08}). The graph studied in  \cite{GRS04} 
connected an open vertex $\bu$ to the vertex $h(\bu)$ with $h(\bu)$ being the 
nearest open vertex in $\{\bw : \bw(d) = \bu(d) + 1\}$, with the vertex being 
chosen with uniform probability in case there are more than one nearest open 
vertex.
This construction  immediately leads to a Markovian
analysis which is exploited in \cite{GRS04} to obtain the tree/forest dichotomy.
However the DSF model considered here has to take care of a `history' set arising from the paths constructed in the past. The Markovian structure is thus  obtained through regeneration times. Moreover to obtain the dichotomy requires information on the tail of the distribution of the regeneration time which we do here through coupling and auxiliary results on renewal processes.

Our paper may  also be viewed as an extension,
albeit in the discrete setting, of the result of
Coupier {\it et al.}\/ \cite{CT11} to any dimension. Our proof
is different from that of \cite{CT11}; while their argument is
percolation theoretic and  crucially depends on the planarity
of $\R^2$, our argument exploits the  Markovian structure of the DSF,
thereby allowing us to extend the result to any dimension. The difficulty of carrying our analysis in the continuous model studied in \cite{CT11} is that there is no obvious extension of regeneration time as considered here.

\begin{theorem}
\label{r-treesforest}
For $d=2$ and $d=3$ the random graph $G$ is connected almost surely
and consists of a single tree
while for $d\geq 4$, it is a disconnected forest containing infinitely many distinct connected components, each connected component
being an infinite tree almost surely.
\end{theorem}

Our second result in this paper
is the convergence of the random graph $G$ for $d=2$,
under a suitable diffusive scaling, to the Brownian web.
The standard Brownian web originated in
the work of Arratia \cite{A79}, \cite{A81}
as the scaling limit of the voter model on $\Z$. It arises naturally as the diffusive
scaling limit of the coalescing simple random walk paths starting from every point
on the space-time lattice. We can thus think of the Brownian web as a collection
of one-dimensional coalescing Brownian motions starting from every point in the space time
plane $\R^2$. Detailed analysis of the Brownian web was carried out in T\'{o}th {\it et al.}\/
\cite{TW98}. Later Fontes {\it et al.}\/ \cite{FINR04} introduced a framework in
which the Brownian web is realized as a random variable taking values in a Polish space.
We recall relevant details from Fontes {\it et al.}\/ \cite{FINR04}.

% I think here recalling from Fontes {\it et al.}\/ \cite{FINR04} is enough
Let $\R^{2}_c$ denote the completion of the space time plane $\R^2$ with
respect to the metric
\begin{equation*}
 \rho((x_1,t_1),(x_2,t_2)) := |\tanh(t_1)-\tanh(t_2)|\vee \Bigl| \frac{\tanh(x_1)}{1+|t_1|}
 -\frac{\tanh(x_2)}{1+|t_2|} \Bigr|.
\end{equation*}
% I think it is good to include this line,
As a topological space $\R^{2}_c$ can be identified with the
continuous image of $[-\infty,\infty]^2$ under a map that identifies the line
$[-\infty,\infty]\times\{\infty\}$ with the point $(\ast,\infty)$, and the line
$[-\infty,\infty]\times\{-\infty\}$ with the point $(\ast,-\infty)$.
%\textcolor{red}{begin KS}Let $\Pi$ be the space of all paths in $\R^{2}_c$ of the form
%$\{(\pi(t),t):t\geq \sigma_{\pi}\}$
%where $\sigma_{\pi}\in [-\infty,\infty]$ denotes the starting time for path
%$\pi$. \textcolor{red}{end KS}
% I think this part we should write in this way,
A path $\pi$ in $\R^{2}_c$ with starting time $\sigma_{\pi}\in [-\infty,\infty]$
is a mapping $\pi :[\sigma_{\pi},\infty]\rightarrow [-\infty,\infty]$ such that
$\pi(\infty)= \ast$ and, when $\sigma_\pi = -\infty$, $\pi(-\infty)= \ast$. 
Also $t\rightarrow (\pi(t),t)$ is a continuous
map from $[\sigma_{\pi},\infty]$ to $(\R^{2}_c,\rho)$.
We then define $\Pi$ to be the space of all paths in $\R^{2}_c$ with all possible starting times in $[-\infty,\infty]$.
The following metric, for $\pi_1,\pi_2\in \Pi$
\begin{equation*}
d_{\Pi} (\pi_1,\pi_2) := |\tanh(\sigma_{\pi_1})-\tanh(\sigma_{\pi_2})|\vee\sup_{t\geq
\sigma_{\pi_1}\wedge
\sigma_{\pi_2}} \Bigl|\frac{\tanh(\pi_1(t\vee\sigma_{\pi_1}))}{1+|t|}-\frac{
\tanh(\pi_2(t\vee\sigma_{\pi_2}))}{1+|t|}\Bigr|
\end{equation*}
%where $\sigma_{\pi}$ denotes the starting point of the path $\pi$.
makes $\Pi$ a complete, separable metric space. Convergence in this
metric can be described as locally uniform convergence of paths as
well as convergence of starting times.
Let ${\cal H}$ be the space of compact subsets of $(\Pi,d_{\Pi})$ equipped with
the Hausdorff metric $d_{{\cal H}}$ given by,
\begin{equation*}
d_{{\cal H}}(K_1,K_2) := \sup_{\pi_1 \in K_1} \inf_{\pi_2 \in
K_2}d_{ \Pi} (\pi_1,\pi_2)\vee
\sup_{\pi_2 \in K_2} \inf_{\pi_1 \in K_1} d_{\Pi} (\pi_1,\pi_2).
\end{equation*}
The space $({\cal H},d_{{\cal H}})$ is a complete separable metric space. Let
$B_{{\cal H}}$ be the Borel  $\sigma-$algebra on the metric space $({\cal H},d_{{\cal H}})$.
The Brownian web ${\cal W}$ is an $({\cal H},B_{{\cal H}})$ valued random
variable.
\ifShowMarginPar
\marginpar{Notation path $\pi$, $\sigma_{\pi}$, ${\cal W}$}
\fi

Ferrari {\it et al.}\/ \cite{FFW05} have shown that, for $d=2$,  the random graph on
the Poisson points introduced by \cite{FLT04},
converges to a Brownian web under a suitable diffusive scaling.
Coletti {\it et al.}\/ \cite{CFD09}
have a similar result for the discrete random graph studied in Gangopadhyay {\it et al.}\/ \cite{GRS04}.
Baccelli {\it et al.}\/ \cite{BB07} have shown that scaled paths of the
successive ancestors in the DSF converges weakly to the Brownian motion
and also conjectured that the scaling limit of the DSF is the  Brownian web. 

Our work here differs from that of
\cite{FFW05} and \cite{CFD09} in that we need to obtain the Brownian web as a 
limit of a Markov process defined through regeneration times, while in the 
earlier work correlation inequalities like the FKG inequality could be used 
because every step of the paths constructing their model had i.i.d. increments. 
The method we employ requires us to control the size of the region 
surveyed to obtain the regeneration time of a process starting from a single 
vertex. Also using a martingale constructed via the joint 
regeneration times of processes starting from two distinct starting points we 
estimate the tail probability of the coalescing time.
This method which we present here can be used in both \cite{FFW05} and 
\cite{CFD09} to obtain their results without invoking correlation inequalities. 
Also for the model considered by \cite{BB07} and
\cite{CT11}, if a 
suitable `pseudo-regeneration time' of joint processes is defined and there is 
a control on the size of the region explored to obtain such pseudo-regeneration 
times, then our approach should yield the convergence to the Brownian web. In 
addition such pseudo-regeneration times should also yield the geometric 
structure of the DSF in dimensions 3 or more.

%\marginpar{check the new definitions of $\chi$ and $\chi_n$}

From a vertex $\bu \in \Z^2$, taking the edges
$\{ \langle h^{k-1}(\bu),h^{k}(\bu) \rangle : k \geq 1)\}$
\ifShowMarginPar
\marginpar{Notation $\pi^{\bu}$}
\fi(with $h^0(\bu) := \bu$ and $h^{k}(\bu) := h(h^{k-1}(\bu))$)
to be straight line segments we parametrize the path formed by these edges as
the piecewise linear function $\pi^{\bu} : [\bu (2), \infty) \to \mathbb R$
such that $\pi^{\bu}(h^{k}(\bu)(2)) := h^{k}(\bu)(1)$ for every $k \geq 0$ and
$\pi^{\bu}(t)$ is linear in the interval $[h^{k}(\bu)(2),h^{k+1}(\bu)(2)]$.
Define ${\cal X} := \{\pi^{\bu}:\bu \in V\}$.
\ifShowMarginPar
\marginpar{Notation ${\cal X}$}
\fi
For given  $\gamma , \sigma > 0$,  a path $\pi$ with starting time $
\sigma_{\pi}$ and for each $n \geq 1$,  the scaled path
$ \pi_n(\gamma,\sigma) : [\sigma_{\pi}/n^2\gamma, \infty] \to
[-\infty, \infty]$ is given by $\pi_n(\gamma,\sigma) (t) :=
\pi(n^2\gamma t)/n \sigma$. Thus, the scaled path $\pi_n(\gamma,\sigma)$ has
 the starting time $\sigma_{\pi_n(\gamma,\sigma) }= \sigma_{\pi}/ n^2\gamma$.
For each $ n \geq 1 $, let ${\cal X}_n(\gamma,\sigma) :=
\{\pi_n^{\bu}(\gamma,\sigma):\bu \in V\}$ be the collection of the scaled paths.
The closure  $\bar{\cal X}_n(\gamma,\sigma)$ of
${\cal X}_n(\gamma,\sigma)$ in $(\Pi,d_{\Pi})$ is a
$({\cal H},{\cal B}_{{\cal H}})$ valued random variable. We have
\ifShowMarginPar
\marginpar{Notation ${\cal X}_n$}
\fi
%
% $\pi_n^{\bu}$ be the scaled path $\pi_n^{\bu}(t) = \frac{\pi(n^2\gamma
% t)}{n\sigma}$  thus $\sigma(\pi_n^{\bu}) = \frac{\sigma_{\pi}}{n^2\gamma}$.

%
% $\pi_n^{\bu}: [\frac{\sigma_{\pi}}{n^2\gamma},
% \infty] \to
% [-\infty, \infty]$  = $

% \textcolor{red}{begin KS}
% Let
% $\pi^{\bu} :=\{(\phi^{\bu}(t),t):t\geq \bu (2)\}$, $\bu \in \mathbb Z^d$
%denote
% the graph of the path $\phi^{\bu}(t)$ and
% ${\cal X} := \bigcup_{\bu \in V}\pi^{\bu}$ the geometric
% representation of our
% random graph.
% For $\gamma, \sigma >0$,  we define the $n$ th order diffusive scaling of
% ${\cal X}$ by
% \begin{equation}
% \label{cont-chi}
% {\cal X} _n(\gamma,\sigma) :=
% \{{\displaystyle \bigcup_{t \geq \bu(t)}}
% (\frac{\phi^{\bu}(t)}{n\sigma},\frac{t}{n^2\gamma}): \bu \in V\}=
% \{(\frac{y(1)}{n\sigma},\frac{y(2)}{n^2\gamma}): (y(1),y(2))\in {\cal
% X}\}.
% \end{equation}
% We have
% \begin{theorem}
% \label{r-BW}
% For $d=2$, there exist $\sigma := \sigma(p)$ and $\gamma := \gamma(p)$ such
% that
% as $n\rightarrow \infty$,  ${\cal X}_n(\gamma,\sigma)$ converges weakly to the
% standard Brownian Web.
% \end{theorem}
% \textcolor{red}{end KS}
%\marginpar{check the new definitions of $\chi$ and $\chi_n$}

% This part I want to write in this way

% For notational simplicity let ${\cal X}$ denote the closure of
% ${\cal X} := \{\phi^{\bu}:\bu \in V\}$ in $(\Pi,d)$.
% For each $n\geq 1$ and for all
% $\gamma, \sigma >0$ the $n$ th order diffusively scaled geometric
% representation of our random graph given by ${\cal X} _n(\gamma,\sigma)$ is an
% $({\cal H},{\cal B}_{{\cal H}})$ valued random variable.

\begin{theorem}
\label{r-BW}
There exist $\sigma := \sigma(p)$ and $\gamma := \gamma(p)$ such
that as $n\rightarrow \infty$,  $\bar{{\cal X}}_n(\gamma,\sigma)$ converges weakly to the
standard Brownian Web ${\cal W}$ as $({\cal H},{\cal B}_{{\cal H}})$ valued random variables.
\end{theorem}
\begin{remark}
 \label{rem:Normconstants}
The scaling property of the Brownian web yields
$\bar{{\cal X}}_n(1,\sigma^\prime) \Rightarrow {\cal W}$ as $n\rightarrow \infty$
for  $\sigma^\prime := \sigma/\sqrt{\gamma}$.
\end{remark}

For the proof of Theorem \ref{r-treesforest}  we obtain a Markovian structure in
our model and define
suitable stopping times for this Markov process.
From these stopping times the process regenerates which allows us
to phrase the problem as a question of recurrence or transience
of the Markov chain. This we do by obtaining a martingale for $d=2$, using
a Lyapunov function technique for $d=3$ and a suitable coupling with a
random walk with independent steps for $d=4$.

The martingale obtained for $d=2$ and the fact that the
distributions of the stopping times have  exponentially decaying tails
are used to prove Theorem \ref{r-BW}.

Finally, although our results are obtained for the random graph constructed
by connecting edges between $L_1$ nearest open vertices, they should also
hold for the model constructed with the $L_2$ metric (see Remark 
\ref{rem:L2Statement} for more details).

The paper is structured as follows -- in the next section we construct the
paths of the graph $G$ starting from $k$ distinct vertices and obtain some properties of these paths.
In Section 3, we derive the martingale (for $d=2$) and also provide a method of approximation
of the paths by independent processes, which is used later to prove Theorem \ref{r-treesforest}
and Theorem \ref{r-BW}.
In Section 4 we prove Theorem \ref{r-treesforest} and in Section 5, we prove Theorem  \ref{r-BW}.

%\input{sec1rf.tex}
%\input{sec2rf.tex}
% !TEX root = ./DSFrf.tex
% !TEX program = XeLaTeX
\section{Construction of the process}
\label{sec:ConstructionOfProcess}
We first  detail a construction of the graph $G$ which is needed to bring out a Markovian
structure. Later we obtain a martingale for $d=2$ which is used in the next two
sections. Before proceeding further we fix some notation:
 for $ \bu \in \Z^d$ and $ r \in \Z $, let $ \mathbb{H} (r ) := \{  \bw \in \Z^d :  \bw (d) \leq r  \} $ be the half-space
and, for $ r > 0 $, let $ S^{+}( \bu, r ) := \{ \bw \in \Z^d : || \bu - \bw 
||_1 \leq r,
\bw(d) > \bu(d) \}  $ be the upper part of closed $L_1$ ball at $ \bu $ having 
radius $r$. As a convention we take $S^{+}( \bu, 0 ) := \emptyset$.

% for $ \bu \in \Z^d$ and $ r > 0 $,
% let $ S( \bu, r ) := \{ \bw \in \Z^d : || \bu - \bw ||_1 \leq r \} $ be the closed $L_1$ ball of
% radius $r$ and  $ \mathbb{H} (r ) := \{  \bw \in \Z^d :  \bw (d) \leq r  \} $ be the half-space.

\begin{figure}[!htb]
\centering
%\begin{center}\leavevmode
% Generated with LaTeXDraw 2.0.8
% Mon Sep 03 10:35:47 IST 2012
% \usepackage[usenames,dvipsnames]{pstricks}
% \usepackage{epsfig}
\scalebox{1} % Change this value to rescale the drawing.
{
\begin{pspicture}(0,-3.6)(15.2,3.60)
\pspolygon[linewidth=0.03,linestyle=dashed,dash=0.16cm 0.16cm,fillstyle=solid,fillcolor=lightgray](1.22,3.2245312)(0.02,2.0845313)(2.42,2.1045313)(2.36,2.1045313)
\pspolygon[linewidth=0.03,linestyle=dashed,dash=0.16cm 0.16cm,fillstyle=solid,fillcolor=lightgray](2.86,2.0845313)(2.88,2.1045313)(4.08,3.2645311)(5.24,2.0845313)
\psdots[dotsize=0.12](11.2,1.6845312)
\psdots[dotsize=0.12](4.82,2.5045311)
\usefont{T1}{ptm}{m}{n}
\rput(11.965625,1.9370313){$g_n(\bu^2)$}
\usefont{T1}{ptm}{m}{n}
\rput(5.782344,2.7570312){$g_n(\bu^1)$}
\psline[linewidth=0.028222222cm,linestyle=dashed,dash=0.17638889cm 0.10583334cm,arrowsize=0.05291667cm 2.0,arrowlength=1.4,arrowinset=0.4]{<-}(1.22,2.7370312)(2.0,3.3370314)
\psline[linewidth=0.028222222cm,linestyle=dashed,dash=0.17638889cm 0.10583334cm,arrowsize=0.05291667cm 2.0,arrowlength=1.4,arrowinset=0.4]{<-}(4.08,2.7370312)(3.5,3.3370314)
\usefont{T1}{ptm}{m}{n}
\rput(2.7745314,3.6370313){$\Delta_{n}(\bu^1,\bu^2)$}
\psline[linewidth=0.06cm,arrowsize=0.05291667cm 2.0,arrowlength=1.4,arrowinset=0.4]{->}(6.44,1.5870312)(6.46,0.02953125)
% \pspolygon[linewidth=0.03,linestyle=dashed,dash=0.16cm 0.16cm,fillstyle=solid,fillcolor=lightgray](11.22,-1.1354687)(10.02,-2.2954688)(12.42,-2.3154688)
% \pspolygon[linewidth=0.03,linestyle=dashed,dash=0.16cm 0.16cm,fillstyle=solid,fillcolor=lightgray](4.1,-1.8954687)(3.66,-2.3354688)(3.68,-2.2954688)(4.52,-2.3354688)
% \pspolygon[linewidth=0.03,linestyle=dashed,dash=0.16cm 0.16cm,fillstyle=solid,fillcolor=lightgray](1.66,-2.3554688)(1.66,-2.3354688)(1.24,-1.9554688)(0.8,-2.3354688)
\psline[linewidth=0.028222222cm,linestyle=dashed,dash=0.17638889cm 0.10583334cm,arrowsize=0.05291667cm 2.0,arrowlength=1.4,arrowinset=0.4]{<-}(11.22,-1.3354688)(9.2,-0.16296875)
\psline[linewidth=0.028222222cm,linestyle=dashed,dash=0.17638889cm 0.10583334cm,arrowsize=0.05291667cm 2.0,arrowlength=1.4,arrowinset=0.4]{<-}(4.1,-2.0954688)(7.2,-0.16296875)
\psline[linewidth=0.028222222cm,linestyle=dashed,dash=0.17638889cm 0.10583334cm,arrowsize=0.05291667cm 2.0,arrowlength=1.4,arrowinset=0.4]{<-}(1.24,-2.1554687)(7.0,-0.16296875)
\usefont{T1}{ptm}{m}{n}
\rput(8.374531,0.13703126){$\Delta_{n+1}(\bu^1,\bu^2)$}
\psdots[dotsize=0.12](11.2,-3.4754686)
\psdots[dotsize=0.12](4.82,-2.6554687)
\usefont{T1}{ptm}{m}{n}
\rput(11.845625,-3.6054688){$g_n(\bu^2)$}
\usefont{T1}{ptm}{m}{n}
\rput(6.9623437,-2.4654686){$g_{n+1}(\bu^1)=g_{n}(\bu^1)$}
\psdots[dotsize=0.12](12.4,-2.3354688)
\psline[linewidth=0.03cm,linestyle=dashed,dash=0.16cm 0.16cm](11.2,-1.1154687)(13.22,-3.0954688)
\psline[linewidth=0.03cm,linestyle=dashed,dash=0.16cm 0.16cm](11.22,-1.1154687)(9.2,-3.0754688)
\psline[linewidth=0.03cm,linestyle=dashed,dash=0.16cm 0.16cm](4.1,-1.8954687)(5.24,-3.0754688)
\psline[linewidth=0.03cm,linestyle=dashed,dash=0.16cm 0.16cm](4.1,-1.8954687)(2.84,-3.0954688)
\psline[linewidth=0.03cm,linestyle=dashed,dash=0.16cm 0.16cm](1.24,-1.9554688)(2.4,-3.0954688)
\psline[linewidth=0.03cm,linestyle=dashed,dash=0.16cm 0.16cm](1.24,-1.9354688)(0.0,-3.0554688)
\psline[linewidth=0.04cm](11.24,-3.4554687)(12.4,-2.3354688)
\usefont{T1}{ptm}{m}{n}
\rput(13.514531,-2.1654687){$g_{n+1}(\bu^2)$}
% \psdots[dotsize=0.12,fillstyle=solid,dotstyle=o](4.805,2.9004688)
% \psdots[dotsize=0.12,fillstyle=solid,dotstyle=o](4.805,3.3004687)
% \psdots[dotsize=0.12,fillstyle=solid,dotstyle=o](4.805,3.7404687)
% \psdots[dotsize=0.12,fillstyle=solid,dotstyle=o](11.205,2.1004686)
% \psdots[dotsize=0.12,fillstyle=solid,dotstyle=o](11.205,2.5204687)
% \psdots[dotsize=0.12,fillstyle=solid,dotstyle=o](11.205,2.9004688)
% \psdots[dotsize=0.12,fillstyle=solid,dotstyle=o](11.185,3.3004687)
% \psdots[dotsize=0.12,fillstyle=solid,dotstyle=o](11.185,3.7004688)
% \psdots[dotsize=0.12,fillstyle=solid,dotstyle=o](4.805,-0.99953127)
% \psdots[dotsize=0.12,fillstyle=solid,dotstyle=o](4.805,-0.59953123)
% \psdots[dotsize=0.12,fillstyle=solid,dotstyle=o](12.365,-1.9395312)
% \psdots[dotsize=0.12,fillstyle=solid,dotstyle=o](12.365,-1.4995313)
% \psdots[dotsize=0.12,fillstyle=solid,dotstyle=o](12.365,-1.0995313)
% \psdots[dotsize=0.12,fillstyle=solid,dotstyle=o](12.365,-0.69953126)
% \psdots[dotsize=0.12,fillstyle=solid,dotstyle=o](12.365,-0.27953124)
\pspolygon[linewidth=0.03,linestyle=dashed,dash=0.16cm 0.16cm,fillstyle=solid,fillcolor=lightgray](11.205,-1.1195313)(10.305,-1.9795313)(12.125,-1.9995313)
\pspolygon[linewidth=0.03,linestyle=dashed,dash=0.16cm 0.16cm,fillstyle=solid,fillcolor=lightgray](1.205,-1.9795313)(0.945,-2.1395314)(1.465,-2.1395311)
\pspolygon[linewidth=0.03,linestyle=dashed,dash=0.16cm 0.16cm,fillstyle=solid,fillcolor=lightgray](4.085,-1.9195312)(3.805,-2.1295312)(4.345,-2.1295312)
% \psdots[dotsize=0.12,fillstyle=solid,dotstyle=o](4.785,-2.1595314)
% \psdots[dotsize=0.12,fillstyle=solid,dotstyle=o](4.805,-1.7995312)
% \psdots[dotsize=0.12,fillstyle=solid,dotstyle=o](4.805,-1.3595313)
\end{pspicture}
}
\caption{The vertices
$ g_{n+1}(\bu^1),g_{n+1}(\bu^2)$ and the history region
$\Delta_{n+1}(\bu^1,\bu^2)$ when $W_n^{\text{move}}
= \{ g_n (\bu^2) \} $, $W_n^{\text{stay}}
= \{ g_n (\bu^1) \}. $}

\label{Fig:Markov3}
%\end{center}
\end{figure}
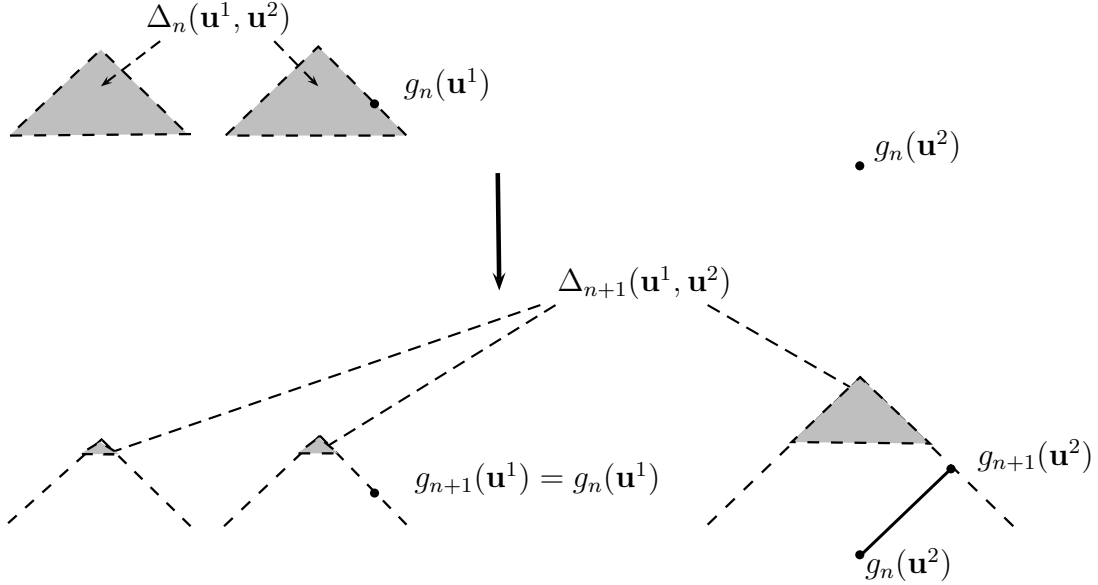

From $k$ ($k \geq 1$) vertices $\bu^{1}, \dotsc , \bu^{k} \in \Z^d$ with  $\bu^1(d) = \cdots = \bu^k(d)$,
we obtain the vertices $\{h^{n}(\bu^{i}): \; n\geq 0, \; 1\leq i \leq k\}$ as a stochastic process.
{\it Note here that the construction described below does not require  the vertices $\bu^{1}, \dotsc ,
\bu^{k} $ to be  open}\/.
The  vertices with the smallest $d$ th co-ordinate are allowed to move, while the others stay put
(see Figure \ref{Fig:TreeStep1} and \ref{Fig:Markov3}).
Each of these vertices explores a region in the half space `above' it to obtain the vertex to
which it moves. During this exploration a vertex may encounter regions which have
been already explored
by other vertices earlier. While the information for the region explored earlier is known, the information
about the freshly explored region is new and is obtained during the exploration process of the
vertices which are moving at that time. The region which has been explored till the $n$ th move
of the entire process and which is needed for the $n+1$ th move is called the {\it history region}\/
and  the information of the uniform random variables in the history region constitutes the {\em history}.
%and denoted by
% $\Delta_n = \Delta_n (\bu^{1}, \dotsc , \bu^{k})$ and $\{(\bw, U_{\bw}): \bw \in \Delta_n\}$
% constitutes the {\it history}\/ $H_n = H_{n} (\bu^{1}, \dotsc , \bu^{k})$.
Formally, let
\begin{itemize}
\item[(i)] $g_0(\bu^{i}) := \bu^{i}$ for all $1 \leq i \leq k$ and $r_0 := \bu^1 (d) $; %)  \min\{g_0 (\bu^{i})(d)
%:1\leq i\leq k\}$;
\item[(ii)]
$W_0^{\text{move}}:=  \{ \bu^{1}, \dotsc , \bu^{k}\} $ and  $W_0^{\text{stay}}:=   \emptyset $; % \{g_0(\bw): \bw \in,\; g_0(\bw)(d) = r_0\}$
%  and
% $W_0^{\text{stay}}:=  \{g_0(\bu^{1}), $ $\dotsc, g_0(\bu^{k})\}\setminus W_0^{\text{move}}$;
\item[(iii)] $ \Delta_0 = \Delta_0  (\bu^{1}, \dotsc , \bu^{k}) := \emptyset $  
and $\Psi_0: \emptyset \to [0,1]$ the empty function (see \cite{HS73}).

% $\Delta_0 = \Delta_0 (\bu^{1}, \dotsc , \bu^{k}):= \{\bw : \bw \in W_0^{\text{stay}}\} $, and
% %\item[(iv)]
% $H_0= H_{0} (\bu^{1}, \dotsc , \bu^{k}):= \{(\bw, x): \bw \in \Delta_0, \; x = U_{\bw}\}$.
\end{itemize}
Having obtained $g_{n}(\bu^{i})$, $r_n$, $W_n^{\text{move}}$, $W_n^{\text{stay}}$, $
\Delta_n$ and $\Psi_n$, for $ 1\leq i \leq k$,  we set
\ifShowMarginPar
\marginpar{Notation $g_n, W_n^{\text{move}}$ and $ W_n^{\text{stay}}$ }
\fi
\begin{itemize}
\item[(i)] $g_{n+1}(\bu) :=h(g_{n}(\bu)) \text{ for all } g_{n}(\bu) \in
W_n^{\text{move}}$ and
$g_{n+1}(\bv) := g_n(\bv) $ for all $g_{n}(\bv) \in W_n^{\text{stay}}$,  $r_{n+1}  :=
\min \{ g_{n+1} (\bu^{i})(d) : 1 \leq i \leq k \}$;
\item[(ii)]
$W_{n+1}^{\text{move}}:= \{g_{n+1}(\bw) : \bw \in \{ \bu^{1}, \dotsc , \bu^{k}\},\;
g_{n+1}(\bw) (d) = r_{n+1}\}$ and
$W_{n+1}^{\text{stay}}:= \{g_{n+1}(\bu^{1}), \dotsc, g_{n+1}(\bu^{k})\} \setminus
W_{n+1}^{\text{move}}$;
\item[(iii)] $\Delta_{n+1} = \Delta_{n+1} (\bu^{1}, \dotsc , \bu^{k}) :=  \bigl( \Delta_n \cup \cup_{ \bu \in  W_n^{\text{move}} }
S^+ ( \bu , || h (\bu)  - \bu ||_1 )  \bigr) \setminus  \mathbb{H} ( r_{n+1}  
)$ and
%\item[(iv)]
%$H_{n+1} = H_{n+1} (\bu^{1}, \dotsc , \bu^{k}) := (\Delta_{n+1}, \Psi_{n+1} )$ 
%where
$ \Psi_{n+1} : \Delta_{n+1} \to [0,1]$ is a map given by $ \Psi_{n+1} ( \bw) := 
U_{\bw } $ for $ \bw \in \Delta_{n+1}$, with $\Psi_{n+1} := \Psi_0$, the 
empty function, when $\Delta_{n+1} = \emptyset$.
%\{(\bw,x): \bw \in \Delta_{n+1}, \; x = U_{\bw}\}$.
\end{itemize}

% !TEX root = ./../DSFrf.tex
% !TEX program = XeLaTeX

\begin{figure}[!htb]
\centering
%\begin{center}\leavevmode
%includegraphics[width=.5\textwidth]{hist1}

% Generated with LaTeXDraw 2.0.8
% Mon Sep 03 09:53:56 IST 2012
% \usepackage[usenames,dvipsnames]{pstricks}
% \usepackage{epsfig}
% \usepackage{pst-grad} % For gradients
% \usepackage{pst-plot} % For axes
\scalebox{1} % Change this value to rescale the drawing.
{
\begin{pspicture}(0,-4.7)(14.070156,4.7)
\pspolygon[linewidth=0.04,linestyle=dashed,dash=0.16cm 0.16cm,fillstyle=solid,fillcolor=lightgray](2.34,4.7925)(0.0,2.4325)(4.76,2.4325)
\pspolygon[linewidth=0.04,linestyle=dashed,dash=0.16cm 0.16cm,fillstyle=solid,fillcolor=lightgray](11.16,3.1925)(10.34,2.4725)(11.96,2.4725)
\psdots[dotsize=0.12](5.16,2.0125)
\psdots[dotsize=0.12](12.36,2.0125)
\psline[linewidth=0.028222222cm,linestyle=dashed,dash=0.17638889cm 0.10583334cm,arrowsize=0.05291667cm 2.0,arrowlength=1.4,arrowinset=0.4]{<-}(3.38,2.755)(7.48,2.755)
\psline[linewidth=0.028222222cm,linestyle=dashed,dash=0.17638889cm 0.10583334cm,arrowsize=0.05291667cm 2.0,arrowlength=1.4,arrowinset=0.4]{<-}(10.98,2.755)(9.38,2.755)
\usefont{T1}{ptm}{m}{n}
\rput(8.4045315,2.955){$\Delta_n(\bu^1,\bu^2)$}
\usefont{T1}{ptm}{m}{n}
\rput(5.3523436,1.7025){$g_n(\bu^1)$}
\usefont{T1}{ptm}{m}{n}
\rput(12.604531,1.7025){$g_n(\bu^2)$}
\psdots[dotsize=0.12](5.14,-4.3875)
\psdots[dotsize=0.12](12.34,-4.3875)
\usefont{T1}{ptm}{m}{n}
\rput(5.3323436,-4.6975){$g_n(\bu^1)$}
\usefont{T1}{ptm}{m}{n}
\rput(12.815625,-4.6975){$g_n(\bu^2)$}
\pspolygon[linewidth=0.03,linestyle=dashed,dash=0.16cm 0.16cm,fillstyle=solid,fillcolor=lightgray](2.34,-1.6275)(1.14,-2.7675)(3.54,-2.7475)(3.48,-2.7475)
\pspolygon[linewidth=0.03,linestyle=dashed,dash=0.16cm 0.16cm,fillstyle=solid,fillcolor=lightgray](3.98,-2.7675)(4.0,-2.7475)(5.2,-1.5875)(6.36,-2.7675)
\psdots[dotsize=0.12](12.32,-3.1675)
\psdots[dotsize=0.12](5.94,-2.3475)
\psline[linewidth=0.04cm](5.14,-4.3475)(5.92,-2.3675)
\psline[linewidth=0.04cm](12.34,-4.3875)(12.32,-3.1875)
\psline[linewidth=0.04cm,linestyle=dashed,dash=0.16cm 0.16cm](5.18,-1.5875)(7.54,-3.9475)
\psline[linewidth=0.03cm,linestyle=dashed,dash=0.16cm 0.16cm](5.18,-1.5875)(2.66,-4.0275)
\psline[linewidth=0.028222222cm,linestyle=dashed,dash=0.17638889cm 0.10583334cm,arrowsize=0.05291667cm 2.0,arrowlength=1.4,arrowinset=0.4]{<-}(2.38,-2.045)(3.18,-1.545)
\psline[linewidth=0.028222222cm,linestyle=dashed,dash=0.17638889cm 0.10583334cm,arrowsize=0.05291667cm 2.0,arrowlength=1.4,arrowinset=0.4]{<-}(5.2,-1.945)(4.38,-1.545)
\usefont{T1}{ptm}{m}{n}
\rput(2.9145312,-1.245){$\Delta_{n+1}(\bu^1,\bu^2)$}
\usefont{T1}{ptm}{m}{n}
\rput(11.465625,-3.0175){$g_{n+1}(\bu^2)$}
\usefont{T1}{ptm}{m}{n}
\rput(7.004531,-2.2975){$g_{n+1}(\bu^1)$}
\psline[linewidth=0.06cm,arrowsize=0.05291667cm 2.0,arrowlength=1.4,arrowinset=0.4]{->}(7.96,0.8525)(7.98,-0.8075)
\psline[linewidth=0.024cm,linestyle=dashed,dash=0.16cm 0.16cm](12.28,-3.1825)(11.18,-4.1625)
\psline[linewidth=0.024cm,linestyle=dashed,dash=0.16cm 0.16cm](12.32,-3.1625)(13.54,-4.1225)
% \psdots[dotsize=0.12,fillstyle=solid,dotstyle=o](5.16,2.4125)
% \psdots[dotsize=0.12,fillstyle=solid,dotstyle=o](5.16,2.8125)
% \psdots[dotsize=0.12,fillstyle=solid,dotstyle=o](5.16,3.1725)
% \psdots[dotsize=0.12,fillstyle=solid,dotstyle=o](5.16,3.6125)
% %\psdots[dotsize=0.12,fillstyle=solid,dotstyle=o](5.16,4.0325)
% \psdots[dotsize=0.12,fillstyle=solid,dotstyle=o](5.16,4.4125)
% \psdots[dotsize=0.12,fillstyle=solid,dotstyle=o](5.16,4.8325)
% \psdots[dotsize=0.12,fillstyle=solid,dotstyle=o](12.38,2.4125)
% \psdots[dotsize=0.12,fillstyle=solid,dotstyle=o](12.38,2.8325)
% \psdots[dotsize=0.12,fillstyle=solid,dotstyle=o](12.38,3.2525)
% \psdots[dotsize=0.12,fillstyle=solid,dotstyle=o](12.38,3.6325)
% \psdots[dotsize=0.12,fillstyle=solid,dotstyle=o](12.38,4.0125)
% \psdots[dotsize=0.12,fillstyle=solid,dotstyle=o](12.38,4.4325)
% \psdots[dotsize=0.12,fillstyle=solid,dotstyle=o](12.38,4.8325)
\psdots[dotsize=0.12,fillstyle=solid,dotstyle=o](5.96,-1.9675)
\psdots[dotsize=0.12,fillstyle=solid,dotstyle=o](5.96,-1.5675)
\psdots[dotsize=0.12,fillstyle=solid,dotstyle=o](5.96,-1.1475)
\psdots[dotsize=0.12,fillstyle=solid,dotstyle=o](5.96,-0.7675)
\psdots[dotsize=0.12,fillstyle=solid,dotstyle=o](12.3,-2.7475)
\psdots[dotsize=0.12,fillstyle=solid,dotstyle=o](12.3,-2.3475)
\psdots[dotsize=0.12,fillstyle=solid,dotstyle=o](12.3,-1.9275)
\psdots[dotsize=0.12,fillstyle=solid,dotstyle=o](12.3,-1.5075)
\psdots[dotsize=0.12,fillstyle=solid,dotstyle=o](12.3,-1.0875)
\rput[r](12.1,-1.5075){$\bigl( g_{n+1}(\bu^2)\bigr)^{\uparrow 4}\to$}
\end{pspicture}
}
\caption{The vertices
$ g_{n+1}(\bu^1),g_{n+1}(\bu^2)$ and the history set
$\Delta_{n+1}(\bu^1,\bu^2) $ when $W_n^{\text{move}}
= \{  g_n (\bu^1), g_n (\bu^2) \} $, $W_n^{\text{stay}}
= \emptyset $.  Note the vertices above $ g_{n+1}(\bu^1)$
and $g_{n+1}(\bu^2) $ are unexplored.}
\label{Fig:TreeStep1}
%\end{center}
\end{figure}
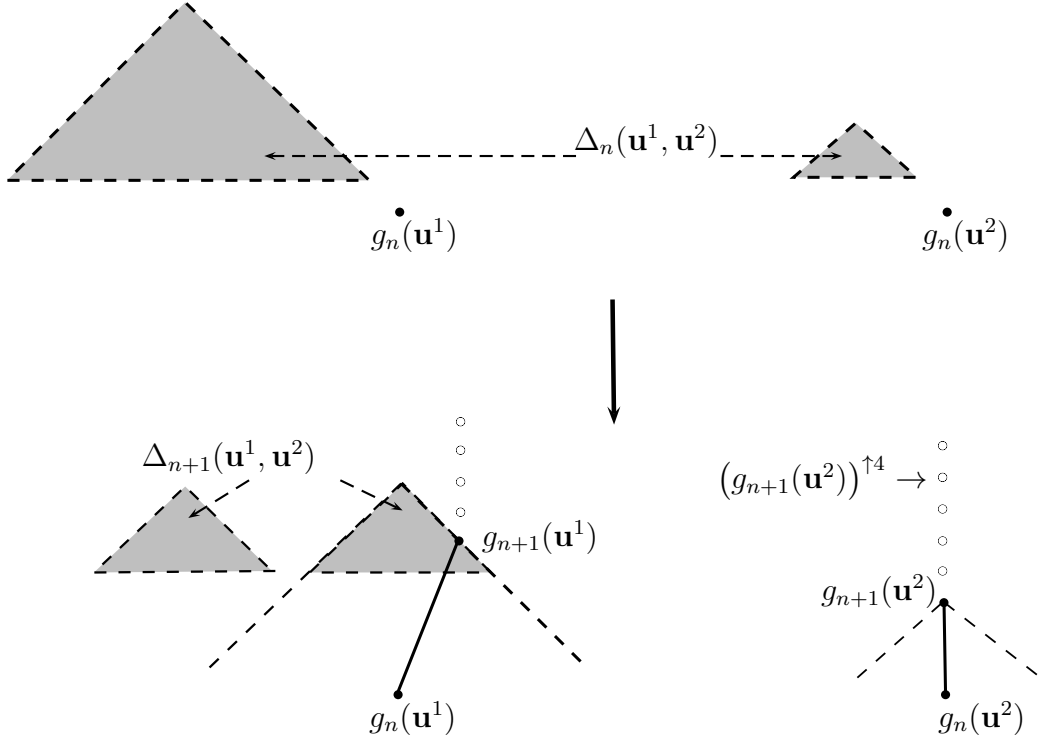

%\end{document}

%\end{document}

Let $ \Delta \subseteq \Z^d $ be a finite union of
$d$-dimensional tetrahedrons, with each tetrahedron in $\Delta$ having
a $(d-1)$-dimensional tetrahedron as a base on the hyperplane $ Q_{r}:= \{\bw
\in \Z^d: \bw(d) = r\}$ for some $ r \in \Z$. In other words, we have  $\Delta =
\cup_{i=1}^{p} S^{+} (\bw^i, t_i) $ for some $p \geq 1$
and $ \bw^{i} \in Q_{r}, t_i \geq 0$ for $ 1 \leq i \leq p$.
We denote this class of subsets by $ \mathcal{D}_r$.
% Further, for any $ \Delta \in \mathcal{D}$, let $ \text{int}(\Delta)
% = \cup_{i=1}^{p} S^{+} (\bw^i, t_i - 1) $ be the {\em interior} of
% $ \Delta $ and $ \delta^{(\text{out})}(\Delta)  = \Delta \setminus
%  \text{int}(\Delta) $ be the outer boundary of $ \Delta$.
Further, for any $ \Delta = \cup_{i=1}^{p} S^{+} (\bw^i, t_i)\in 
\mathcal{D}_r$, let 
$ \Psi : \Delta \to [0,1] $ be a mapping such
that $ \Psi (\bw) \geq  p $ for all $\bw \in \Delta^0 = \cup_{i=1}^{p} S^{+} 
(\bw^i, t_i - 1) $,
with $\Psi$ being the empty function when $\Delta = \emptyset$.
Let 
\begin{quote}
$ \mathfrak{S}^{(k)}:=  \bigl\{ \mathfrak{s} = (\bv^1, \dotsc, \bv^k,  
\Delta, \Psi)
: \bv^i \in \Z^d \mbox{  for } i = 1, \dotsc, k $, 
$\Delta \in \mathcal{D}_r$ for $r:= \min\{\bv^1(d), \dotsc, \bv^k(d)\}$, 
 $ \bv^i  \in \Delta \setminus \Delta^0 \mbox{  for all } \bv^i \text{ 
with } \bv^i(d) > r$ and $\Psi$ satisfies the conditions above$\bigl\}$.
\end{quote}

% 
% Let us denote $ \mathfrak{S}^{(k)}:=  \{ \mathfrak{s} = (\bv^1, \bv^2, \dotsc, \bv^k,  \Delta, \Psi)
% : \bv^i \in \Z^d$ for all $ i = 1, \dotsc, k, \Delta \in \mathcal{D}$  and $ \Psi : \Delta \to [0,1] $ is a mapping such
% that $ \Psi (\bw) \geq  p $ for all $\bw \in  \text{int}(\Delta)\}$.

% At the start of the process a vertex $\bu^{i}$ may be closed. When it's turn comes to
% move it explores the region above it like any other vertex.
%
% If $\max\{ \bu^{i}(d) - \bu^{j}(d): \; 1 \leq i, j \leq k\} = m_0$, then by the
% $m_0+1$ th move all the vertices would have moved from their initial positions
% $\bu^{1}, \dotsc , \bu^{k}$.
%
%For $n > m_0$,
\begin{remark}
\label{r-mopen} We have the following observations:
\begin{itemize}
 \item[(a)] From the definition of the history region $\Delta_n $,
formed at the $n$ th step, is either empty or  an element of $ \mathcal{D}$
with the bases of  tetrahedrons being contained in $Q_{r_n+1}$.
\item[(b)] Clearly, by definition of $r_n$,  $W_n^{\text{move}} \subseteq Q_{r_n} $.
\item[(c)] From the definition of $ \Delta_n$, all vertices in the set $ \Xi_n := \Z^d \setminus
\bigl( \Delta_n \cup   \mathbb{H}  (r_n)  \bigr) $ are unexplored until the $n+1$ th step, for each $ n \geq 0$.
\end{itemize}
\end{remark}

% Furthermore, for each $ n > m_0  $, all vertices in the set ~$ \Xi_n := \Z^d
% \setminus
% \Bigl( \Delta_n \cup   \mathbb{H}  (r_n)  \Bigr) $ are unexplored until the $n+1$ th step.
% Thus, each $ n > m_0  $,  $ i = 1, \dotsc, k $,  and $ m \geq 1 $, the vertex
% $ g^{\uparrow m }_n ( \bu^i) $ does never belong to $\Delta_n$ and  {\em always unexplored},
% where $ g^{\uparrow m }_n ( \bu^i) $ is defined by
% \begin{equation*}
% \label{r-mdef}
% g^{\uparrow m }_n ( \bu^i)  (j)  :=
% \begin{cases}
% g_n ( \bu^i)  (j) & \text{ for } 1 \leq j \leq d-1 ,\\
% m+g_n ( \bu^i) (d)  & \text{ for  } j = d.
% \end{cases}\eqno{\Box}
% \end{equation*}
% %For each $ n > m_0 $, the vertices $\{(g_n(\bu^{i}))_m :  m
% %\geq 1, i = 1, \dotsc, k \}$
% %are {\em never} in $\Delta_n$ and hence unexplored until the $n+1$ th step.
% %\qed

We now obtain the Markov process implicit in our construction.
For each $ n \geq 0$, set $ \mathcal{Z}^{(k)}_n := ( g_n ( \bu^{1}) ,\dotsc,  
g_n ( \bu^{k}), \Delta_n, \Psi_n )  $.
Clearly, $  \mathcal{Z}^{(k)}_n \in  \mathfrak{S}^{(k)} $.
% $  \displaystyle{(\Z^d)^k \times \{(\bw,x):
% x \in [0,1], \; \bw \in \Delta, \; \Delta \subseteq \Z^d, \; \Delta \text{ finite}\}} $.
Let $  {\bf Y} := \{ V_{\bw} : \bw \in \Z^{d},   \bw(d) > 0 \} $ be an independent
collection of i.i.d. uniform $ [0,1]$-valued random variables.  For any $ n \geq 
1 $, suppose
$ \mathcal{Z}_n^{(k)} = \mathfrak{s} \; (= ( \bv^1, \dotsc , \bv^k, \Delta, \Psi))  $ for some
$ \mathfrak{s}  \in \mathfrak{S}^{(k)}$.
 % ( \bv^1, \dotsc , \bv^k, H) \in {\cal S} $,
We define the collection of random variables $ \tilde{\bf Y} := \{   \tilde{V}_{\bw} : \bw \in \Z^{d},
\bw(d) > r \} $ for $r= \min \{\bv^1(d), \dotsc , \bv^k(d)\}$
as follows:
% for $\Delta$ the associated history region of the history $H$,
\begin{equation*}
  \tilde{V}_{\bw} :=
\begin{cases}
\Psi (\bw) & \text{ if } \bw \in \Delta; \\
V_{\bw^{\prime}} & \text{ if } \bw \not \in \Delta,  \bw(j) = \bw^{\prime} (j), 
j \neq d \text{ and } \bw(d)
=  \bw^{\prime} (d) + r .
\end{cases}
\end{equation*}
The above definition implies that $ \tilde{\bf Y} $  is a function of $ {\bf Y} $ and $ \mathfrak{s}  $,
say $ \tilde{\bf Y} = f ( {\bf Y} , \mathfrak{s} ) $
where $ f $ is a function from $ [0,1]^{ \Z^d \setminus \mathbb{H}  (0)  } \times \mathfrak{S}^{(k)}$ to
$ [0,1]^{ \Z^d \setminus \mathbb{H}  (r)  }$.
From the above definition and the fact that the vertices in $ \Xi_n = \Z^d \setminus
\bigl( \Delta \cup   \mathbb{H}  (r)   \bigr)  $ are unexplored, and hence
can be replaced by another set of i.i.d. uniform random variables,  for the family
$ {\bf X} := \{ U_{\bw}: \bw \in \Z^{d},   \bw(d) > r  \} $, we have
\begin{equation*}
{\bf X}  \mid  \mathcal{Z}_n^{(k)}
\stackrel{d}{=}   f \bigl( {\bf Y} ,  \mathcal{Z}_n^{(k)} \bigr) .
\end{equation*}

From the definition of the process, we obtain that  $ g_{n+1} ( \bu^{1}), \dotsc,
g_{n+1} ( \bu^{k}) $, $\Delta_{n+1}$ and $ \Psi_{n+1}  $ is a function of  $  
\mathcal{Z}_n^{(k)} = \bigl( g_{n}
( \bu^{1}), \dotsc, g_{n} ( \bu^{k}), \Delta_{n} , \Psi_{n} \bigr)$
and $ {\bf X} $, i.e.,
\begin{equation*}
\mathcal{Z}_{n+1}^{(k)}
= f_1 \bigl( \mathcal{Z}_n^{(k)} , {\bf X} \bigr)
\end{equation*}
where $ f_1 $ is a function on $ \mathfrak{S}^{(k)} \times [0,1]^{ \Z^d \setminus \mathbb{H}
(r_n)  } \to \mathfrak{S}^{(k)} $.
Therefore, from the above observation,
the conditional distribution of $  \mathcal{Z}_{n+1}^{(k)}  $,  given $ \{ \mathcal{Z}_j^{(k)} : 0 \leq j \leq n\} $,
is the same as that of $ f_1 \bigl( \mathcal{Z}_n^{(k)}, f ( {\bf Y} , \mathcal{Z}_n^{(k)} )  \bigr) $.
Hence, the process $ \{  \mathcal{Z}_n^{(k)} : n \geq 1 \} $ admits a random mapping representation,
which proves the Markov property (see, for example, Levin {\it et al.}\/ 
\cite{LPW08}).

% For vertices  $ \bu^{1},   \dotsc , \bu^{k}$ with $m_0:=
% \max\{\bu^{i}(d)- \bu^{j}(d):\; 1 \leq i, j \leq k\} $
\begin{prop}
\label{k-Markov}
The process $\{ \mathcal{Z}_n^{(k)} = (g_{n}(\bu^{1}),g_{n}(\bu^{2}), \dotsc , 
g_{n}(\bu^{k}), \Delta_{n},\Psi_n ) : n \geq 0\}$
is Markov with state space $\mathfrak{S}^{(k)} $.
\end{prop}
% := \displaystyle{(\Z^d)^k \times \{(\bw,x):
% x \in [0,1],
% \; \bw\in S, \; S
% \subseteq \Z^d, \; S \text{ finite}\}}$.

%\subsection{Regeneration and martingale}

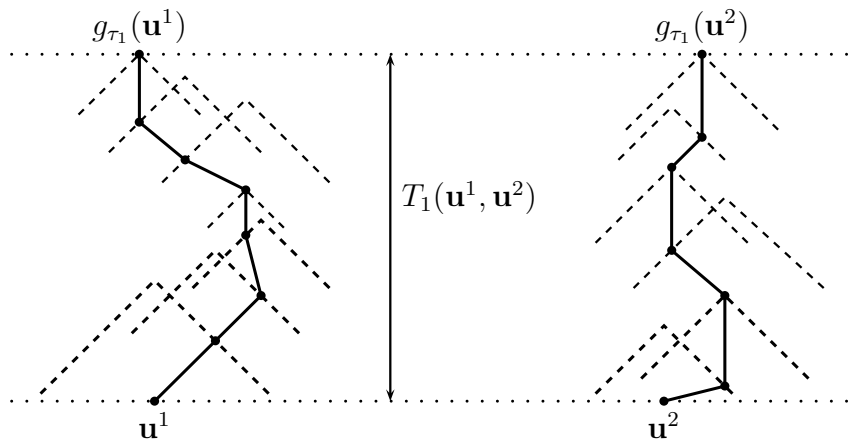
\begin{figure}[!htb]
\centering
%\begin{center}\leavevmode
%includegraphics[width=.5\textwidth]{hist1}

\begin{pspicture}(0,-2.7)(10.76,2.70)

%The point \bu_1
\rput[t](1.4,-2.65){$\bu^1$}
%Path From {\bu^1}
\psline[linewidth=0.04cm](1.4,-2.5)(2.2,-1.7)(2.8,-1.1)(2.6,-0.3)(2.6,0.3)(1.8,0.7)(1.2,1.2)(1.2,2.1)
%vertices for path from \bu_1
\psdots[dotsize=0.12](1.4,-2.5)(2.2,-1.7)(2.8,-1.1)(2.6,-0.3)(2.6,0.3)(1.8,0.7)(1.2,1.2)(1.2,2.1)
% the point g_{\tau_1}(\bu^1)
\rput[b](1.2,2.25){$g_{\tau_1}(\bu^1)$}
% History region - for path from \bu_1
\psline[linewidth=0.04cm,linestyle=dashed,dash=0.1cm 0.1cm](-0.1,-2.4)(1.4,-0.9)(2.9,-2.4)
\psline[linewidth=0.04cm,linestyle=dashed,dash=0.1cm 0.1cm](1.1,-1.6)(2.2,-0.5)(3.3,-1.6)
\psline[linewidth=0.04cm,linestyle=dashed,dash=0.1cm 0.1cm](1.9,-1.0)(2.8,-0.1)(3.7,-1.0)
\psline[linewidth=0.03cm,linestyle=dashed,dash=0.1cm 0.1cm](2.1,-0.2)(2.6,0.3)(3.1,-0.2)
\psline[linewidth=0.03cm,linestyle=dashed,dash=0.1cm 0.1cm](1.5,0.4)(2.6,1.5)(3.7,0.4)
\psline[linewidth=0.03cm,linestyle=dashed,dash=0.1cm 0.1cm](0.8,0.8)(1.8,1.8)(2.8,0.8)
\psline[linewidth=0.03cm,linestyle=dashed,dash=0.1cm 0.1cm](0.4,1.3)(1.2,2.1)(2.0,1.3)

%The point \bu_2
\rput[t](8.1,-2.65){$\bu^2$}
% From {\bu^2}
\psline[linewidth=0.04cm](8.1,-2.5)(8.9,-2.3)(8.9,-1.1)(8.2,-0.5)(8.2,0.6)(8.6,1.0)(8.6,2.1)
%vertices for path from \bu_2
\psdots[dotsize=0.12](8.1,-2.5)(8.9,-2.3)(8.9,-1.1)(8.2,-0.5)(8.2,0.6)(8.6,1.0)(8.6,2.1)
% the point g_{\tau_1}(\bu^1)
\rput[b](8.6,2.25){$g_{\tau_1}(\bu^2)$}

% History region - for path from \bu_2
\psline[linewidth=0.04cm,linestyle=dashed,dash=0.1cm 0.11cm](7.2,-2.4)(8.1,-1.5)(9.0,-2.4)
\psline[linewidth=0.04cm,linestyle=dashed,dash=0.1cm 0.1cm](7.8,-2.2)(8.9,-1.1)(10.0,-2.2)
\psline[linewidth=0.03cm,linestyle=dashed,dash=0.1cm 0.1cm](7.7,-1.0)(8.9,0.2)(10.2,-1.0)
\psline[linewidth=0.03cm,linestyle=dashed,dash=0.1cm 0.1cm](7.2,-0.4)(8.2,0.6)(9.2,-0.4)
\psline[linewidth=0.03cm,linestyle=dashed,dash=0.1cm 0.1cm](7.5,0.7)(8.2,1.4)(8.9,0.7)
\psline[linewidth=0.03cm,linestyle=dashed,dash=0.1cm 0.1cm](7.6,1.1)(8.6,2.1)(9.6,1.1)

% bottom dashed line -- starting line of vertices
\psline[linewidth=0.04cm,linestyle=dotted,dotsep=0.16cm](-0.5,-2.5)(10.7,-2.5)
% Top line - regeneration line
\psline[linewidth=0.04cm,linestyle=dotted,dotsep=0.16cm](-0.5,2.1)(10.7,2.1)
% The length T_1 (\bu^1, \bu^2)
\psline[linewidth=0.03cm]{<->}(4.5,-2.5)(4.5,2.1)
\rput[l](4.65,0.2){$T_1(\bu^1,\bu^2)$}

\end{pspicture}
%\end{pspicture}
%\input{construction.tex}
\caption{At regeneration step $\tau_1 (\bu^1,\bu^2)$ of the process
$ g_{\tau_1}(\bu^1) (d) = g_{\tau_1}(\bu^2)(d) $ and $ \Delta_{\tau_1} = \emptyset $}
\label{Fig:Regeneration1}
%\end{center}
\end{figure}

For the remainder of this section we fix $ \bu^{1},  \dotsc , \bu^{k}$ with
$\bu^{1}(d)=   \cdots = \bu^{k}(d) $.
Set $  \tau_0 =  \tau(\bu^{1},\dotsc, \bu^{k}) := 0 $ and, for $ l \geq 1
$, define
\begin{equation}
\label{r-tau}
\begin{split}
\tau_l  = \tau_l(\bu^{1},\dotsc, \bu^{k}) & := \inf \{ n > \tau_{l-1} : \Delta_n = \emptyset\}; \\
\sigma_l =  \sigma_l (\bu^{1},\dotsc, \bu^{k}) & := \tau_l(\bu^{1},\dotsc, \bu^{k})  - \tau_{l-1}(\bu^{1},\dotsc,
\bu^{k}).
\end{split}
\end{equation}
We call $ \tau_l$ the step at which the $ l$ th {\em simultaneous regeneration of $ k $ joint paths} occurs.
We note here that $\tau_l$ denotes the number of steps (in the above construction) required
for the joint process to regenerate   (i.e., to reach a state of empty history for the $l$ th time) and
$ \sigma_l$ denotes the total number of steps  (again in the above construction) between the $ l-1$ th and $l$ th
simultaneous regeneration of $ k $ joint paths.
This is not the same as the time
(measured as the distance in the $d$ th co-ordinate) for regeneration,
which we later denote by $ T_l$ (see Figure
\ref{Fig:Regeneration1}).
Also at each regeneration step $\tau_l$,  the paths must be at the same level in
terms of their $d$ th co-ordinate, i.e.,
$g_{{\tau}_l}(\bu^{1})(d) = \cdots = g_{{\tau}_l}(\bu^{k})(d)$.
%(see Figure \ref{Fig:Regeneration1}).

% \documentclass[a4paper,12pt]{article}
% \usepackage{pst-all}
% \begin{document}

%\end{document}

%\end{document}

Our first task is to show that the Markov process, defined in Proposition
\ref{k-Markov}, regenerates almost surely.
In fact, we prove the much stronger statement that the number of inter-regeneration steps
has exponentially decaying tail probabilities.
%More precisely,
\begin{prop}
\label{prop:RecTimeExpTail}
For any $ l \geq 1 $ and $ \bu^{1},\dotsc, \bu^{k} \in \Z^d $ with $ \bu^{1} (d) =
\cdots = \bu^{k} (d) $, we have
\begin{equation}
\label{eqn:expDecay}
 \P ( \sigma_l \geq n ) \leq  C_1^{(k)} \exp ( - C_2^{(k)} n )
\end{equation}
for all $ n \geq 1 $, where $ C_1^{(k)} $ and $ C_2^{(k)} $ are positive constants, depending on $k$, but not
 on $l, n $ or $ \bu^{1},\dotsc, \bu^{k} $.
\end{prop}

Since $ \tau_l < \infty $ almost surely, we obtain that
\begin{cor}
\label{cor:EmptyHistMC}
For $ \bu^1, \dotsc, \bu^k $ as above, the process
\begin{equation*}
\bigl\{ \bigl(g_{\tau_l}(\bu^1), \dotsc , g_{\tau_l}(\bu^k) \bigr) :  l \geq 0  \bigr\} 
\text{ is a Markov chain on } (\Z^d)^k.
\end{equation*}
\end{cor}

For $ \bw \in \Z^d$, let $ ( \bw)^{\uparrow m } $ be defined by
\begin{equation*}
\label{r-mdef}
( \bw)^{\uparrow m }  (j)  :=
\begin{cases}
\bw (j) & \text{ for } 1 \leq j \leq d-1 ,\\
m+ \bw (d)  & \text{ for  } j = d.
\end{cases}
\end{equation*}
%The set $ \{ ( \bw)^{\uparrow m }  : m \geq 1 \} $ represents the line above 
the vertex $\bw$.
Note
\begin{itemize}
\item[] if $ \bw \in \Z^d $,$ \bv \in V$ are such that
$ \bv (d) > \bw(d) $, then $   ||   h(\bw)  - \bw  ||_1 \leq  ||  \bv  - \bw 
||_1 $.
\end{itemize}
The main idea behind regeneration is contained in the above 
simple
observation and the following Proposition \ref{prop:StrucDelta}. 
\begin{prop}
\label{prop:StrucDelta} For $ \bu^1, \dotsc, \bu^k $ as in the construction of the process with the history region $\Delta_n$ and the quantity $r_n$ as described there, we have,
for any $n \geq 0,  m \geq 1$ and $ 1 \leq i \leq k$, the vertex
$ \bigl( g_n(\bu^i) \bigr)^{\uparrow m } \not \in \Delta_n \cup \mathbb{H}(r_n) $.
\end{prop}

\noindent{\bf Proof }: Fix $ n \geq 0,  m  \geq 1$ and $ 1 \leq i \leq k$.
We have $  \bigl( g_n(\bu^i) \bigr)^{\uparrow m } (d) = g_n(\bu^i)(d) + m > g_n(\bu^i)(d)
\geq r_n $, so that $ \bigl( g_n(\bu^i) \bigr)^{\uparrow m } \not \in  \mathbb{H}  (r_n)$.
It is enough to show $ \bigl( g_n(\bu^i) \bigr)^{\uparrow m }  \not \in \cup_{j=0}^n \Delta_{j} $.
First $  \bigl( g_n(\bu^i) \bigr)^{\uparrow m }  \not \in \Delta_{0} $. So take 
$ n \geq 1$ and inductively assume that $ \bigl( g_n(\bu^i) \bigr)^{\uparrow m } 
\not\in  \Delta_{j}$ for some $0 \leq j \leq n-1$.
% 
% So take $ n \geq 1$ and
% fix any $ j = 0, 1, \dotsc, n-1 $
%  and  
%  
Let $ \bw \in W_j^{\text{move}}$.
Since the  vertices with the smallest $d$ th co-ordinate  are in  $W_j^{\text{move}}$ and, in the next step their $d$ th co-ordinate increases by at least one unit, we
have $ \bw(d) = r_j < r_n \leq g_n (\bu^i) (d) $. Since
$ g_n ( \bu^i) \in V$,  clearly $ ||   h(\bw)  - \bw  ||_1 \leq  ||  g_n(\bu^i)  - \bw ||_1 $. Further,
we have $ || \bigl( g_n(\bu^i) \bigr)^{\uparrow m } - \bw  ||_1 = m +  ||  g_n(\bu^i)  - \bw ||_1
> ||  g_n(\bu^i)  - \bw ||_1 $ and
hence $ \bigl( g_n(\bu^i) \bigr)^{\uparrow m }  \not\in S^{+} ( \bw ,
||  h(\bw)  - \bw  ||_1) $. Thus, $ \bigl( g_n(\bu^i) \bigr)^{\uparrow m } \not\in  \Delta_{j+1}$. \qed

In order to prove Proposition \ref{prop:RecTimeExpTail}, we define a random variable
$ L_n$ which represents the {\em height} of the history region
$\Delta_n$, measured along the $d$ th co-ordinate from the lowest vertex among
$ g_n (\bu^{1}) , \dotsc, g_n (\bu^{k}) $. Using Proposition \ref{prop:StrucDelta}, for any $ 1 \leq i \leq k$,
the set of vertices $ \{ \bigl( g_n(\bu^i) \bigr)^{\uparrow m } : m \geq 1\} $ 
is not
explored till the $n$ th step (see Figure \ref{Fig:TreeStep1}). This provides an upper bound on the size of
the next step, and hence on the increment of the height of the history region.
The height of the first vertex of $ \{ \bigl( g_n(\bu^i) \bigr)^{\uparrow m } : 
m \geq 1\} $ which is open from $g_n(\bu^i)$ is a geometric random variable 
irrespective of the history  carried.
Using these geometric random variables, we construct a coupling with a Markov chain $ M_n $
which dominates the height random variable. Hence, the Markov chain's return time to $ 0 $
will dominate the return time of $ L_n$ to $ 0 $. The Markov chain is 
constructed
so that it uses an independent sequence of random variables when $ L_n $ has already
returned to $ 0 $ but $ M_n $ is positive.

\begin{remark}
\label{rem:L2Statement}
Note here that  Proposition \ref{prop:StrucDelta} remains valid even if we consider
the  $L_2$ norm and define regeneration in the same way, i.e., when the history
set becomes empty. Though the geometric structure of the history region changes,
 we can still provide a bound on the increment of the height  and construct a 
dominating Markov chain with the same properties as above.  
The geometry of the $L_1$
structure has been used very minimally, and wherever they have been used we may see that the results also hold 
when it is replaced with the $L_2$ norm.
%However, in the Poisson directed spanning forest, apart from the history being never
%empty (which may be handled by restricting ), we are unable to provide
\end{remark}

To prove Proposition \ref{prop:RecTimeExpTail}, we need an auxiliary lemma
on Markov chains, whose proof is given in the Appendix.
Let $ \{ \theta_n : n \geq 1 \} $ be a sequence of i.i.d.
positive integer valued random variables with $ \P ( \theta_1 = 1 ) > 0 $ and $
\P (
\theta_ 1 \geq n ) \leq C_3 \exp ( -C_4 n ) $ for all $ n \geq 1 $ where $ C_3, C_4 $
are positive constants.
Define a sequence of random variables as follows: $ M_0 := 0 $ and for $ l \geq 0, M_{l+1} :=
\max \{ M_l , \theta_{l+1} \} - 1 $. Let
$\tau^M := \inf \{l \geq 1: M_l = 0 \}$ be the first return time  of $ M_l$ to $ 0 $.
\begin{lemma}
\label{lem:MCAuxRes}
For $ n \geq 1$, we have
\begin{equation*}
\P ( \tau^M \geq n ) \leq  C_5 \exp ( - C_6 n )
\end{equation*}
where $ C_5 $ and $ C_6 $ are positive constants.
\end{lemma}

\noindent {\bf Proof of Proposition \ref{prop:RecTimeExpTail}: }
We first observe that by the Markov property (Proposition \ref{k-Markov}) it is enough to
show the result for $ l = 1 $. In order to study that, we define,
\begin{equation}
\label{eqn:def:MaxLevel}
L_n :=
\begin{cases}
\max \{ \bw (d) : \bw \in \Delta_n\} - r_n   & \text{ if } \Delta_n
\neq \emptyset \\
0 & \text{ if }  \Delta_n = \emptyset
\end{cases}
\end{equation}
where $ r_n = \min\{ g_n (\bu^{i}) (d) : i = 1, \dotsc, k\} $.
We set,
\begin{equation*}
\label{r-tau_L}
\tau^L  := \inf \{ n \geq 1 : L_n = 0 \}
\end{equation*}
and observe that $ \tau_1 = \tau^L. $

% We construct random variables $M_n$ and $X_n$ such that $M_n$ is constructed
% as in Lemma \ref{lem:MCAuxRes} using the property of the history
% region $ \Delta_{n} $ that
For any fixed $ n \geq 0$, using Proposition \ref{prop:StrucDelta},
we have that $  \bigl(g_n (\bu^i) \bigr)^{\uparrow m }  $ is unexplored for $ m
\geq 1 $ and $1 \leq i \leq k $ (see Figure \ref{Fig:TreeStep1}). We now define the collection of random variables
\begin{equation}
\label{r-geome}
\Bigl\{J_{n+1}(\bw) := \inf\{m \geq 1:   \bigl( \bw \bigr)^{\uparrow m }   \in V\}:
\bw \in W_{n}^{\text{move}} \Bigr\},
\end{equation}
where $ V $ is the set of all open points. This is a collection of i.i.d. geometric random
variables with parameter $p$, i.e. each of the random variables takes the  value $m$ with
probability $p (1-p)^{m-1}$ for $m = 1,2,\dotsc$. Also,
\begin{equation}
\label{r-heightgeom}
||g_{n}(\bw) - g_{n+1}(\bw)||_1 \leq J_{n+1}(\bw) \text{ for all $\bw$ with }
 g_{n}(\bw) \in W_{n}^{\text{move}}.
\end{equation}
Let $\{G^{i,1}_{n}:\;  1\leq i\leq k, n\geq 0\}$ be another family of i.i.d.
geometric random variables with parameter $p$,
independent of $\{U_{\bw}: \bw \in \Z^d\}$.

Now given $g_{n}(\bu^{1}),\dotsc, g_{n}(\bu^{k})$ and $H_{n}$,  we
define $ \{ M_n := M_n(\bu^{1}, \dotsc,\bu^{k}), X_n := X_n(\bu^{1}, \dotsc,\bu^{k}) : n \geq 0 \} $
as follows:
$$
\text{set } M_0 = 0 = X_0 \text{ and } M_{n+1} = \max \{ M_n, J^{1}_{n+1} \} - 1
\text{ for } n \geq 0
$$
where
\begin{align}
\label{eqn:DefJn}
J^{1}_{n+1} :=
\begin{cases}
 \max\{J_{n+1}(\bu): g_{n}(\bu) \in W_{n}^{\text{move}} \}   & \text{ if }
\# W_{n}^{\text{move}} = k \text{ and } X_n = 0,\\
 \max\{G^{i,1}_{{n+1}}, J_{n+1}(\bu): & \\
\quad g_{n}(\bu) \in W_{n}^{\text{move}},
i = 1, \dotsc, k-k' \} & \text{ if } \# W_{n}^{\text{move}} = k'<k\text{ and } X_n = 0,\\
\max\{G^{i,1}_{n+1}: 1\leq i\leq k\} & \text{ if } X_n = 1,
% g_{n}(\bu) = g_{n}(\bv),\\
% \max\{J_{n}(\bv),G_n\} \text{ if } g_{n}(\bu)(d) >
%g_{n}(\bv)(d).
\end{cases}
\end{align}
and
\begin{align}
% \label{eqn:DefJn}
X_{n+1} :=
\begin{cases}
1   & \text{ if }
X_n = 0, L_{n+1} = 0\\
X_n   & \text{ otherwise. }
% g_{n}(\bu) = g_{n}(\bv),\\
% \max\{J_{n}(\bv),G_n\} \text{ if } g_{n}(\bu)(d) >
%g_{n}(\bv)(d).
\end{cases}
\end{align}
From (\ref{eqn:DefJn}) it follows that $\{ J^{1}_{n+1} : n\geq 0\}$ is a family of
i.i.d. copies of $J$ where for any $m \geq 1$,
\begin{align}
\label{eqn:DistJn}
\P(J \geq m )  = 1 -  (1- (1-p)^{m-1})^k \leq k (1-p)^{m-1}
\end{align}
and hence the sequence
$ \{ J^{1}_n : n \geq 1 \} $ satisfies the conditions of Lemma \ref{lem:MCAuxRes}.

Further, we
claim that $ 0 \leq L_n \leq M_n$  for all $ 0 \leq n \leq \tau^L $. Indeed, this
holds for $ n = 0 $, and assume
that it holds for some $ 0\leq n < \tau^L$. If $\Delta_{n+1} = \emptyset$ then we have $ 0 = L_{n+1} \leq M_{n+1}$.
Otherwise if $ \bw \in \Delta_{n+1} $, then, from the
definition of $ \Delta_{n+1} $,
either $ \bw \in \Delta_n $
or $ \bw \in S^{+} ( \bu, || \bu - h(\bu) ||_1 ) $ for some $ \bu \in
W_n^{\text{move}}$.
Therefore, from (\ref{r-heightgeom}) and  (\ref{eqn:DefJn}),
$ \bw (d) \leq \max \{ \max \{ \bu (d) : \bu \in \Delta_n \} , \min \{
g_n (\bu^{i})(d) , \; 1 \leq i \leq k \} +
|| \bu - h(\bu) ||_1 : \bu \in W_n^{\text{move}} \} \leq \max \{ L_n +
r_n , r_n + J_{n+1} \} =
\max \{ L_n , J_{n+1} \} + r_n $.
Also $ r_{n+1} = \min \{ g_{n+1} (\bu^{i})(d)  , \; 1 \leq i \leq k  \}
\geq \min \{ g_n (\bu^{i})(d)  , \; 1 \leq i \leq k  \} + 1 = r_n +
1$. Thus $ L_{n+1} \leq \max \{
L_n ,  J_{n+1} \} - 1 \leq  \max \{ M_n ,  J_{n+1} \} - 1 = M_{n+1} $.

Define,
\begin{equation*}
\tau^M = \tau^M(\bu^{1},\dotsc, \bu^{k}) := \inf \{n \geq 1: M_n = 0 \} .
\end{equation*}
Note that the distribution of $\tau^M(\bu^{1},\dotsc, \bu^{k})$ does not depend on
$\bu^{1},\dotsc, \bu^{k}$.
From the above observation that $ 0 \leq L_n \leq M_n$ for $0\leq n \leq \tau_1$,
we obtain that
\begin{equation*}
\tau_1 = \tau^L \leq \tau^M .
\end{equation*}
Using Lemma \ref{lem:MCAuxRes}, we obtain Proposition
\ref{prop:RecTimeExpTail}. \qed

\medskip

The following lemma will be used to show that the inter-regeneration times as 
well as the width of the explored regions
during a regeneration have exponentially decaying tail probabilities. Let $ \{ \theta_i : i \geq 1 \} $ be i.i.d.
random variables and $ N $ be any random variable taking values in $ \{ 0, 1, 2, \dotsc \} $.
We define the random sum $ S $ as follows:
\begin{equation*}
S := \begin{cases} 0  & \text{  if } N = 0, \\
\sum_{i=1}^{n} \theta_i & \text{ if } N = n.
\end{cases}
\end{equation*}
Then, we have following lemma.
\begin{lemma}
\label{lem:RandomSumHavingExpTail}
Suppose that for some $ \beta > 0 $ and $ \alpha > 0 $, $ \E \bigl( \exp( \beta \theta_1) \bigr) < \infty $
and $ \E \bigl( \exp( \alpha N) \bigr) < \infty $.  Then, there exists $ \gamma > 0 $ such that
$ \E \bigl( \exp( \gamma S) \bigr) < \infty $.
\end{lemma}
We note here that no assumption of independence or the structure of dependence between
$ N $ and $ \theta_i$'s have been imposed here.
The proof uses the Cauchy-Schwartz inequality and has been relegated to the Appendix.

We now consider the width of the explored region between the
$l-1$ and $ l $ th regenerations. For the process starting from
$ \bu^1,  \dotsc , \bu^k$ with
$\bu^1(d)=   \cdots = \bu^k(d) $ define
\begin{equation}
 \label{r-W2}
W_l = W_l (\bu^1, \dotsc,  \bu^k) := \sum_{n= \tau_{l-1}  + 1}^{\tau_{l}
}\sum_{i=1}^k
||g_n(\bu^i)-g_{n-1}(\bu^i)||_1.
\end{equation}

Using $\{G^{i,l+1}_{n}:\;  1\leq i\leq k, n\geq 0\}$,
another family of i.i.d. geometric random variables with parameter $p$ and
independent of both $\{U_{\bw}: \bw \in \Z^d\}$ and $\{G^{i,j}_{n}:\;  1\leq 
i\leq k, 1 \leq j \leq l, n\geq 0\}$,
our construction ensures that
% construct $\{M_n(g_{\tau_l}(\bu^1), \dotsc ,
% g_{\tau_l}(\bu^k)),X_n(g_{\tau_l}(\bu^1), \dotsc , $ $
% g_{\tau_l}(\bu^k)):n\geq 0\}$ such that
$\sigma_{l+1}(\bu^1, \dotsc , \bu^k)
\leq \tau^M(g_{\tau_l}(\bu^1), \dotsc ,
g_{\tau_l}(\bu^k))$ and $\tau^M(g_{\tau_l}(\bu^1), \dotsc ,
g_{\tau_l}(\bu^k))$ is an i.i.d. copy of $\tau^M(\bu_1,\dotsc,\bu_k)$.
Also, for $\tau_l \leq n < \tau_{l+1}$, we have
\begin{equation*}
 \sum_{i=1}^k ||g_{n+1}(\bu^i)-g_{n}(\bu^i)||_1 \leq
\sum_{g_{n}(\bu^i) \in W^n_{\text{move}}} J_{n+1}(\bu^i) \leq k
J^{l+1}_{(n-\tau_l)+1}
\end{equation*}
where the last sum is over distinct elements of $W_n^{\text{move}}$ to avoid
double counting and $J^{l+1}_{i}$ is defined as in (\ref{eqn:DefJn}) using
$\{G^{i,l+1}_n:1\leq i\leq k, n\geq 0\}$ instead of $\{G^{i,1}_n:1\leq i\leq k, n\geq 0\}$.
Further it follows that
$W_{l+1}\leq \sum_{i=1}^{\tau^M(g_{\tau_l}(\bu_1),\dotsc,g_{\tau_l}(\bu_k))} kJ^{l+1}_{i}$
and $\sum_{i=1}^{\tau^M(g_{\tau_l}(\bu_1),\dotsc,g_{\tau_l}(\bu_k))} kJ^{l+1}_{i}$ is an
i.i.d. copy of $W^M := \sum_{i=1}^{\tau^M(\bu_1,\dotsc,\bu_k)} kJ^{1}_{i}$.
From Lemma \ref{lem:RandomSumHavingExpTail}, we conclude that for
some $ \alpha > 0 $, $ \E \bigl( \exp (\alpha W^M ) \bigr) < \infty$.

The time for the $l$ th regeneration, measured by the distance travelled by
process in the $d $ th co-ordinate (see Figure \ref{Fig:Regeneration1}),  is defined as
\begin{equation}
\label{def:TimekthReg}
T_l = T_l  ( \bu^1, \dotsc,  \bu^k ) := g_{\tau_l } ( \bu^1) (d) - \bu^1 (d) =
 g_{\tau_l } ( \bu^i) (d) - \bu^i (d) \text{ for } 1 \leq i \leq k .
\end{equation}
% For all $1 \leq i \neq j \leq d$, by our choice of the starting vertices,
% $\bu^i(d) = \bu^j(d)$; also at any regeneration step $g_{\tau_l } ( \bu^i)(d) =
% (d) g_{\tau_l } ( \bu^j) (d)$, so $T_l = g_{\tau_l } ( \bu^i) (d) - \bu^i (d)$.
%Clearly $T_l - T_{l-1} \leq W_l$ and we have
% Now, at the $n$ th stage only the vertices in $W^n_{\text{move}}$
% move, so from (\ref{r-heightgeom}) and with $J_{n+1}$ as in (\ref{eqn:DefJn}) we
% have
Clearly $T_l - T_{l-1} \leq W_l$ and from the fact that, for
some $ \alpha > 0 $, $ \E \bigl( \exp (\alpha W^M ) \bigr) < \infty$, we conclude
the following Proposition.
\begin{prop}
\label{prop:WidthTail}
 For any $ l \geq 1 $ and $ \bu^1,  \dotsc , \bu^k$ with
$\bu^1(d)=   \cdots = \bu^k(d) $, we have
\begin{equation}
\label{eqn:momentFinite}
 \P  ( T_l - T_{l-1}\geq n) \leq  \P  ( W_l  \geq n)  \leq \P (W^M \geq n) \leq C_7^{(k)} \exp ( - C_8^{(k)} n )
\end{equation}
for all $ n \geq 1 $, where $ C_7^{(k)} $ and $ C_8^{(k)} $ are positive constants, depending on $k$ but not
$ l  $ or $ \bu^1,  \dotsc , \bu^k$.
\end{prop}

Now we specialize to the case $k=1$, i.e., the process starting from just one vertex $\bu$ and consider
the inter-regeneration step process $ \sigma_l (\bu) $, the width process $ W_l (\bu) $
and the inter-regeneration time process $  T_l (\bu) - T_{l-1} (\bu)$, defined in (\ref{r-tau}),
(\ref{r-W2}) and (\ref{def:TimekthReg}) respectively. Using the Markov property 
of $ \{ \mathcal{Z}_l^{(1)}: l \geq 0 \} $
and  the translation invariance of the model, we have that $ \{  \bigl( \sigma_l (\bu), W_l (\bu) ,
(T_l (\bu)  - T_{l-1} (\bu) ) \bigr) : l \geq 1\}$ is an i.i.d.
family of random vectors taking values in $ \{ 1, 2, 3, \dotsc \}^3 $. Using the translation invariance of
our model,  the distribution of these random variables  does not depend on the 
choice of the starting vertex $\bu$.
Furthermore, each of the marginals of this random vector has exponentially decaying  tail probability.
Let
\begin{equation}
\label{not:bar}
 \overline{\bw} := (\bw(1), \dotsc , \bw({d-1}) ) \text { for } \bw = (\bw(1), \dotsc , \bw(d))
 \end{equation}
denote the first $d-1$ co-ordinates of $\bw \in \Z^d$. We observe the following:
\begin{remark}
\label{rem:FromSinglePoint}
Define, for any $ l \geq 1$,
\begin{equation}
\label{Ydef}
Y_{l}^{(\bu)} :=  \overline{g}_{\tau_{l} (\bu)} (\bu)  - \overline{g}_{\tau_{l-1} (\bu) } (\bu) .
\end{equation}
We have
\begin{enumerate}
\item[(a)] $ \{  Y_{l}^{(\bu)} : l \geq 1\} $ is a sequence of i.i.d. $ 
\Z^{d-1} $-valued random vectors,
whose distribution does not depend on $ \bu$. For $\bu =\mathbf{0}$, we denote
 $Y_{l}^{(\mathbf{0})}$ by $Y_{l}$.

\item[(b)] Denoting the $L_1$ norm in $ (d-1)$ dimensions by $ || \cdot ||_{1, d-1} $, we also observe that $  || Y_1  ||_{1, d-1} \leq
W_1 ({\mathbf 0}) $, so that we also have
\begin{equation}
 \label{eqn:EstYInc}
\P(  || Y_1  ||_{1, d-1}  \geq n ) \leq C_7^{(1)} \exp ( - C_8^{(1)} n )
\end{equation}
for all $ n  \geq 1 $, where $ C_7^{(1)}$ and $ C_8^{(1)} $ are as in  (\ref{eqn:momentFinite}).

\item[(c)] By reflection symmetry of the model, about any of the first $ (d-1)$
co-ordinates, we have  that  each co-ordinate of $Y_1 $ is marginally symmetric.
Further, the rotational symmetry of the model in the first $ (d-1)$ co-ordinates implies that  the
marginal distributions are same for $ i = 1,2, \dotsc, d-1$. In other words,
$ \P ( Y_1 (1) = +m   )  =  \P ( Y_1 (i) = +m   ) =   \P ( Y_1 (i) = -m   )  $ for all $ m \geq 1 $ and
$ 1 \leq i \leq d-1$ where $ Y_1 (i) $ denotes  the $i$ th  co-ordinate of $ Y_1$.
%Hence, for
%any $ 1 \leq i \leq d-1$, $  \E \bigl[ (Y_1 (i) )^2 \bigr]= \sigma^2 $ for some $ \sigma^2 > 0$.

\item[(d)] For any $ i, j \in \{ 1,2, \dotsc, d-1\}, i \neq j $, let $ s =  \E  \bigl(  Y_1 (i)  Y_1 (j)  \bigr)  $.
By reflection symmetry along the $i$ th co-ordinate, with other co-ordinates being fixed, we observe
that the joint distribution of $(Y_1 (i),   Y_1 (j) )$ remains unchanged. This implies that
$ s =  \E  \bigl(  (- Y_1 (i) ) Y_1 (j)  \bigr)  $ and hence, $ \E  \bigl(  Y_1 (i)  Y_1 (j)  \bigr) = 0 $.  The same
argument holds to prove $  \E  \bigl[  \bigl( Y_1 (i) \bigr)^{m_1} \bigl(  Y_1 (j)  \bigr)^{m_2} \bigr] = 0  $
for $ m_1, m_2 \geq 1$  with at least one of them being odd. \qed

\end{enumerate}

%  is a sequence of i.i.d. ,
%
% It should be noted, from the
% rotational symmetry in the first
% $ (d-1)$ co-ordinates and
% that the distribution of the co-ordinates of $Y_{l+1}$ are uncorrelated and
% marginally each of them is symmetric about $0$.   Further, for all
%  and
%
% \begin{equation}
% \label{r-Y_i}
% \text{ and } = 0
% \end{equation}
% for some $ \sigma^2 > 0 $  and  if

\end{remark}

%\input{2ASRR17.07.tex}
% !TEX root = ./DSFrf.tex
% !TEX program = XeLaTeX
\section{Martingale and independent processes}
In this section, we study the joint evolution of paths starting from 
two vertices $ \bu $ and $ \bv $ with
$\bu(d)=    \bv(d) $. The process $ \{(g_{\tau_l}(\bu),
g_{\tau_l}(\bv)): l \geq 0\} $ is a  Markov chain on $(\Z^d)^2$ by Corollary \ref{cor:EmptyHistMC}.
Further, our model is translation invariant. Hence,  the process $ \{ g_{\tau_l}(\bu) -
g_{\tau_l}(\bv): l \geq 0\}$ is also a Markov chain on $ \Z^d $.
However, as  observed earlier, (see Figure \ref{Fig:Regeneration1}), $ g_{\tau_l}(\bu) (d) =
g_{\tau_l}(\bv) (d) $ for every $ l \geq 1 $. Thus, using notation as in (\ref{not:bar}),
\begin{equation}
\label{def:Z_l}
\{Z_l = Z_l ( \bu, \bv) := \overline{g}_{\tau_l}(\bu)  - \overline{g}_{\tau_l}(\bv): l \geq 0\}
\end{equation}
is a  $\Z^{d-1}$ valued Markov chain. Further, we observe that $ \overline{{\bf 0}} = (\underbrace{0,0, \dotsc, 0}_{d-1}) $
is the only absorbing state of this Markov chain.

In Subsection \ref{subSec:Martingale}, we prove that the process $Z_l$ is a martingale for $d=2$.
Later, in Subsection \ref{subSec:IndPr} and Subsection \ref{Coupling},
we study,  for any $ d \geq 2$,  how $ Z_l $ can be approximated by a process having i.i.d.
increments when the starting point of paths are sufficiently far apart.

\subsection{Martingale}
\label{subSec:Martingale}
In this subsection we restrict ourselves to $d=2$ and fix two vertices $\bu, \bv \in \Z^2$ with $\bu(2) =  \bv(2)$.
We first observe that, for $ l \geq 0$, the regeneration time  $T_l = T_l(\bu, \bv)$ is a stopping time
with respect to the filtration $\{{\cal F}_t : t \geq  0\}$ where ${\cal F}_{t} :=
\sigma\{U_{\bw}: \bw(2) \leq \bu(2) + t\}$.
By our construction, $ g_{\tau_{l}}(\bu)$ is ${\cal F}_{T_{l}}$
measurable. Therefore, $ \overline{g}_{\tau_{l}}(\bu)$, given by the projection from $ \Z^2
\to \Z $, is also ${\cal F}_{T_{l}}$ measurable.
\begin{prop}
\label{prop:RWMartingale}
For $d = 2$ and $\bu,\bv \in \Z^2$ with $\bu(2)=\bv (2)$, the process
$ \bigl\{ \overline{ g}_{\tau_{l}}(\bu) : l \geq 0\bigr\}$ is a martingale with respect to the filtration
$\{{\cal F}_{T_{l}} : l \geq 0 \}$.
\end{prop}
The above proposition also holds  for $ \overline{g}_{\tau_{l}}(\bv) $, so we obtain
\begin{cor}
\label{cor:ZMartingaleD=2}
For $d = 2 $ and any $ \bu, \bv \in \Z^2 $ with $ \bu (2) = \bv (2) $,
the process $ \{ Z_l= Z_l ( \bu, \bv) : l \geq 0  \} $ is a martingale with respect to
the filtration $\{{\cal F}_{T_{l}} : l \geq 0 \}$.
\end{cor}
\noindent {\bf Proof of Proposition \ref{prop:RWMartingale}}:
Consider the process $ (g_n (\bu), g_n (\bv), \Delta_n ( \bu, \bv) , \Psi_n)$ starting
from $ \bu, \bv $ with $\bu(2)=\bv (2)$,  and the process $  (g_n (\bu),  \Delta_n ( \bu) , \Psi_n) $
starting from $ \bu $ with the same set of uniform random variables
$\{U_{\bw}: \bw \in \Z^2\}$.
Observe that every joint regeneration of the paths from a pair of vertices $ \bu, \bv  $ is also a regeneration of
the single path from $ \bu $, i.e., for every $ l \geq 0$, we have
\begin{equation}
\label{defn:TNL}
T_l ( \bu, \bv ) = T_{N_l}( \bu)
\end{equation}
for some sequence $ N_l= N_l(\bu, \bv ) $ (see Figure \ref{Fig:Regeneration1} where $ N_1 = 2$ for
$ \bu^1$ and $  N_1 = 3$ for $ \bu^2$). Therefore, we have,
\begin{equation*}
  \overline{g}_{\tau_{l}( \bu, \bv )}(\bu) =  \overline{g}_{\tau_{N_l}( \bu)}(\bu)
  = \sum_{ i = 1}^{N_l} Y_i^{({\bu})}
\end{equation*}
where $  Y_{i}^{({\bu})} :=  \overline{g}_{\tau_{i}( {\bu})} ( {\bu}) -
\overline{g}_{\tau_{i-1}( {\bu})} ({\bu}) $ is as in Remark \ref{rem:FromSinglePoint}.
Since $ N_l \leq T_l( \bu, \bv ) $,  and each of
$ T_i(\bu, \bv ) - T_{i-1}(\bu, \bv) $ and  $ Y_i^{({\bu})} $ has an exponentially decaying tail probability
(see Proposition \ref{prop:WidthTail} and equation (\ref{eqn:EstYInc})),  for every $ l \geq 0$, we have that
$ \E ( |  \overline{g}_{\tau_{l} (\bu, \bv )}(\bu)   | ) < \infty $.

Further we need to show that
\begin{equation}
\label{k-condexp}
\E \bigl[ \overline{g}_{\tau_{l+1}(\bu, \bv )}(\bu)  - \overline{g}_{\tau_{l}(\bu, \bv )}(\bu)   |
 {\cal F}_{T_l  (\bu, \bv)} \bigr]
= \E \bigl[  \sum_{ i = N_l +1}^{ N_{l+1} } Y^{({\bu})}_i  | {\cal F}_{T_{N_l} ( {\bf u}) } \bigr]
= 0 \text{ a.s.}
\end{equation}
Denoting $ {\cal G}_i :=  {\cal F}_{T_{i}(\bu)} $, we have that $ Y^{({\bu})}_{i+1} $
is independent of $ {\cal G}_i  $. We also observe that $ N_{l} $ is $ \{ {\cal G}_i  : i \geq 0 \} $
adapted for each $ l \geq 0 $, i.e., $ \{ N_{l} \leq m \} \in {\cal G}_m $. Therefore, using ${\mathbf 1}$ 
to denote the indicator function, for any
$ A \in {\cal F}_{T_{N_l}(\bu)} = {\cal G}_{N_l} $, we have
\begin{align*}
\lefteqn{ \E \bigl[ {\mathbf 1} (A)  \sum_{ i = N_l +1}^{ N_{l+1} } Y^{({\bu})}_i \bigr] } \\
& = \E \bigl[ {\mathbf 1} (A)  \sum_{ n_l = 1}^{\infty } \sum_{ m = 1}^{\infty }
{\mathbf 1} ( N_l = n_l) {\mathbf 1} ( N_{l+1} = n_l +m )
 \sum_{ i = 1}^{ m  }  Y^{({\bu})}_{n_l +i} \bigr]   \\
& = \E \bigl[ {\mathbf 1} (A)  \sum_{ n_l = 1}^{\infty } \sum_{ m = 1}^{\infty }
{\mathbf 1} ( N_l = n_l) {\mathbf 1} ( N_{l+1} \geq  n_l +m )
  Y^{({\bu})}_{n_l +m} \bigr]  \\
& = \sum_{ n_l = 1}^{\infty } \sum_{ m = 1}^{\infty } \E \bigl[ {\mathbf 1} (A)
{\mathbf 1} ( N_l = n_l)[ 1 - {\mathbf 1} ( N_{l+1} \leq  n_l +m -1 )  ]
  Y^{({\bu})}_{n_l +m} \bigr]  \\
& = \sum_{ n_l = 1}^{\infty } \sum_{ m = 1}^{\infty } \E \Bigl[ \E \bigl[ {\mathbf 1} (A)
{\mathbf 1} ( N_l = n_l)[ 1 - {\mathbf 1} ( N_{l+1} \leq  n_l +m -1 )  ]
  Y^{({\bu})}_{n_l +m} \mid   {\cal G}_{n_l+m -1 } \bigr]  \Bigr]\\
& = \sum_{ n_l = 1}^{\infty } \sum_{ m = 1}^{\infty } \E \Bigl[  {\mathbf 1} (A)
 {\mathbf 1} ( N_l = n_l)[ 1 - {\mathbf 1} ( N_{l+1} \leq  n_l +m -1 )  ]
 \E \bigl[ Y^{({\bu})}_{n_l +m} \mid   {\cal G}_{n_l+m -1 } \bigr]  \Bigr]\\
& = \sum_{ n_l = 1}^{\infty } \sum_{ m = 1}^{\infty } \E \Bigl[  {\mathbf 1} (A)
 {\mathbf 1} ( N_l = n_l)[ 1 - {\mathbf 1} ( N_{l+1} \leq  n_l +m -1 )  ]
 \E \bigl[ Y^{({\bu})}_{n_l +m} \bigr]  \Bigr] \\
 & = 0 .
\end{align*}
In the above, we have used Fubini's theorem to interchange the expectation and summation.
% Observe that $ \bigl| {\mathbf 1} (A)  \sum_{ i = N_l +1}^{ N_{l+1} } Y^{({\bu})}_i \bigr|
% \leq   \sum_{ i = N_l + 1}^{ N_{l+1} } | Y^{({\bu})}_i |  \leq W_{l+1} (\bu, \bv)$. 
Observe that $  {\mathbf 1} (A)  \sum_{ i = N_l +1}^{ N_{l+1} } |Y^{({\bu})}_i |
\leq   W_{l+1} (\bu, \bv)$.
Hence,  using Proposition \ref{prop:WidthTail}, 
$ \E ( {\mathbf 1} (A)  \sum_{ i = N_l +1}^{ N_{l+1} }| Y^{({\bu})}_i |)
< \infty $. 
% Using Lemma \ref{lem:RandomSumHavingExpTail},
% we conclude that for some $ \alpha > 0$, $ \E \bigl( \exp (\alpha  \sum_{ i = 1}^{ N_{l+1} } | Y^{({\bu})}_i |) \bigr)
% < \infty$.   
In the above we have also used property (c) of Remark  \ref{rem:FromSinglePoint}
and the fact that, since   $ A \in {\cal G}_{N_l} $,
 $ A \cap \{  N_l = n_l \} \in {\cal G}_{n_l} \subseteq
{\cal G}_{n_l+m -1 }  $ for all $ m \geq 1 $. Also, $ \{ N_{l+1} \leq  n_l +m -1 \}
\in {\cal G}_{n_l+m -1 }  $ and $  Y^{({\bu})}_{n_l +m}
$ is independent of $ {\cal G}_{n_l+m -1 }  $. \qed

\subsection{Independent processes}
\label{subSec:IndPr}
In this subsection, we describe  {\em simultaneous regenerations of two independent paths}.
This will be used to approximate the paths at simultaneous regenerations of joint paths when
the starting points are far apart. We start with  a result (Lemma \ref{lem:IndRegen}) about renewal processes,
which is proved in the Appendix.
% The idea of the proof
% is similar to that of Proposition \ref{prop:RecTimeExpTail}, with the maximum residual life for the components
% playing the role of the height in the proof of Proposition \ref{prop:RecTimeExpTail}.

Let $ \{ \xi_n^{(1)} : n \geq 1 \} $ and $ \{ \xi_n^{(2)} : n \geq 1 \} $ be two
independent collections of i.i.d. inter-arrival times (positive integer
valued random variables) with $ \P (  \xi_1^{(1)} = j ) =  f_j^{(1)}  $ and $ \P (  \xi_1^{(2)} = j ) = f_j^{(2)}  $.
 We assume that, for any $ m \geq 1$,
$  \max \{ \P ( \xi_1^{(1)} \geq m) , \P ( \xi_1^{(2)} \geq m) \} \leq C_9 \exp ( -C_{10} m ) $
where $C_9$ and $ C_{10} $ are positive constants and $ f_1^{(1)} f_1^{(2)} > 0 $.
Let $ S_0^{(1)}  := 0 =:  S_0^{(2)} $ and, for $ n \geq 1, S_n^{(1)} := \sum_{j=1}^n
\xi_j^{(1)} $ and $S_n^{(2)} := \sum_{j=1}^n \xi_j^{(2)} $.
For any $ n \geq 1 $ and $ i = 1,2$,  define the residual life of the $i$ th component at time $n$ by
\begin{equation}
R_n^{(i)} := \inf \{  S_k^{(i)} : S_k^{(i)}  \geq n \} - n.
\end{equation}
We consider the joint residual process $ ( R_n^{(i)}: i = 1,2 ) $ and
define
\begin{equation*}
\tau^R := \inf \{ n \geq 1 : R_n^{(1)} = R_n^{(2)}= 0  \}.
\end{equation*}
\begin{lemma}
\label{lem:IndRegen}
For any $ n \geq 1 $, we have
\begin{equation*}
 \P ( \tau^R \geq n ) \leq C_{11} \exp ( - C_{12}  n )
\end{equation*}
where $ C_{11} $ and $ C_{12} $ are positive constants, depending
on the distribution of $ \xi_n^{(i)} $'s only.
\end{lemma}

%\begin{remark}
%It should be noted that the condition $ \P ( \xi_1^{(i)} = 1 ) > 0 $ is not necessary and
%can be replaced by the condition that that the renewal processes, given by the inter-arrival
%times $\xi_n^{(i)} $, are aperiodic.
%\end{remark}

Now we fix two vertices $ \bu$ and $\bv$ with $ \bu (d) = \bv (d) $ and consider two independent constructions of the
marginal paths, one starting from $ \bu$ and the other from $ \bv $.
More precisely, we start with two independent collections of uniform i.i.d. random variables,
$\{U^{\bu}_{\bw}: \bw \in \Z^d\}$ and $\{U^{\bv}_{\bw}: \bw \in \Z^d\}$.  Now, as in Section \ref{sec:ConstructionOfProcess},
we start two paths, one from $ \bu$, constructed using the collection $\{U^{\bu}_{\bw}: \bw \in \Z^d\}$, and
the other from $ \bv $, constructed using the collection $\{U^{\bv}_{\bw}: \bw \in \Z^d\}$.
We denote these respective paths  by $  \{ g_n^{(\text{Ind})} (\bu) : n \geq 0 \} $ and $  \{ g_n^{(\text{Ind})} (\bv) : n \geq 0 \} $.
The above processes being independent, we have two independent collections of  
regeneration times,  $ \{ T_l^{(\text{Ind})} ({\bu}):
l \geq 0 \} $ and  $ \{ T_l^{(\text{Ind})} ({\bv}):
l \geq 0 \} $ (see equation  (\ref{def:TimekthReg}) for definition).
% Since the collection $ \{ T_l^{(\text{Ind})} ({\bu^i}):
% l \geq 0 \} $ uses only the random variables $\{U^i_{\bw}: \bw \in \Z^d\}$, these families
% are independent.
As mentioned in  Remark \ref{rem:FromSinglePoint},
for a single starting point, the distribution of the collection $ \{ T_l^{(\text{Ind})} ({\bu}): l \geq 0 \} $ or
$ \{ T_l^{(\text{Ind})} ({\bv}): l \geq 0 \}$ does not depend on the starting point, and  is an independent
copy of $ \{ T_l({\bf 0}): l \geq 0 \} $.

Take  $
R_n^{(1)} = \inf \{ T_l^{(\text{Ind})} ({\bu})   : T_l^{(\text{Ind})} ({\bu})  \geq n \} - n $
and $ R_n^{(2)} = \inf \{ T_l^{(\text{Ind})} ({\bv})   : T_l^{(\text{Ind})} ({\bv})  \geq n \} - n $.
Note here that
$ \{ T_{l+1}^{(\text{Ind})} ({\bu})  - T_l^{(\text{Ind})} ({\bu})  : l \geq 0 \} $ and
$ \{ T_{l+1}^{(\text{Ind})} ({\bv})  - T_l^{(\text{Ind})} ({\bv})  : l \geq 0 \} $
are two independent collections of i.i.d. random variables, which play the respective roles of $ \{ \xi_l^{(1)} : l \geq 1 \} $
and  $ \{ \xi_l^{(2)} : l \geq 1 \} $  of Lemma \ref{lem:IndRegen}.
Set, $ T_0^{(\text{Ind})} (\bu,  \bv) = 0 $ and, for $ l \geq 0 $,
\begin{equation}
\label{def:indRegenTime}
T_{l+1}^{(\text{Ind})}  (\bu,  \bv)   := \inf \{ n > T_{l}^{(\text{Ind})}
 (\bu, \bv)  : R_n^{(1)}  =  R_n^{(2)}  = 0  \}.
\end{equation}
We call $ T_{l}^{(\text{Ind})}  (\bu,\bv) $, the time for the $ l$ th
{\em simultaneous regeneration time of two independent paths}.

Applying Lemma \ref{lem:IndRegen}  and observing that each
$ T_{l}^{(\text{Ind})}  (\bu, \bv)  $ represents the occurrence of a renewal event,
we obtain the following proposition.
\begin{prop}
\label{prop:IndIIDOcc}
$ \{ T_{l+1}^{{\rm (Ind)}}   (\bu, \bv)  -
T_{l}^{{\rm (Ind)}}   (\bu, \bv) : l \geq 0\} $
is an i.i.d. sequence of random variables
taking values in $ \{ 1, 2, 3, \dotsc \} $ and, for all $ n \geq 1 $,
\begin{equation}
 \label{eqn:TStarexpest}
\P ( T_1^{{\rm (Ind)}}   (\bu, \bv) \geq n ) \leq
C_{13}  \exp ( - C_{14} n )
\end{equation}
where $ C_{13} $ and $ C_{14} $ are positive constants.
\end{prop}

% Now, for $ \bu, \bv \in Z^d $ with $ \bu (d) = \bv (d)$,  let $\{g_n^{(\text{Ind})} ({\bu^i}): n \geq 0\}$, $i=1,2$,
% be the independent versions of the paths  starting from $\bu $ and $ \bv$ respectively, as described above.
% Also, let $\{g_{\tau_l ({\bu^i})}^{(\text{Ind})} ({\bu^i}): l \geq 0\}$ be the
% $i$ th  process evaluated at its regeneration steps $\tau_l (\bu^i) $, $l \geq 0$.
% Let $T_{l}^{(\text{Ind})}   (\bu, \bv)$ be the $l$ th simultaneous regeneration time as defined in (\ref{def:indRegenTime}).

Any simultaneous regeneration time of two independent paths $T_{l}^{(\text{Ind})}   (\bu, \bv)$
is also a regeneration time for each of  the marginal processes. Therefore, as in (\ref{defn:TNL}), 
for every $ l \geq 0$, we have
\begin{equation*}
T_{l}^{(\text{Ind})}   (\bu, \bv )  = T^{(\text{Ind})}_{N_l^{\bu}}( \bu) =  T^{(\text{Ind})}_{N_l^{\bv}}( \bv)
\end{equation*}
for some sequences $ N_l^{\bu} (= N_l^{\bu} (\bu, \bv )) $ and $ N_l^{\bv} (= N_l^{\bv} (\bu, \bv )) $ with $
N_0^{\bu} = N_0^{\bv} = 0$. 
%  . suppose that
%  Suppose $N_l(i)$ ($i = 1, 2$, $l \geq 0$)  is such that the $l$ th joint
% regeneration coincides with the $N_l(i)$ th regeneration of the $i$ th process, i.e.
% \begin{equation*}
% T_{l}^{(\text{Ind})}   (\bu, \bv) = T^{(\text{Ind})}_{N_l(i)}({\bu^i}), \; i = 1, 2
% \end{equation*}
% where $ T^{(\text{Ind})}_{N_l(i)}({\bu^i})$ is the $N_l(i)$ th regeneration time of the
% $i$ th process starting only from ${\bu^i}$.

As in (\ref{r-W2}) consider the width of the explored region between the
$l-1$ and $ l $ th regenerations of both the independent  processes,
\begin{equation}
\label{r-Wind}
W_l^{{\rm (Ind)}}  (\bu,  \bv ) :=
\sum_{ t = N_{l-1}^{\bu} + 1}^{ N_l^{\bu} } W_t^{\rm (Ind)} ( \bu) +
 \sum_{ t = N_{l-1}^{\bv} + 1}^{ N_l^{\bv} } W_t^{\rm (Ind)} ( \bv)
% ||g_t({\bu^i})-g_{t-1}({\bu^i})||_1,
\end{equation}
where $ \{ W_l^{\rm (Ind)} ( \bu) : l \geq 1 \}  $  and $  \{ W_l^{\rm (Ind)} ( \bv) : l \geq 1 \}  $ are
the explored width processes associated with the processes $ \{ g_l^{(\text{Ind})} (\bu) : l \geq 0 \} $
and  $ \{ g_l^{(\text{Ind})} (\bv) : l \geq 0 \} $ respectively. 
% From Proposition \ref{prop:IndIIDOcc} and
% the fact that $ \{ g_l^{(\text{Ind})} (\bu) : l \geq 0 \} $  and $ \{ g_l^{(\text{Ind})} (\bv) : l \geq 0 \} $ are
% independent families, 
Observe that $ \{ (N_{l+1}^{\bu} - N_l^{\bu}, N_{l+1}^{\bv} - N_l^{\bv} ) : 
l \geq 0 \} $ is a sequence of i.i.d. random variables,
and hence, $ \{ W_l^{{\rm (Ind)}}  (\bu,  \bv )  : l \geq
1 \} $ is also a sequence of i.i.d. random variables.

Since $ \max \{ N_{1}^{\bu}, N_{1}^{\bv} \} \leq T_{1}^{(\text{Ind})}  (\bu, \bv) $ and
$T_{1}^{(\text{Ind})}   (\bu, \bv) - T_{0}^{(\text{Ind})} (\bu, \bv) =  T_{1}^{(\text{Ind})}   (\bu, \bv)$ 
satisfies (\ref{eqn:TStarexpest}),
using Lemma \ref{lem:RandomSumHavingExpTail},    we conclude, for any $ l \geq 1$
\begin{equation}
\label{eqn:IndWidthTailBound}
\P(W_l^{{\rm (Ind)}}  (\bu, \bv )  \geq n) \leq C_{15}  \exp ( - C_{16} n )
\end{equation}
where $ C_{15} $ and $ C_{16} $ are positive constants.

From (\ref{Ydef}), we can write, for any $ l \geq 1$,
\begin{equation*}
\overline{g}^{(\text{Ind})}_{\tau_{N_l^{\bu}}}({\bu}) = \overline{\bu} +
\sum_{ t = 1}^{ N_l^{\bu} } Y_t^{(\bu)} \quad \text{ and }\quad  \overline{g}^{(\text{Ind})}_{\tau_{N_l^{\bv}}}({\bv}) = \overline{\bv} +
\sum_{ t = 1}^{ N_l^{\bv} } Y_t^{(\bv)}.
\end{equation*}
At each simultaneous regeneration time $T_{l}^{(\text{Ind})}  (\bu, \bv)  $,
the $d$ th co-ordinates of both the paths coincide  and equal $T_{l}^{(\text{Ind})}
(\bu, \bv)  $. We consider the first $d-1$ co-ordinates
of these paths and define
\begin{equation}
  \label{ZDef}
\begin{split}
\psi_l^{\bu} & := \overline{g}^{(\text{Ind})}_{\tau_{N_l^{\bu}}}({\bu}) -
 \overline{g}^{(\text{Ind})}_{\tau_{N_{l-1}^{\bu}}}({\bu})
=  \sum_{ t = N_{l-1}^{\bu} + 1}^{ N_l^{\bu}} Y_t^{(\bu)};  \\
\psi_l^{\bv} & := \overline{g}^{(\text{Ind})}_{\tau_{N_l^{\bv}}}({\bv}) -
 \overline{g}^{(\text{Ind})}_{\tau_{N_{l-1}^{\bv}}}({\bv})
=  \sum_{ t = N_{l-1}^{\bv} + 1}^{ N_l^{\bv}} Y_t^{(\bv)}.
\end{split}
\end{equation}
The process $ (\psi_l^{\bu}, \psi_l^{\bv}) $ represents the increment between the $ (l-1)$ th 
and  $ l$ th simultaneous regeneration times in the first $(d-1) $ co-ordinates of the independent
paths starting from $ \bu $ and $ \bv $ respectively. We list  the properties of $ (\psi_l^{\bu}, 
\psi_l^{\bv})$ in the following proposition.

\begin{prop}
\label{prop:PropertyIncrementsForIndPaths}
The family $ \{ (\psi_l^{\bu}, \psi_l^{\bv} )   : l \geq 1\} $ is an i.i.d. collection of random variables
taking values in $\Z^{2(d-1)}$ satisfying the following properties.
 \begin{itemize}
\item[(a)]  For any $ n \geq 1$,
\begin{equation*}
\P( || \psi_1^{\bu} ||_{1, d-1} + || \psi_1^{\bv} ||_{1, d-1}  \geq n ) \leq C_{15}  \exp ( - C_{16} n )
\end{equation*}
where $ C_{15} $ and $ C_{16} $ are as in (\ref{eqn:IndWidthTailBound}).

\item[(b)] The marginal distribution of each co-ordinate of $  \psi_1^{\bu} $ as well as
$   \psi_1^{\bv} $ is symmetric. Further, they are all same. More precisely, with $\psi_l^{\bu}(j)$ and
$\psi_l^{\bv}(j)$ being the $j$ th co-ordinate of $\psi_l^{\bu}$ and $\psi_l^{\bv}$ respectively,
$ \P \bigl( \psi_1^{\bu}(j) = r \bigr) = \P \bigl( \psi_1^{\bu}(j) = -r \bigr) =
\P \bigl( \psi_1^{\bv}(j) = r \bigr) = \P \bigl( \psi_1^{\bv}(j) = -r \bigr) = \P \bigl( \psi_1^{\bu}(1) = r \bigr)
$  for all  $ j = 1, 2, \dotsc, d-1$ and  $ r \geq 1$.

% Furthermore,
% $ \P \bigl( \psi_l^{(i)}(j) = r \bigr) $ is independent of $ i = 1, \dotsc, k$ and
% $ j = 1, \dotsc, d-1$ .
\item[(c)]  $\E\Bigl[ \bigr( \psi_1^{\bu}(i)\bigr)^{m_1}  
\bigr( \psi_1^{\bv}(j)\bigr)^{m_2}\Bigl] $ is
 independent of $ i,j$ and depends only on $m_1, m_2$ and
 $\E\Bigl[ \bigr( \psi_1^{\bu}(i)\bigr)^{m_1}  \bigr( \psi_1^{\bv}(j)\bigr)^{m_2}\Bigl] = 0 $
if at least one of $ m_1, m_2 $ is odd.
\end{itemize}
\end{prop}

\noindent{\bf Proof}: Noting  that   $T_{l}^{(\text{Ind})}    (\bu, \bv)  $ represents the occurrence
of a  renewal event and using the fact that the families of i.i.d. random variables $ \{  Y_t^{(\bu)} : t \geq 1 \} $ and
$ \{  Y_t^{(\bv)} : t \geq 1 \} $ are independent, we conclude that $ \{  (\psi_l^{\bu}, \psi_l^{\bv} )  : l \geq 1 \}$
is a family of i.i.d. random variables. Further,
$ || \psi_1^{(\bu)} ||_{1, d-1} +  ||  \psi_1^{(\bv)} ||_{1, d-1} \leq W_1^{{\rm (Ind)}}  (\bu, \bv ) $. Using
(\ref{eqn:IndWidthTailBound}), we conclude part (a).

Using the fact that $ Y_t^{(\bu)} $ is symmetric in each
co-ordinate, we conclude $ \P \bigl( \psi_1^{\bu}(j) = r \bigr) = \P \bigl( \psi_1^{\bu}(j) = -r \bigr) =
\P \bigl( \psi_1^{\bv}(j) = r \bigr) = \P \bigl( \psi_1^{\bv}(j) = -r \bigr) $ for any $ r \geq 1$ for all $ j = 1,
\dotsc, d-1$. Further,
$ Y_t ^{(\bu)} $ is rotation invariant along any of the first $d-1$ co-ordinates and
hence $  \P \bigl( \psi_1^{\bu}(j) = r \bigr) =  \P \bigl( \psi_1^{\bu}(1) = r \bigr)$
for $ j = 2, 3, \dotsc, d-1$. This proves (b).

To study the joint distribution of $ ( \psi_1^{\bu}(i),  \psi_1^{\bv}(j) ) $, we observe
that we may apply the rotation operator independently for both families $ \{ U_{\bw}^{\bu} \} $ 
and $ \{ U_{\bw}^{\bv} \} $ of
uniform random variables, so that the $ i$ th co-ordinate after rotation becomes the first 
co-ordinate for $ \{ U_{\bw}^{\bu} \} $
and $ j $ th co-ordinate after rotation becomes the first co-ordinate for $ \{ U_{\bw}^{\bv} \} $.
The distribution of $ (Y_t^{(\bu)} (i), Y_t^{(\bv)}(j)) $ after rotation remains unchanged,
which implies that the joint distribution of $  ( \psi_1^{\bu}(i),  \psi_1^{\bv}(j) )  $  is 
independent of $ i,j $.
Thus, $ \E\Bigl[ \bigr( \psi_1^{\bu}(i)\bigr)^{m_1}  \bigr( \psi_1^{\bv}(j)\bigr)^{m_2}\Bigl]  $
is independent of choice of $ i,j$.
If we fix the realizations for one family of uniform random variables and reflect the realizations 
of the other family along some co-ordinate,
the distribution of $  (Y_t^{(\bu)} (i), Y_t^{(\bv)}(j) ) $ remains unaltered. Therefore, 
we have $ \P \bigl(  (\psi_1^{\bu}(i),  \psi_1^{\bv}(j))
= (l, k) \bigr) = \P \bigl(  (\psi_1^{\bu}(i),  \psi_1^{\bv}(j)) = (l, - k) \bigr) $ for
any $ l,k \in \Z $. This proves (c). \qed

\begin{remark}
\label{rem:DiffsMakeRW}
Let $  \xi_l := \psi_l^{\bv}  -  \psi_l^{\bu} $ for $l \geq 1$ and set
$ S_0  := \overline{\bv} - \overline{\bu} $ and for $ n \geq 1,  S_n := S_0 + \sum_{l=1}^n  \xi_l  $.
From the above Proposition, we conclude that the
process $ \{  S_n : n \geq 0 \} $ is an isotropic, symmetric random walk starting from $\overline{\bv} - \overline{\bu}  $
on $ \Z^{d-1}$
with i.i.d. steps, each step having exponentially decaying tail probability. 
Clearly, $ \P (  \xi_l  = {\bf 0}) \geq p^2 $ so that the random walk is also aperiodic.
This will be used in the next section.
\end{remark}

%
%
% \item[(ii)] the increments of the $\overline{g}^i_n$ between its regenerations given as in (\ref{Ydef}) by
% $Y_{l+1}(i) := \overline{g}^i_{\tau_{l+1}} -  \overline{g}^i_{\tau_{l}}$, $ l \geq 0$ are also i.i.d. both
% as a process in $l$ as well as a collection of $k$ processes,
% \item[(iii)] $\psi_1^{(i)} = \sum_{j=1}^{N_1(i)} Y_j(i)$.
% \end{itemize}
% Using the inherent symmetries of the marginals of $\psi_l^{(i)}$ (and carrying out calculations,
% similar to that in the proof of Proposition \ref{prop:RWMartingale},
% we obtain the following.
% \begin{itemize}
% \end{itemize}

\subsection{Coupling of joint process and independent process}
\label{Coupling}
In this subsection, we describe a coupling
of two independent paths with the joint paths starting from $ \bu, \bv$ with
$ \bu (d)  = \bv (d) $.
Without loss of generality, we may assume $  \bu (d) = 0 $. Set
$ d_{\min} :=  ||  \bu -  \bv ||_1 $.
As in the last subsection, we start with two independent  collections of i.i.d.
uniform random variables $ \{ U_{\bw}^{\bu} : \bw \in \Z^d, \bw (d) > 0 \} $
and $ \{ U_{\bw}^{\bv} : \bw \in \Z^d, \bw (d) > 0\} $, and
construct two independent paths, $\{g_n^{(\text{Ind})} ({\bu}) : n \geq 0\}$ and
$\{g_n^{(\text{Ind})} ({\bv}) : n \geq 0\}$,  starting from  $ \bu $ and $ \bv $ respectively.

Fix $ r < d_{\min}/2 $ and another independent collection of uniform random variables
$ \{ U_{\bw} : \bw \in \Z^d, \bw (d) > 0  \} $. We define a new collection
of uniform random variables $ \{ {\tilde{U}}_{\bw} : \bw \in \Z^d,
\bw(d) >  0 \}  $ by
\begin{equation*}
 {\tilde{U}}_{\bw}   := \begin{cases}
 U_{\bw}^{\bu}  & \text{ if } || \overline{\bw} - \overline{\bu} ||_{1, d-1} < r \\
 U_{\bw}^{\bv}  & \text{ if } || \overline{\bw} - \overline{\bv } ||_{1, d-1} < r \\
 U_{\bw} & \text{ otherwise.}
  \end{cases}
  \end{equation*}
Using the collection $ \{ \tilde{U}_{\bw} : \bw \in \Z^d, \bw (d) > 0  \} $, we construct the
joint paths (as in Section 2)
from the points $ \bu $ and $  \bv $ until its first simultaneous regeneration time
$T_1(\bu, \bv)$ of joint paths from $\bu$ and $ \bv$.

\begin{figure}[htb]

\begin{center}\leavevmode

% Generated with LaTeXDraw 2.0.8
% Mon Nov 19 21:14:50 IST 2012
% \usepackage[usenames,dvipsnames]{pstricks}
% \usepackage{epsfig}
% \usepackage{pst-grad} % For gradients
% \usepackage{pst-plot} % For axes

\begin{pspicture}(2.5,-1.65)(16.5,1.75)
\psframe[linewidth=0.025999999,dimen=outer,fillstyle=solid,fillcolor=lightgray](9.075469,1.493828)(3.4954686,-1.3261719)
\psframe[linewidth=0.025999999,dimen=outer,fillstyle=solid,fillcolor=lightgray](9.075469,-1.3261719)(9.075469,-1.3461719)
\psframe[linewidth=0.025999999,dimen=outer,fillstyle=solid,fillcolor=lightgray](16.115469,1.4738281)(9.895469,-1.3461719)
\psdots[dotsize=0.12](6.1192183,-1.3396094)
\psdots[dotsize=0.12](12.759218,-1.3196094)
\psline[linewidth=0.04cm,linestyle=dashed,dash=0.16cm 0.16cm](6.1392183,0.26039064)(4.6992188,-1.0996093)
\psline[linewidth=0.04cm,linestyle=dashed,dash=0.16cm 0.16cm](6.1592183,0.26039064)(7.5192184,-0.9196093)
\psline[linewidth=0.04cm,linestyle=dashed,dash=0.16cm 0.16cm](12.739219,-0.3396094)(11.939218,-1.0996093)
\psline[linewidth=0.04cm,linestyle=dashed,dash=0.16cm 0.16cm](12.739219,-0.3596094)(13.479218,-1.1196094)
\psline[linewidth=0.04cm,linestyle=dashed,dash=0.16cm 0.16cm](6.9592185,0.6603907)(5.9392185,-0.3196094)
\psline[linewidth=0.04cm,linestyle=dashed,dash=0.16cm 0.16cm](6.9592185,0.6603907)(7.8992186,-0.2996093)
\psline[linewidth=0.04cm,linestyle=dashed,dash=0.16cm 0.16cm](13.5192175,0.0403906)(12.579219,-0.8796094)
\psline[linewidth=0.04cm,linestyle=dashed,dash=0.16cm 0.16cm](13.5192175,0.0203906)(14.279219,-0.7396094)
\psline[linewidth=0.04cm,linestyle=dashed,dash=0.16cm 0.16cm](7.5792184,1.0603907)(6.779218,0.30039066)
\psline[linewidth=0.04cm,linestyle=dashed,dash=0.16cm 0.16cm](7.5792184,1.0403906)(8.319218,0.2803907)
\psline[linewidth=0.04cm](6.1592183,-1.2796094)(6.9592185,-0.4196093)
\psline[linewidth=0.04cm](7.5392184,0.0603907)(6.9992185,-0.3996093)
\psline[linewidth=0.04cm](7.5392184,0.0603907)(7.3592186,0.8403906)
\psline[linewidth=0.04cm](7.3592186,0.8803907)(7.3592186,1.4803905)
\psdots[dotsize=0.12](6.9592185,-0.4396093)
\psdots[dotsize=0.12](7.5592184,0.0603907)
\psdots[dotsize=0.12](7.3592186,0.8603907)
\psdots[dotsize=0.12](7.3592186,1.4603906)
\psdots[dotsize=0.12](13.459219,-1.0996093)
\psdots[dotsize=0.12](13.5192175,0.0203906)
\psline[linewidth=0.04cm](12.819218,-1.3196094)(13.439218,-1.0996093)
\psline[linewidth=0.04cm](13.479218,-1.0596095)(13.5192175,0.0403906)
\psline[linewidth=0.03cm,linestyle=dashed,dash=0.16cm 0.16cm](7.299218,1.4003907)(6.799218,0.9803906)
\psline[linewidth=0.03cm,linestyle=dashed,dash=0.16cm 0.16cm](7.3992186,1.4403907)(7.7392187,1.0803907)
\psdots[dotsize=0.12](12.959219,0.8603907)
\psline[linewidth=0.03cm,linestyle=dashed,dash=0.16cm 0.16cm](13.499219,1.4003907)(12.239219,0.20039064)
\psline[linewidth=0.03cm,linestyle=dashed,dash=0.16cm 0.16cm](13.539219,1.4403907)(14.719218,0.50039065)
\psline[linewidth=0.03cm,linestyle=dashed,dash=0.16cm 0.16cm](14.719218,0.50039065)(14.959219,0.2803907)
\psline[linewidth=0.05cm](13.499219,0.0603907)(12.959219,0.8203907)
\usefont{T1}{ptm}{m}{n}
\rput(5.90, -1.70){$\bu$}
\usefont{T1}{ptm}{m}{n}
\rput(12.55,-1.70){$\bv$}
\psdots[dotsize=0.12](12.953438,1.493828)
\psline[linewidth=0.04cm,linestyle=dashed,dash=0.16cm 0.16cm](12.933437,1.4538281)(12.533438,1.0738281)
\psline[linewidth=0.04cm,linestyle=dashed,dash=0.16cm 0.16cm](13.013437,1.493828)(13.413438,1.0938281)
\psline[linewidth=0.04cm](12.973437,0.9138281)(12.973437,1.493828)
\usefont{T1}{ptm}{m}{n}
%\rput(5.2875,0.5638281){$D_{0}^{\mathbf{0},\bx}(\mathbf{0})$}
\usefont{T1}{ptm}{m}{n}
%\rput(11.041406,0.3438281){$D_{0}^{\mathbf{0},\bx}(\bx)$}
\usefont{T1}{ptm}{m}{n}
\rput(4.80,0.95){$T_1^{(\text{Ind})} (\bu, \bv)$}
\usefont{T1}{ptm}{m}{n}
\rput(8.0,1.90){$g^{(\text{Ind})}_{\tau_{N_1^{\bu}}}(\bu)$}
\usefont{T1}{ptm}{m}{n}
\rput(13.5,1.90){$g^{(\text{Ind})}_{\tau_{N_1^{\bv} }}(\bv)$}
\psline[linewidth=0.025999999cm,linestyle=dashed,dash=0.16cm 0.16cm,arrowsize=0.1cm 2.0,arrowlength=1.4,arrowinset=0.4]{<->}(12.835468,-1.6222657)(15.995468,-1.5422657)
\usefont{T1}{ptm}{m}{n}
\rput(14.396094,-1.9){$r$}
\psline[linewidth=0.025999999cm,linestyle=dashed,dash=0.16cm 0.16cm,arrowsize=0.1cm 2.0,arrowlength=1.4,arrowinset=0.4]{<->}(6.135468,-1.6222657)(9.095468,-1.5422657)
\usefont{T1}{ptm}{m}{n}
\rput(7.80,-1.9){$r$}
\psline[linewidth=0.024cm,linestyle=dashed,dash=0.16cm 0.16cm,arrowsize=0.1cm 2.0,arrowlength=1.4,arrowinset=0.4]{<->}(3.25,1.4672656)(3.25,-1.3527344)
\end{pspicture}
\caption{The shaded regions represent part of the cylinders (up to $T_1^{(\text{Ind})} (\bu, \bv)$) of width $ r $ around $ \bu $ and $ \bv$. In the left cylinder
we use the collection $ \{ U^{\bu}_{\bw} \} $, on the right cylinder
we use the collection $ \{ U^{\bv}_{\bw} \} $ and in the remaining region, we use $ \{ U_{\bw} \} $. }
\label{Width}

\end{center}
\end{figure}
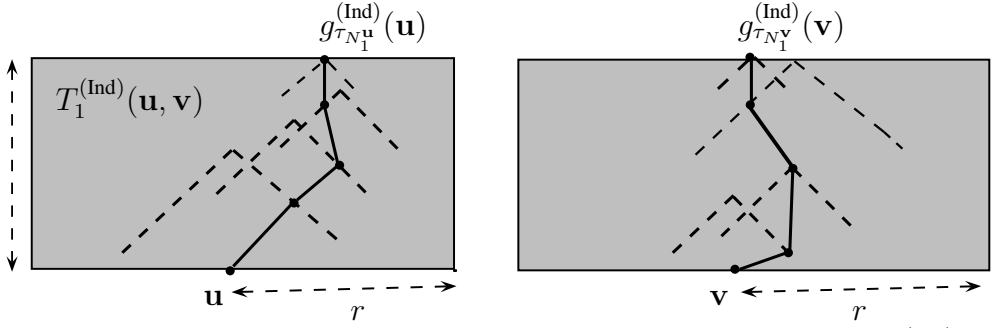

Now,  as defined in  (\ref{def:indRegenTime}), let $T_{1}^{(\text{Ind})} (\bu, \bv)
$ be the first simultaneous
regeneration time of two independent processes and $N_1^{\bu}$ and $N_1^{\bv}$ be the number of individual
regenerations of the independent paths starting from $\bu$ and $ \bv$ respectively.
With the width of the explored region for two independent processes,
as defined in (\ref{r-Wind}),
we consider the event where the total width of the explored region until the first simultaneous
regeneration time of the two independent paths is less than $r$ (see Figure \ref{Width}).
More precisely, we consider the event
\begin{align*}
 A^{\text{Good}} (r)  := \Bigl\{ W_1^{{\rm (Ind)}} (\bu, \bv) \leq r \Bigr\}.
\end{align*}

% \documentclass[a4paper,12pt]{article}
% \usepackage{pst-all}
% \begin{document}

% \begin{figure}[h]
% \begin{center}\leavevmode
%includegraphics[width=.5\textwidth]{hist1}

% Generated with LaTeXDraw 2.0.8
% Tue Aug 28 13:33:03 IST 2012
% \usepackage[usenames,dvipsnames]{pstricks}
% \usepackage{epsfig}
% \usepackage{pst-grad} % For gradients
% \usepackage{pst-plot} % For axes
% Generated with LaTeXDraw 2.0.8
% Sun Sep 30 11:03:29 IST 2012
% \usepackage[usenames,dvipsnames]{pstricks}
% \usepackage{epsfig}
% \usepackage{pst-grad} % For gradients
% \usepackage{pst-plot} % For axes

%\end{pspicture}
%\input{construction.tex}

%\end{document}

%\end{document}
On the event $A^{\text{Good}} (r) $ consider the following segments of paths:-\\
(i) the joint path process $(g_n(\bu) ,
 g_n (\bv) ) $
started simultaneously from $(\bu, \bv)$, using the collection $ \{ {\tilde{U}}_{\bw}  : \bw \in \Z^d,
\bw(d) >  0 \} $, until the first simultaneous regeneration time
$ T_1 (\bu, \bv) $ of joint paths; and \\
(ii) the independent paths
$ \{ g^{(\text{Ind})} _n({\bu}) : n \geq 0 \}$, constructed using only the
collection $  \{ U_{\bw}^{\bu} : \bw \in \Z^d, \bw(d) >  0 \} $,
and $ \{ g^{(\text{Ind})} _n({\bv}) : n \geq 0 \}$, constructed using only the
collection $  \{ U_{\bw}^{\bv} : \bw \in \Z^d, \bw(d) >  0 \} $, until the first
simultaneous regeneration of the independent paths $T_{1}^{(\text{Ind})}
(\bu, \bv)  $.\\
These two segments coincide as geometric paths (i.e. line segments in $\R^2$)
although they may be indexed differently.
% in a different fashion than that in the independent paths, but same
%paths (as line segments in $\R^2$) will be traced in both of them.
Therefore,  we must have
\begin{equation*}
T_{1}^{(\text{Ind})} (\bu, \bv)  = T_1 ( \bu,  \bv)
\end{equation*}
and hence we have,
\begin{align}
\label{eqn:CouplingJtAndInd}
& \P\bigl[  \bigl( \overline{g}_{\tau_1 ( \bu,  \bv) }  ( \bu) , \overline{g}_{\tau_1 ( \bu,  \bv) }  ( \bv) \bigr)
=  \bigl( \overline{\bu}+ \psi_1^{\bu},  \overline{\bv}+ \psi_1^{\bv} \bigr) \bigr]  \nonumber  \\
& \geq  \P \bigl(A^{\text{Good}} (r)\big) \nonumber \\
 &   \geq  1 -  C_{15}  \exp ( - C_{16} r   ).
\end{align}
Finally,  using the Markov property, we can use this coupling for each subsequent
joint regeneration step. The new value of $ d_{\min} $ for the $l$ th regeneration
has to be computed from the position of the processes at the $l-1$ th joint regeneration
and the value of $ r $ has to be chosen accordingly.

%\input{3ASRR1707.tex}
% !TEX root = ./DSFrf.tex
% !TEX program = XeLaTeX
\section{Trees and Forest}
In this section we prove Theorem \ref{r-treesforest}. For
$ d = 2, 3 $, we need to prove that for any $ \bu, \bv \in V$, the paths $ \pi^{\bu} $
and $ \pi^{\bv} $ coincide eventually, i.e., $  \pi^{\bu}  (t) =  \pi^{\bv}  (t) $
for all $ t \geq t_0 $ for some $ t_0 < \infty  $.

First, we claim that it is enough to prove that
\begin{equation}
\label{r-meet}
\pi^{\bu} \text{
and } \pi^{\bv} \text{  coincide eventually for } \bu, \bv \in   V \text{ with }
\bu(d) = \bv (d) .
\end{equation}
Indeed, for $\bu, \bv \in   V$ with $\bu(d) < \bv (d) $ we have, from (\ref{r-meet}),
\begin{align*}
 \P \bigl[   \bigcap_{ \bw \in V, \bu (d) = \bw (d) } \{ \text{the paths } \pi^{\bu} \text{ and }
 \pi^{\bw} \text{  coincide eventually} \} \bigr] = 1; \\
 \P  \bigl[    \bigcap_{ \bw^\prime \in V, \bw^\prime (d) = \bv (d) }  \bigl\{ \text{the paths } \pi^{\bv} \text{ and }
 \pi^{\bw^\prime} \text{  coincide eventually} \} \bigr] = 1.
\end{align*}
Further, $ \P \bigl[$there exist  $\bw, \bw^\prime \in V $ with $ \bw (d) = \bu (d), \bw^\prime (d)
= \bv (d),  h(\bw ) = \bw^\prime  \bigr] = 1 $. Since, the intersection of these three events
has probability $1$,  $ \pi^{\bu} $ and $ \pi^{\bv} $ meet.

Now, to prove that for any two vertices $ \bu $ and $ \bv $ with $ \bu (d)
= \bv (d) $, the paths coincide eventually, we show that $  \P ( Z_l (\bu, \bv) = 0 $
for some $l \geq 0) = 1 $ where $Z_l$ is as in (\ref{def:Z_l}). Recall, at the beginning of Section 3,
we had observed that $ \{ Z_l (\bu, \bv) : l \geq 0 \} $ is a Markov chain taking values in $ \Z^{d-1} $
with $ \overline{\bf 0} \in \Z^{d-1}$ being its only  absorbing state.
\subsection{d=2}
\begin{prop}
\label{prop:NonCrossingPaths}
The paths of $ \chi $ do not cross each other almost surely.
\end{prop}

\begin{figure}[!htb]
\centering
% Generated with LaTeXDraw 2.0.8
% Mon Jul 14 10:17:47 IST 2014
% \usepackage[usenames,dvipsnames]{pstricks}
% \usepackage{epsfig}
% \usepackage{pst-grad} % For gradients
% \usepackage{pst-plot} % For axes
% Generated with LaTeXDraw 2.0.8
% Mon Jul 28 17:15:29 IST 2014
% \usepackage[usenames,dvipsnames]{pstricks}
% \usepackage{epsfig}
% \usepackage{pst-grad} % For gradients
% \usepackage{pst-plot} % For axes
\scalebox{.95} % Change this value to rescale the drawing.
{
\begin{pspicture}(0,-4.2)(14.43,4.2)
\definecolor{color117}{rgb}{0.36470588235294116,0.21568627450980393,0.21568627450980393}
\definecolor{color117b}{rgb}{0.39215686274509803,0.2901960784313726,0.2901960784313726}
\pspolygon[linewidth=0.01cm,linestyle=dashed,dash=0.16cm 0.16cm](1.6,2.6)(0.0,1)(3.2,1)
\pspolygon[linewidth=0.01cm,linestyle=dashed,dash=0.16cm 0.16cm](0.8,1.8)(3.2,4.1)(5.6,1.8)
\pspolygon[linewidth=0.01cm,linestyle=dashed,dash=0.16cm 0.16cm](9.6,2.6)(8.0,1)(11.2,1)
\pspolygon[linewidth=0.01cm,linestyle=dashed,dash=0.16cm 0.16cm](10.8,3.8)(7.2,0.2)(14.4,0.2)
\psdots[dotsize=0.14](1.62,0.95375)
\psdots[dotsize=0.14](3.24,1.77375)
\psdots[dotsize=0.14](10.82,0.19375)
\psdots[dotsize=0.14](9.64,0.99375)
\psdots[dotsize=0.14](1.5,2.43375)
%\usefont{T1}{ptm}{m}{n}
\rput(1.3,2.88375){$h(\bu)$}
\psdots[dotsize=0.14](9.48,2.43375)
\rput(9.25,2.75){$h(\bu)$}
\psline[linewidth=0.02cm](1.48,2.39375)(1.6,0.97375)
\psline[linewidth=0.02cm](9.48,2.39375)(9.64,1.01375)
\pspolygon[linewidth=0.01cm,linestyle=dashed,dash=0.16cm 0.16cm](1.64,-1.37)(0.02,-2.97)(3.24,-2.97)
\pspolygon[linewidth=0.01cm,linestyle=dashed,dash=0.16cm 0.16cm](2.42,-1.01)(0.82,-2.61)(4.0,-2.61)
\pspolygon[linewidth=0.01cm,linestyle=dashed,dash=0.16cm 0.16cm](9.22,-1.41)(7.6,-3.02)(10.82,-3.02)
\psdots[dotsize=0.14](1.6,-2.98625)
\psdots[dotsize=0.14](9.22,-3.00625)
\psdots[dotsize=0.14](10.02,-3.42625)
\pspolygon[linewidth=0.01cm,linestyle=dashed,dash=0.16cm 0.16cm](10.02,-1.03)(7.62,-3.4)(12.4,-3.4)
\psdots[dotsize=0.14](9.06,-1.54625)
\rput(8.85,-1.15){$h(\bu)$}
\psdots[dotsize=0.14](1.48,-1.54625)
\rput(1.2,-1.3){$h(\bu)$}
\psline[linewidth=0.02cm](1.46,-1.56625)(1.58,-2.94625)
\psline[linewidth=0.02cm](9.06,-1.54625)(9.22,-2.96625)
% \psline[linewidth=0.02cm,linestyle=dashed,dash=0.16cm 0.16cm](0.26,1.17375)(2.98,1.21375)
% \psline[linewidth=0.02cm,linestyle=dashed,dash=0.16cm 0.16cm](1.08,2.03375)(5.38,2.03375)
% \psline[linewidth=0.02cm,linestyle=dashed,dash=0.16cm 0.16cm](8.28,1.21375)(10.98,1.23375)
% \psline[linewidth=0.02cm,linestyle=dashed,dash=0.16cm 0.16cm](7.46,0.41375)(14.18,0.45375)
% \psline[linewidth=0.02cm,linestyle=dashed,dash=0.16cm 0.16cm](7.86,-3.26625)(12.2,-3.18625)
% \psline[linewidth=0.02cm,linestyle=dashed,dash=0.16cm 0.16cm](7.86,-2.82625)(10.64,-2.82625)
% \psline[linewidth=0.02cm,linestyle=dashed,dash=0.16cm 0.16cm](0.28,-2.78625)(3.02,-2.76625)
% \psline[linewidth=0.02cm,linestyle=dashed,dash=0.16cm 0.16cm](1.1,-2.38625)(3.76,-2.38625)
\psdots[dotsize=0.14](2.42,-2.62625)
\psline[linewidth=0.1cm,linecolor=blue](1.38,2.35375)(1.1,2.05375)
\psline[linewidth=0.1cm,linecolor=blue](1.38,-2.02625)(1.02,-2.36625)
\psline[linewidth=0.1cm,linecolor=blue](9.32,2.33375)(8.22,1.23375)
\psline[linewidth=0.1cm,linecolor=blue](9.02,-2.04625)(8.26,-2.82625)
\usefont{T1}{ptm}{m}{n}
\rput(5.709219,-0.29625){(a)}
\usefont{T1}{ptm}{m}{n}
\rput(5.8992186,-3.99625){(b)}
\usefont{T1}{ptm}{m}{n}
\rput(1.6545312,0.68375){$\bu$}
\usefont{T1}{ptm}{m}{n}
\rput(3.3945312,1.58375){$\bv$}
\usefont{T1}{ptm}{m}{n}
\rput(10.894531,-0.01625){$\bv$}
\usefont{T1}{ptm}{m}{n}
\rput(9.6545315,0.74375){$\bu$}
\usefont{T1}{ptm}{m}{n}
\rput(1.6145313,-3.35625){$\bu$}
\usefont{T1}{ptm}{m}{n}
\rput(4.494531,-2.71625){$\bv$}
\usefont{T1}{ptm}{m}{n}
\rput(10.214531,-3.67625){$\bv$}
\usefont{T1}{ptm}{m}{n}
\rput(13.034532,-2.97625){$\bu$}
\psline[linewidth=0.015cm,arrowsize=0.05291667cm 2.0,arrowlength=1.4,arrowinset=0.4]{<-}(9.52,-2.96625)(12.44,-2.92625)
\psline[linewidth=0.015cm,arrowsize=0.05291667cm 2.0,arrowlength=1.4,arrowinset=0.4]{<-}(2.54,-2.68625)(4.14,-2.68625)
\pspolygon[linewidth=0.022](3.2,4.1)(1.05,2)(5.38,2)
\pspolygon[linewidth=0.022](1.6,2.6)(0.2,1.2)(3.0,1.2)
\pspolygon[linewidth=0.022](1.64,-1.37)(0.24,-2.78)(3.0,-2.78) 
\pspolygon[linewidth=0.022](2.42,-1.01)(1.06,-2.39)(3.78,-2.39)   
\pspolygon[linewidth=0.022](9.6,2.6)(8.2,1.23)(11.0,1.23)     
\pspolygon[linewidth=0.022](10.8,3.8)(7.42,0.42)(14.18,0.42)  
\pspolygon[linewidth=0.022](10.02,-1.03)(7.8,-3.2)(12.22,-3.2)
\pspolygon[linewidth=0.022](9.22,-1.41)(7.8,-2.83)(10.64,-2.83)
\end{pspicture} 
}

\caption{For the edges from $\bu$ and $\bv$ to cross, the location of $h(\bv)$ must be on the bold part of the appropriate triangle.}

\label{noncrossfig}
%\end{center}
\end{figure}
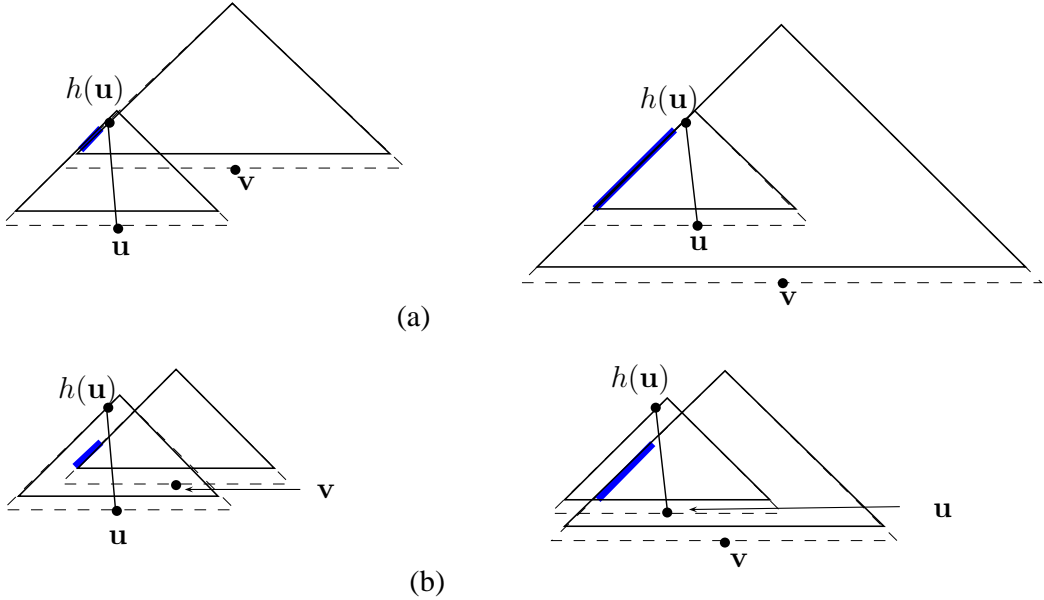

\noindent{\bf Proof}: We present a proof based on Figure \ref{noncrossfig}, a formal, though
cumbersome, proof may be written on these lines.

Let $\bu, \bv \in \Z^2$ be as in Figure \ref{noncrossfig} with $\bu(1)
\leq \bv(1)$ (otherwise we just interchange the roles of $\bu$ and $\bv$). We
show that the edges $\langle \bu, h(\bu)\rangle$ and $\langle \bv, h(\bv)
\rangle$ do not cross.

Let $S^o(\bu, ||\bu - h(\bu)||_1) = \{ \bw \in \Z : ||\bu - \bw||_1 \leq  ||\bu - h(\bu)||_1
\mbox{ and } \bw(2) \geq \bu(2)\}$ be the triangle $S^+(\bu, ||\bu - h(\bu)||_1)$ (defined at the
beginning of Section 2) extended at the base by an extra set of edges on the
horizontal line containing $\bu$. Also let $\bar{S}^o(\bu),\bar{S}^+(\bu)
\subseteq \R^2$ be the smallest simply connected closed triangle in $\R^2$
containing $S^o(\bu, ||\bu - h(\bu)||_1), {S}^+(\bu, ||\bu - h(\bu)||_1)$ respectively.
The three linear segments of the triangle $\bar{S}^o(\bu)$ are appropriately
called the base, the left boundary and the right boundary of it. 
Consider the case when $h(\bu)$ is a vertex on the left boundary of $\bar{S}^+(\bu)$.
A similar argument may be given when $h(\bu)$ is a vertex on the right boundary of $\bar{S}^+(\bu)$.

In case  $\bar{S}^+(\bu) \cap \bar{S}^+(\bv) = \emptyset$, then $\langle \bu,
h(\bu)\rangle$  and $\langle \bv,
h(\bv)\rangle$ do not cross because they lie in  $\bar{S}^o(\bu)$ and
$\bar{S}^o(\bv)$ respectively.

On the other hand if  $\bar{S}^+(\bu) \cap \bar{S}^+(\bv) \neq \emptyset$
then we note that --
\begin{itemize}

\item[(a)] if the left boundary of $\bar{S^o}(\bu)$ has non-empty
intersection with
$\bar{S^o}(\bv)$, then, $h(\bu)$ being open,
the left boundaries of $\bar{S^o}(\bu)$ and $\bar{S^o}(\bv)$ must have an
overlap (see Figure \ref{noncrossfig} (a)). For the edges $\langle \bu,
h(\bu)\rangle$  and $\langle \bv,
h(\bv)\rangle$ to cross, $h(\bv)$ must necessarily be a lattice point, which
is not $h(\bu)$, but lies on the {\it interior}  of this overlap and below
$h(\bu)$. Note that $h(\bv)$ being in the interior of this overlap violates the minimality
of $U_{h(\bu)}$.
\item[(b)] if the left boundary of $\bar{S^o}(\bu)$ has empty intersection
with
$\bar{S^o}(\bv)$ (see Figure \ref{noncrossfig} (b)), then,  for the edges
$\langle \bu,
h(\bu)\rangle$  and $\langle \bv,
h(\bv)\rangle$ to cross, $h(\bv)$ must lie on the part of the left boundary
of $\bar{S^o}(\bv)$ which is in the interior of $\bar{S^o}(\bu)$.
However, this is not possible because the vertices
in the interior of $\bar{S^o}(\bu)$ are closed.
\end{itemize}
This completes the proof of the proposition. \qed

For the proof of Theorem \ref{r-treesforest} in the case $ d = 2$ we consider the process constructed
from the two vertices $ \bu $ and $ \bv $ with $ \bu (2)
= \bv (2) $. Without loss of generality, we may
assume that $ \bu (1) > \bv (1)$.
Since the paths $\{g_n(\bu):\;n\geq 0\}$ and $\{g_n(\bv):\;n\geq 0\}$ do not
cross each other, from Corollary \ref{cor:ZMartingaleD=2} we have that
 $\{ Z_l (\bu, \bv)  = g_{\tau_{l}(\bu,\bv )}(\bu )(1) -
g_{\tau_{l}(\bu,\bv)}(\bv)(1)  : l\geq 0\}$ is a
non-negative martingale.
By the martingale convergence theorem, there exists a random variable
$Z_{\infty}$ such that $ Z_l  (\bu, \bv)
\rightarrow Z_{\infty}$ a.s. as $ l \to \infty$. Also, $0$ being the only
absorbing state of the Markov chain
$ \{ Z_l  (\bu, \bv)  : l \geq 0 \}$ we have
$Z_{\infty}=0$ a.s. and hence $Z_l (\bu, \bv) = 0$ for  some $ l  $ a.s.
This proves Theorem \ref{r-treesforest} for $d=2$.
\qed

%\medspace
%\vspace{.6cm}

\subsection{$d=3$}
We show that Foster's criterion (see Asmussen \cite{A03}, Proposition 5.3 of
Chapter I, page 18), used in Gangopadhyay {\it et al.}\/ \cite{GRS04}, is
applicable here. We start with the process constructed from the vertices $\bu, \bv \in \Z^3$
with $ \bu (3) = \bv (3) $ and consider the process $Z_l = Z_l(\bu, \bv)$
where $Z_l$ is as defined in (\ref{def:Z_l}).
Also, changing the transition probability of $Z_l$ from the state $ {\bf 0} = (0,0) $ in any 
reasonable way, so that the state  $ {\bf 0}$ is no longer absorbing, we  make the Markov chain $ \{ Z_l   ( \bu,
\bv) : l \geq 0 \} $ irreducible. With a slight abuse of notation, we continue to
denote the modified chain by $ \{ Z_l   ( \bu, \bv) : l \geq 0 \} $ and it is enough to
 show that the modified chain is recurrent,

We now show that if the points are far apart, $  Z_l   ( \bu, \bv) $ can be approximated
in expectation by the independent process.
\begin{prop}
\label{prop:ExpecCompareIndAndDep}
For any $\bu, \bv $ as above and $ {\rm x} \in \Z^2$ and $ m \geq 1$, we have
\begin{align*}
\lefteqn{ \biggl| \E\Bigl[  \bigl(  || Z_{l+1} (\bu, \bv) ||_2^2 -
 || {\rm x}||_2^2  \bigr)^m  \mid Z_{l} (\bu, \bv) = {\rm x} \Bigr]  } \\
&  - \E\Bigl[   \bigl(  ||  ( \overline{\bu} +\psi_1^{\bu} )  -
(  \overline{\bv} +\psi_1^{\bv} )  ||_2^2 -  || {\rm x}||_2^2 \bigr)^m \Bigr] \biggr|
\leq C_{17}^{(m)} \exp ( - C_{18}^{(m)} ||  {\rm x} ||_2 )
\end{align*}
where $ C_{17}^{(m)} $ and $ C_{18}^{(m)} $ are positive constants, depending on $m$, and $ ||  \cdot  ||_2 $
denotes the standard Euclidean norm.
\end{prop}

% As in  \cite{GRS04}, to % we show that the estimates (3), (4) and (5)
% %of  \cite{GRS04} hold in this case too, i.e.,
% we need appropriate bounds of  $ \E\Bigl[  \bigl(  || Z_{l+1} (\bu, \bv) ||_2^2 -
%  || Z_{l} (\bu, \bv) ||_2^2 \bigr)^m  \mid Z_{l} (\bu, \bv) = {\rm x} \Bigr] $
%  for $ {\rm x} \in \Z^2 $ and $ m = 1,2, 3  $.
\noindent{\bf Proof}: Since our model is  spatially  translation invariant and
$ Z_{l} (\bu, \bv)  $ is a time homogeneous Markov chain, we may
take  $ \bv = {\bf 0} = (0, 0, 0) $ and $ \bu = ( {\rm x}, 0 ) $ and $ l = 0 $.

Now, we use the coupling described in Subsection \ref{Coupling}, with $ k = 2 $ and $ r = d_{\min}/3 =
( |{\rm x} (1)| + |{\rm x}(2)|)/3$. Observe that $ ||Z_{1} (\bu, \bv) -
{\rm x}||_2 \leq  ||Z_{1} (\bu, \bv) - {\rm x}||_1 \leq  W_1 (\bu,  \bv) $
and $ || \psi_1^{\bu} -  \psi_1^{\bv}  ||_2 \leq || \psi_1^{\bu} -  \psi_1^{\bv}  ||_1 \leq W_1^{(\text{Ind})}
 (\bu, \bv) $, where $ W_1 (\bu,  \bv)$ and $ W_1^{(\text{Ind})}  (\bu, \bv)  $
are as defined in  (\ref{r-W2}) and (\ref{r-Wind}) respectively.
Also, on the event $ (A^{\text{Good}} (r) )^c $, we have  $ W_1 (\bu,  \bv)
> d_{\min}/3 = ( |{\rm x} (1)| + |{\rm x}(2)|)/3  $ and $  W_1^{(\text{Ind})}  (\bu, \bv)   > d_{\min}/3 =
( |{\rm x} (1)| + |{\rm x}(2)|)/3 $. Thus, with $ \psi_1^{\bu}$ and $  \psi_1^{\bv}$ as in (\ref{ZDef}), from the definition
of $ A^{\text{Good}} (r) $ and equation (\ref{eqn:CouplingJtAndInd}), we have
\begin{align*}
& \biggl| \E\Bigl[  \bigl(  || Z_{1} (\bu, \bv) ||_2^2 -
 || {\rm x}||_2^2  \bigr)^m  \mid Z_{0} (\bu, \bv) = {\rm x} \Bigr]   \\
&  - \E\Bigl[   \bigl(  ||  ( \overline{\bu} +\psi_1^{\bu} )  -
(  \overline{\bv} +\psi_1^{\bv} )  ||_2^2 -  || {\rm x}||_2^2 \bigr)^m \Bigr] \biggr| \\
 = & \; \biggl| \E\Bigl[  \bigl(  || Z_{1} (\bu, \bv) ||_2^2 -
 || {\rm x}||_2^2  \bigr)^m  {\mathbf 1} \bigl( (A^{\text{Good}} (r) )^c  \bigr)
 \mid Z_{0} (\bu, \bv) = {\rm x} \Bigr]   \\
&  - \E\Bigl[   \bigl(  ||  ( \overline{\bu} +\psi_1^{\bu} )  -
(  \overline{\bv} +\psi_1^{\bv} )  ||_2^2 -  || {\rm x}||_2^2 \bigr)^m
{\mathbf 1} \bigl( (A^{\text{Good}} (r) )^c  \bigr)   \Bigr] \biggr| \\
\leq & \;  \E \bigl[  2^m  \bigl( || Z_{1} (\bu, \bv) ||_2^{2m} +
|| {\rm x} ||_2^{2m}  \bigr) {\mathbf 1} \bigl( (A^{\text{Good}} (r) )^c  \bigr) \bigr]   \\
&  \qquad +    \E \bigl[ 2^m  \bigl(   ||  ( \psi_1^{\bu}   -  \psi_1^{\bv} ) +  {\rm x}   ||_2^{2m}
+  || {\rm x} ||_2^{2m} \bigr) {\mathbf 1} \bigl( (A^{\text{Good}} (r) )^c  \bigr) \bigr]    \\
\leq & \; 2^{m}\E  \Bigl[  \Bigl(  || {\rm x} ||_2^{2m}
+ 2^{2m}  \bigl[   || Z_{1} (\bu, \bv) -  {\rm x} ||_2^{2m} +  || {\rm x} ||_2^{2m}   \bigr] \Bigr) {\mathbf 1}
\bigl( (A^{\text{Good}} (r) )^c  \bigr) \Bigr] \\
& + 2^{m}\E  \Bigl[  \Bigl(  || {\rm x} ||_2^{2m}
+ 2^{2m}  \bigl[   || \psi_1^{\bu}   -  \psi_1^{\bv} ||_2^{2m} +  || {\rm x} ||_2^{2m}   \bigr] \Bigr) {\mathbf 1}
\bigl( (A^{\text{Good}} (r) )^c  \bigr) \Bigr] \\
\leq & \; 2^{4m} \Bigl[    || {\rm x} ||_2^{2m} \P \bigl( (A^{\text{Good}} (r) )^c  \bigr)
 + \E \bigl[  \bigl( W_1 (\bu,  \bv) \bigr)^{2m}
{\mathbf 1} ( W_1 (\bu,  \bv)  > d_{\min}/ 3 ) \bigr] \\
& \qquad +  \E \bigl[  \bigl(  W_1^{(\text{Ind})}(\bu,  \bv)  \bigr)^{2m} {\mathbf 1}
( W_1^{(\text{Ind})}(\bu,  \bv)   > d_{\min}/ 3 )   \bigr] \Bigr] \\
\leq & \; C_{17}^{(m)} \exp ( - C_{18}^{(m)} ||  {\rm x} ||_2 )
 \end{align*}
for a proper choice of $ C_{17}^{(m)} ,  C_{18}^{(m)} > 0$. \qed

 Using the properties of $   \psi_1^{\bu} $ and $\psi_1^{\bv} $ from Proposition \ref{prop:PropertyIncrementsForIndPaths},
we can compute the moments of $   ||  ( \overline{\bu} + \psi_1^{\bu} )  -  (  \overline{\bv} +\psi_1^{\bv} )  ||_2^2 -
 || {\rm x}||_2^2 $. It is easy to check that
\begin{align}
\label{eqn:FirstMomentInd}
 \E\Bigl[    \bigl(  ||  ( \overline{\bu} +\psi_1^{\bu} )  -
(  \overline{\bv} +\psi_1^{\bv} )  ||_2^2 -  || {\rm x}||_2^2 \bigr)  \Bigr] & = 4 \E \bigl[  \bigl( \psi_1^{\bu}  (1)\bigr)^2 \bigr] = \alpha
~ (\text{say});\\
\label{eqn:SecondMomentInd}
 \E\Bigl[   \bigl(  ||  ( \overline{\bu} +\psi_1^{\bu} )  -
(  \overline{\bv} +\psi_1^{\bv} )  ||_2^2 -  || {\rm x}||_2^2 \bigr)^2  \Bigr] & \geq 8  \E \bigl[  \bigl( \psi_1^{\bu}  (1)\bigr)^2 \bigr] || {\rm x}||_2^2
= 2 \alpha  || {\rm x}||_2^2; \\
\label{eqn:ThirdMomentInd}
 \E\Bigl[  \bigl(  ||  ( \overline{\bu} +\psi_1^{\bu} )  -
(  \overline{\bv} +\psi_1^{\bv} )  ||_2^2 -  || {\rm x}||_2^2 \bigr)^3  \Bigr] & = O ( || {\rm x}||_2^2) .
\end{align}
The proofs follow from straightforward calculations and have been relegated to the Appendix.

Now, we consider $ f : \Z^2 \to [0, \infty) $ defined by $ f ( {\rm x}) = \sqrt{ \log ( 1
+ || {\rm x}||_2^2) } $. Clearly, $  f ( {\rm x})  \to \infty $ as $  || {\rm x}||_2 \to \infty $.
Using Taylor's expansion of the function $ h (t) = \sqrt{ \log (1+t) } $ and observing that the
fourth derivative of $ h$ is always negative, we have
\begin{align*}
\lefteqn{  \E \bigl[  f \bigl( Z_{1} (\bu, \bv) \bigr) -   f  \bigl( Z_{0} (\bu, \bv) \bigr)
 \mid Z_{0} (\bu, \bv) =  {\rm x} \bigr] } \\
&  = \E \bigl[  h ( || \bigl( Z_{1} (\bu, \bv) \bigr) ||^2_2 ) - h ( ||  Z_{0} (\bu, \bv) ||^2_2) \mid Z_{0} (\bu, \bv) =  {\rm x} \bigr]  \\
& \leq \sum_{ m= 1}^3 \frac{ h^{(m)} (|| {\rm x} ||_2^2  ) }{ m !}  \E\Bigl[  \bigl(  || Z_{1} (\bu, \bv) ||_2^2 -
 || {\rm x}||_2^2  \bigr)^m  \mid Z_{0} (\bu, \bv) = {\rm x} \Bigr] \\
& \leq \sum_{ m= 1}^3 \frac{ h^{(m)} (|| {\rm x} ||_2^2  ) }{ m !} \E\Bigl[   \bigl(  ||  ( \overline{\bu} +\psi_1^{\bu} )  -
(  \overline{\bv} +\psi_1^{\bv} )  ||_2^2 -  || {\rm x}||_2^2 \bigr)^m \Bigr]  \\
& \quad + \sum_{ m= 1}^3 \frac{ \bigl| h^{(m)} (|| {\rm x} ||_2^2  )  \bigr|}{m!} C_{17}^{(m)} \exp ( - C_{18}^{(m)} ||  {\rm x} ||_2 )
\end{align*}
where $ h^{(m)} $ represents the $m$ th derivative of $ h $. Plugging in the expressions for $ h^{(m)} $ and
 the moments given in (\ref{eqn:FirstMomentInd}), (\ref{eqn:SecondMomentInd}) and
(\ref{eqn:ThirdMomentInd}), we have that the first sum above is bounded by
$  - \alpha  || {\rm x}||_2^2 / \bigl[  8 ( 1 +  || {\rm x}||_2^2)^2 \bigl( \log ( 1 +  || {\rm x}||_2^2) \bigr)^{3/2} \bigr]$
for all large $  || {\rm x}||_2 $ and the second sum is bounded by $ C_{19} \exp ( - C_{20} || {\rm x} ||_2) $ 
for a proper choice of the constants $ C_{19} $ and $C_{20}$.  % in conjunction with Proposition \ref{prop:ExpecCompareIndAndDep},
This yields that
\begin{equation*}
 \E \bigl[  f \bigl( Z_{l+1} (\bu, \bv) \bigr) -   f  \bigl( Z_{l} (\bu, \bv) \bigr)  \mid Z_{l} (\bu, \bv) =  {\rm x} \bigr] < 0
\end{equation*}
for $ ||  {\rm x}||_2 $ large enough. This implies that $ Z_l  (\bu, \bv)  $ is recurrent and completes the proof for $ d = 3$.  \qed

% Now, the estimates (3), (4) and (5) of \cite{GRS04} follow from direct
% computations of the moments of the marginals $ \psi_1^{(1)}$ and $ \psi_1^{\bv}$.  For example,
% when $ m = 1$, with $ \psi_1^{(i)} (j) $ denoting the
% the $ j $ th co-ordinate of $ \psi_1^{(i)} $ for $ i = 1, 2 $ and $ j = 1, 2$, using
% the observations made about the marginals  in Section 3, we have $
% \E\Bigl[    ||  ( \overline{\bu} + \psi_1^{\bu} )  -  (  \overline{\bv} +\psi_1^{\bv} )  ||_2^2 -
%  || {\rm x}||_2^2  \Bigr] = \E\Bigl[  \bigl( {\rm x} (1) + \psi_1^{\bu} (1) - \psi_1^{\bv}(1) \bigr)^2 +
%    \bigl( {\rm x} (2) + \psi_1^{\bu} (2) - \psi_1^{\bv}(2) \bigr)^2 -   \bigl( {\rm x} (1) \bigr)^2
%    -  \bigl( {\rm x} (2) \bigr)^2  \Bigr] = 4 \E \bigl[  \bigl( \psi_1^{\bu}  (1)\bigr)^2 \bigr]  $ which
%    yields (3). Similar calculations yield (4) and (5).

\subsection{$d\geq 4$}
We present the proof for $d=4$; the argument being similar for $d>4$. To show that $ \P(G $ has infinitely many distinct trees$)=1$,
it is enough to prove that
$ \P (G $ has at least $m$ distinct trees$) = 1 $ for any $ m \geq 2$. The probability
measure $ \P $ is ergodic as it is a product measure given by i.i.d. uniform random variables
on each vertex of $ \Z^d$. Clearly, for any $ m \geq 2$, the event $ \{ G $ has at least $m$ distinct trees$\}$ is
translation invariant under the group of translations and hence $ \P(G $ has at least
$m$ distinct trees$) $ is either $0$ or $1$. So, it suffices to show that $\P(G $ has at least $m$ distinct trees$)>0$.

We first show the above for $m=2$. It is enough to exhibit two open vertices %, $ \bu$ and $ \bv$,
such that the paths from those two vertices do not meet with positive probability.
We follow the same ideas as in  \cite{GRS04} to achieve this; however there is one crucial change.
In \cite{GRS04}, each unit increase in the fourth axis represents an unit increase in time co-ordinate. For our model,
the time taken for the joint regeneration of paths starting from two vertices is taken to be a unit of time. More precisely, starting with two open vertices $ \bu$ and $ \bv$ having the same fourth co-ordinate,
at the first joint regeneration time of the paths from $ \bu$ and $ \bv$, we think of the time  co-ordinate as having increased by
a unit and  at each joint regeneration thereafter, the time co-ordinate increases by  one unit. At these joint regenerations,
the fourth co-ordinates for both paths are equal. At the $l$ th regeneration, the paths have not yet met if $  Z_l (\bu,\bv) \neq (0,0,0)$.
Since the paths coalesce once they meet, it is enough to prove that  with positive probability they do not meet at any of the
joint regeneration times, i.e., $ \P \bigl( Z_l (\bu,\bv) \neq (0,0,0) $ for all $ l \geq 0\bigr) > 0$
for some pair of open vertices $ \bu $ and $ \bv $. %, i.e.,
%the Markov chain $\{ Z_l (\bu,\bv) : l \geq 0 \}$ is not absorbed at $(0,0,0)$ has a positive probability.

Our strategy is to let  the  joint paths, from $ \bu $ and $\bv $, evolve for $ n^4 $ joint regeneration times
where $  \bu $ and $\bv $ are sufficiently far apart (of order $ n $).
Then, with a very high probability the  paths  have travelled further away (of  order $n^2$).
Using the Markov property at the regeneration times, we may now start the paths from these new
vertices and continue this process. We make this more precise in the following Proposition.
For $\epsilon>0$, define the event
 \begin{align*}
  A_{n, \epsilon} (\bu, \bv )  :=
 \Bigl\{  & Z_{n^4} ( \bu,\bv)   \in  D_{n^{2(1+\epsilon)}}
 \setminus D_{n^{2(1-\epsilon)}}\Bigr\},
%%% &  Z_{j} ( \bu,\bv) \neq (0,0,0)  \text{ for all }j=1,\dotsc,n^4 \Bigr\},
 \end{align*}
  where $ D_r   := \{ \bx \in \Z^3: || \bx||_1\leq r\}$.
We show
\begin{prop}
\label{prop:D4MainProp}
 For $0<\epsilon<\frac{1}{3}$, there exist constants $C_{21},\beta>0$ and $n_0\geq 1$ such that,
 for all $n\geq n_0$,
 \begin{equation*}
 \inf_{ \overline{\bv }   \in \overline{ \bu }  +D_{n^{1+\epsilon}}
 \setminus D_{n^{1-\epsilon}}  }
 \P \bigl( A_{n,\epsilon}(\bu,\bv)  \mid \bu, \bv \in V \bigr)\geq 1 - C_{21} n^{-\beta}.
  \end{equation*}
\end{prop}

Assuming Proposition \ref{prop:D4MainProp}, we first prove the result.
Fix $ \epsilon < 1/3 $ and  choose $ \bu^0 = (0,0,0,0)$ and $ \bu^1 = (n_0,0,0,0) $ where $ n_0 $
 is as above. Let $ E_2 $ be the event that
both $ \bu^0, \bu^1 \in V $, so that $ \P (E_2) = p^2 > 0$.
Clearly,  $ n_0^{ 1-\epsilon} < ||  \bu^0 -  \bu^1 ||_1  < n_0^{1+\epsilon} $.
We show that $ \P ( Z_l ( \bu^0, \bu^1) \neq (0,0,0) $ for all $ l \geq 1 | E_2 ) > 0 $.
For $ j \geq 1$, set $ r_j = \sum_{ i = 0}^{j-1} \bigl( n_0^{2^{i}}\bigr)^{4} $.
Since $ (0,0,0) $ is an absorbing state, we have $  \P \bigl(  Z_l (\bu^0,\bu^1) \neq (0,0,0) 
\text{ for all } l \geq 1 \mid E_2  \bigr)
=   \lim_{ j \to \infty} \P \bigl(  Z_{r_j} (\bu^0,\bu^1) \neq (0,0,0) \mid E_2  \bigr)  
\geq \lim_{ j \to \infty} \P \bigl(  Z_{r_i} (\bu^0,\bu^1)
 \in   D_{n_0^{(2^i)(1+\epsilon)}}
\setminus D_{n_0^{(2^i)(1-\epsilon)}} \text{ for all } 1 \leq i \leq j \mid E_2 \bigr) $ -- 
the last inequality follows because
$Z_{r_i} (\bu^0,\bu^1)
 \in   D_{n_0^{(2^i)(1+\epsilon)}}
\setminus D_{n_0^{(2^i)(1-\epsilon)}} $
imposes further restrictions  on the Markov chain.
 For any $ j \geq 1$, we have
\begin{align*}
&  \P \bigl(  Z_{r_i} (\bu^0,\bu^1) \in   D_{n_0^{(2^i)(1+\epsilon)}}
\setminus D_{n_0^{(2^i)(1-\epsilon)}} \text{ for all } 1 \leq i \leq j \mid E_2 \bigr)  \\
& =  \P \bigl(  Z_{r_j} (\bu^0,\bu^1) \in  D_{n_0^{(2^j)(1+\epsilon)}}
\setminus D_{n_0^{(2^j)(1-\epsilon)}}  \mid Z_{r_{j-1}} (\bu^0,\bu^1) \in  D_{n_0^{(2^{j-1})(1+\epsilon)}}
\setminus D_{n_0^{(2^{j-1})(1-\epsilon)}}\bigr) \\
& \quad   \times  \P \bigl(  Z_{r_j} (\bu^0,\bu^1) \in  D_{n_0^{(2^i)(1+\epsilon)}}
\setminus D_{n_0^{(2^i)(1-\epsilon)}} \text{ for all } 1 \leq i \leq j-1 \mid E_2  \bigr)
\end{align*}
%where we have used the Markov property in writing the conditional probability.
Now, using the Markov property of
$ Z_{l} (\bu^0,\bu^1) $,  Proposition \ref{prop:D4MainProp}
and the translation invariance of our model, we see that the conditional probability in the second line is at least as large as
$ \inf _{   \overline{\bv }   \in \overline{ \bu }  + D_{n_0^{(2^{j-1})(1+\epsilon)}}
\setminus D_{n_0^{(2^{j-1})(1-\epsilon)}}  }
 \P \bigl( A_{n_0^{2^{j-1}} ,\epsilon}(\bu,\bv)   \mid \bu, \bv \in V  \bigr) \geq  1 - C_{21} (n_0^{2^{j-1}})^{-\beta}
$. By repeating this argument, we conclude that 
$$  \P \bigl(  Z_{r_j} (\bu^0,\bu^1) \neq (0,0,0) \bigr)  \geq
\prod_{i=1}^{j} \bigl[  1 -  C_{21} (n_0^{2^{i-1}})^{-\beta} \bigr]  \to \prod_{i=1}^{\infty} \bigl[
1 -  C_{21} (n_0^{2^{i-1}})^{-\beta} \bigr].  $$
Therefore, $ \P (G $ has at least $2$ distinct trees$) \geq p^2  \prod_{i=1}^{\infty} \bigl[
1 -  C_{21} (n_0^{2^{i-1}})^{-\beta} \bigr] >0 $.

The above calculations hold for any pair of points which satisfy the initial condition.
We now use this to prove that $ \P (G $ has at least $m$ distinct trees$) > 0 $ for any $ m \geq 2$.
Fix $ \delta > 0 $ such that $ m(m-1)\delta / 2 < 1 $.  Note that  $ \prod_{i=1}^{\infty} \bigl[1 -
C_{21} (n^{2^{i-1}})^{-\beta} \bigr]  \to 1 $ as $ n \to \infty$. Now, we choose $ n_1 > n_0 $ so large
that $  \prod_{i=1}^{\infty} \bigl[1 -  C_{21} (n_1^{2^{i-1}})^{-\beta} \bigr]  > 1 - \delta $
and $ m < n_1^{\epsilon} $ where $ \epsilon $ and $ n_0 $ are as above. Now, consider the points
$ \bu^{i} = ((i-1) n_1, 0, 0, 0) $ for $ i = 1,2, \dotsc, m $. Clearly, all of them have the same fourth
co-ordinate and $ n_1^{1-\epsilon} < n_1 \leq  || \bu^{i} - \bu^{j} ||_1 \leq mn_1 < n_1^{1+\epsilon}$.
Let $ E_m $ be the event that all the points $ \bu^{1}, \dotsc, \bu^{m} $ are open. Then $ \P ( E_m ) = p^m > 0$.
We now consider the event $ A_{i,j} $ that the paths
from $ \bu^i $ and $ \bu^{j} $ do not  meet for $ i > j $. From above calculations and by our choice of $ n_1 $,
 we have $ \P ( A_{i,j} | E_m ) > 1 - \delta $. Further, we consider the intersection of all the events
$ A_{ i, j} $ for $ 1 \leq j < i \leq m $. Clearly, $ \P ( \cap_{ 1 \leq j < i \leq m } A_{i,j}  \mid E_m ) \geq 1 - m(m-1)\delta / 2
$ so that $ \P (G $ has at least $m$ distinct trees$) \geq p^m (  1 - m(m-1)\delta / 2 ) > 0 $.

To prove Proposition \ref{prop:D4MainProp}, we define a new event where paths are constructed by using
independent uniform random variables of their own and then use the coupling described in Subsection \ref{Coupling}.
Consider the event
 \begin{align*}
  A^{(\text{Ind})}_{n, \epsilon} (\bu, \bv )  :=
 \Bigl\{  & \overline{\bv} + \sum_{l = 1}^{n^4} \psi_l^{\bv}   \in  \overline{\bu} +
 \sum_{l = 1}^{n^4} \psi_l^{\bu} + D_{n^{2(1+\epsilon)}}
 \setminus D_{n^{2(1-\epsilon)}} , \\
 &   \overline{\bv} +  \sum_{l = 1}^{j} \psi_l^{\bv}  \not \in   \overline{\bu} +
 \sum_{l = 1}^{j} \psi_l^{\bu} +
 D_{ K \log n }   \text{ for all } j=1,\dotsc,n^4 \Bigr \},
 \end{align*}
 where $ K $ is a suitably chosen large constant. We have
% This event corresponds to the event (19)
% defined in  \cite{GRS04} for which it was shown that
\begin{prop}
\label{prop:D4IndVersionMain}
 For $0<\epsilon<\frac{1}{3}$,  there exists $ n_0  $ such that
\begin{equation*}
 \inf_{ \overline{\bv }   \in \overline{\bu }  +D_{n^{1+\epsilon}}
 \setminus D_{n^{1-\epsilon}}  }
 \P \bigl( A^{(\rm{Ind})}_{n,\epsilon}(\bu,\bv) \bigr) \geq 1 - C_{22} n^{-\alpha},
  \end{equation*}
for some constant $ C_{22}, \alpha > 0 $ and for all $ n \geq n_0$.
\end{prop}

First we prove Proposition \ref{prop:D4MainProp} assuming the result of Proposition \ref{prop:D4IndVersionMain}.

\noindent{\bf Proof of Proposition \ref{prop:D4MainProp}}: We employ the coupling described earlier in
Subsection \ref{Coupling}, on the event $ A^{(\text{Ind})}_{n, \epsilon} (\bu, \bv )  $ defined
above. This time we will continue the coupling step by step for $ n^4 $ simultaneous regeneration steps of independent paths.
At each step we choose $ r = K \log n / 3 $
and say that the coupling is  successful at step $ j $ if  the event $ A^{\text{Good}} (r)  $
occurs. We do the coupling at step $ j+1 $ if the coupling is successful at step $ j$. Note,
if the coupling is successful at every step $ j = 1, \dotsc, n^4 $, we have,  for $ j = 1, 2, \dotsc, n^4$,
\begin{equation*}
  \overline{\bu^0} +  \sum_{l = 1}^{j} \psi_l^{\bu^0} = \overline{g}_{\tau_j(\bu^0,\bv^0)} ( \bu^0)\text{ and }
  \overline{\bv^0} +   \sum_{l = 1}^{j} \psi_l^{\bv^0} = \overline{g}_{\tau_j(\bu^0,\bv^0)} ( \bv^0).
\end{equation*}
Therefore, we get
\begin{equation*}
 \P \bigl( A_{n,\epsilon}(\bu^0,\bv^0) \bigr)  \geq  \P \bigl(
A^{(\text{Ind})}_{n, \epsilon} (\bu^0, \bv^0 )  \cap \{\text{Coupling is successful for } j = 1, 2,
\dotsc, n^4\} \bigr).
\end{equation*}
Using the Markov property and the estimate of the coupling being successful, given in
(\ref{eqn:CouplingJtAndInd}), we obtain, for all sufficiently large $ n $,
 \begin{equation*}
 \P \bigl( A_{n,\epsilon}(\bu^0,\bv^0) \bigr)  \geq 1 -  C_{22} n^{-\alpha} -
 C_{15} n^4 \exp ( - C_{16} K \log n / 3   )\geq 1 - C_{21} n^{- \beta}
\end{equation*}
for suitable choice of $ \beta > 0 $ and $ C_{21}$. This proves the proposition \ref{prop:D4MainProp}. \qed

Finally, we indicate the steps for proving Proposition \ref{prop:D4IndVersionMain}. 
By Remark \ref{rem:DiffsMakeRW}, $  \{ S_j = \overline{\bv} - \overline{\bu} + \sum_{ l=1}^{j}
\psi_l^{\bv} - \psi_l^{\bu} : j \geq 1 \}$ is an aperiodic, isotropic, symmetric random walk on $ \Z^3$ 
starting from $  \overline{\bv}- \overline{\bu} $.
The event $ \P \bigl( A_{n,\epsilon}(\bu^0,\bv^0) \bigr) $ is not
satisfied if any of the following occurs : (a)  the random walk travels too far, i.e.,  $ \{ S_{n^4} \not \in D_{n^{2(1+\epsilon)}} \}$
 or the random walk travels too little, i.e., $  \{ S_{n^4} \in   D_{n^{2(1-\epsilon)}} $) or (c) it comes too close to a given point
at distance of order $ n $, i.e., $ \{ S_j \in - (\overline{\bv}- \overline{\bu}) + D_{K \log n}, $ for some $ 1 \leq j \leq n^4 \}$.
For an aperiodic, isotropic, symmetric random walk, it can be shown that each of these events have small probability.
For more details, we refer the reader to Lemma 3.3 of \cite{GRS04}. \qed
%\input{4ASRR1807.tex}
%\input{sec5rf.tex}
%\input{sec5rf.AS.RR.2014.07.08.tex}
% !TEX root = ./DSFrf.tex
% !TEX program = XeLaTeX
\section{Brownian Web}
In this section we prove Theorem \ref{r-BW}.
% We mostly follow techniques used
% by Ferrari {\it et al.}\/ \cite{FFW05}, but one should note that in case of
% DSF we have long range interactions and that involves nontrivial
% modifications
% of the argument.
We begin by recalling that the Brownian web takes values in the metric
space ${\cal H}$ equipped with the Hausdorff metric
$d_{\cal H}$ where ${\cal H}$ is the space of compact
subsets of the path space $(\Pi, d_{\Pi})$ (see the discussion in the paragraphs
after the statement of Theorem \ref{r-treesforest} in Section 1).
As introduced earlier, for any $n \geq 1$,
the collection of scaled paths ${\cal X}_n(\gamma,\sigma)$ is obtained from $G$
with normalization constants
$\gamma,\sigma$ and we had remarked that
the closure of ${\cal X}_n(\gamma,\sigma)$ in $(\Pi, d_{\Pi})$ denoted by
$\overline{{\cal X}}_n(\gamma,\sigma)$ is a $({\cal H}, {\cal B}_{\cal H})$
valued random variable.
% Note that the closure of ${\cal X}_1$ is obtained from
% ${\cal X}$ by adding all the paths of the form $\pi\equiv -\infty$  or
% $\pi\equiv \infty$ with starting time $\sigma_{\pi} \in \Z\cup
% \{-\infty,\infty\}$.
% From now onwards we denote $\overline{{\cal X}}_n(\gamma,\sigma)$ also by
% ${\cal X}_n(\gamma,\sigma)$.

We need some more notation.
For a subset $\Gamma \subseteq \Pi$ of paths  and for $t\in \R$ let
$\Gamma^{t} := \{\pi\in\Gamma  : \sigma_{\pi} \leq t\}$ be the set of paths which start
`below' $t$.
For $t>0$ and $t_0,a,b\in \R$ with $a<b$, we define two counting random variables as
follows
\begin{align*}
\eta_\Gamma(t_0,t;a,b) & := \#\{\pi(t_0+t):\pi \in \Gamma^{t_0}\text{ and }
\pi(t_0)\in [a,b]\} \text{ and }\\
\hat{\eta}_\Gamma(t_0,t;a,b) & := \#\{\pi(t_0+t):\pi \in \Gamma^{t_0}\text{ and }
\pi(t_0 + t)\in [a,b]\}.
\end{align*}
%
% $$
%
% be the
% $\{0,1,2,\dots,\infty\}$  valued random variable giving the number of distinct
% points on $\R \times \{t_0 + t\}$ intersected by some path in $K^{t_0^{-}}$.
% Define another counting random variable
% \begin{align*}
% \hat{\eta}_K(t_0,t;a,b):=& \#\{\pi(t_0+t):\pi \in K^{t_0^{-}}\text{ and }
% \pi(t_0 + t)\in [a,b]\}.
% \end{align*}
Theorem 2.2 in Fontes {\it et al.} \cite{FINR04} provided a  criteria
for a sequence of $({\cal H}, {\cal B}_{\cal H})$ valued random variables
with non-crossing paths to converge weakly to the Brownian web.
In the following we denote, the standard Brownian motion starting
from $\bx$ by $B^{\bx}$ and standard  coalescing Brownian
motions starting from $\bx^1,
\dotsc, \bx^k$ respectively, by $ (W^{\bx^1}, \dotsc, W^{\bx^k})$.

\begin{theorem}
\cite{FINR04} Suppose $\xi_1,\xi_2,\dotsc$ are $({\cal H},B_{{\cal H}})$
valued random variables with non-crossing paths. Assume that the following
conditions hold.
\begin{itemize}
\item[$(I_1)$] For all $\by \in \R^2$, there exist $\zeta_n^{\by} \in \xi_n$ such that 
for any finite set of points $\bx^1, \dotsc, \bx^k$ from a deterministic countable dense set ${\cal D}$ of $\R^2$, 
as $n \to \infty$,
$
( \zeta^{\bx^1}_n, \dotsc, \zeta^{\bx^k}_n) \Rightarrow (W^{\bx^1}, \dotsc, W^{\bx^k} ).
$
% where $ W^{\bx^1}, \dotsc, W^{\bx^k}$ are standard coalescing
% Brownian motions starting from $\bx^1,\dotsc, \bx^k$ respectively.

\item[$(B_1)$] For all $t>0,\limsup_{n\rightarrow \infty}\sup_{(a,t_0)\in
\R^2}\P(\eta_{\xi_n}(t_0,t;a,a+\epsilon)\geq 2)\rightarrow 0$  as $
\epsilon\downarrow 0$.
\item[$(B_2)$] For all $t>0,\frac{1}{\epsilon}\limsup_{n\rightarrow
\infty}\sup_{(a,t_0)\in
\R^2}\P(\eta_{\xi_n}(t_0,t;a,a+\epsilon)\geq 3)\rightarrow 0$  as $
\epsilon\downarrow 0$.
\end{itemize}
Then $\xi_n$ converges in distribution to the standard Brownian web ${\cal W}$.
\end{theorem}
The convergence in $ (I_1) $ occurs in the space $ \Pi^k $. Note that the
convergence in $ \Pi $ implies that the starting points converge as points
in $ \R^2 $ and the paths converge uniformly on the compact sets
of time.

In Theorem 1.4 and Lemma 6.1 of  Newman {\it et al.} \/ \cite{NRS05}, it was further proved that the
condition $(B_2)$ can be replaced by $(E_1^{\prime})$ where

\begin{itemize}
\item[$(E_1^{\prime})$] if ${\cal Z}^{t_0}$ is any subsequential limit of
$\{{\cal X}_n^{t_0}: n \geq 1\}$ for  $ t_0 \in \R$, then
for all $t,a,b\in \R$ with $t > 0$ and $a<b$,
$\E[\hat{\eta}_{{\cal Z}^{t_0}}(t_0,t;a,b)]\leq
\E[\hat{\eta}_{{\cal W}}(t_0,t;a,b)]=\frac{b-a}{\sqrt{t\pi}}$.
\end{itemize}

It is worthwhile mentioning here that for a sequence of $({\cal H}, {\cal
B}_{{\cal H}})$ valued random variables $\xi_n$ with non-crossing paths, property $(I_1)$ implies
tightness (see Proposition B.2 in the Appendix of \cite{FINR04}) and hence
such a subsequential limit ${\cal Z}^{t_0}$ exists. Thus, to prove Theorem  \ref{r-BW}
we need to show that for some $\gamma(p)>0$
and $\sigma(p)>0$ the sequence
$\overline{{\cal X}}_n(\gamma,\sigma)$ satisfies the conditions $(I_1)$, $(B_1)$ and
$(E_1^{\prime})$ and
hence converges to the standard Brownian web.

\subsection{Verification of condition $(I_1)$}

We proceed by the method of induction and follow a mixture of argument of Ferrari {\it et al.}\/ \cite{FFW05}
and Coletti {\it et al.}\/ \cite{CFD09}.
In Section 2, we proved a regeneration
property for any single path, which we use in Proposition \ref{lemma:Bmotion} to
show the convergence of this path (appropriately scaled) to a Brownian
motion. For showing the joint convergence of more than one path, we use the fact that the paths behave
(almost) independently when they are separated by a large distance and, when they come close to each other,
they coalesce very quickly. This idea was initially introduced by Ferrari {\it et al.}\/ \cite{FFW05}.
It should be noted here that the dependency structure of our model is quite different from that of \cite{FFW05}
where paths are independent when they are separated by a fixed distance. Later Coletti {\it et al.}\/ \cite{CFD09} modified it
to deal with long range interactions and we use a similar approach to prove the joint convergence of paths.

We first recall that for a path $ \pi^{\bu}$ and $ \gamma , \sigma > 0$, the scaled path is
defined by  $ \pi_n^{\bu} = \pi_n^{\bu} (\gamma,\sigma) : [ \bu(2) /n^2\gamma, \infty] \to
[-\infty, \infty]$ such that  $\pi_n^{\bu}  (t)=
\pi(n^2\gamma t)/n \sigma$. Note that the distribution of the path $ \pi^{\bu} $
depends only on uniform random variables in $ \{ y > \bu(2) \}$ and is
independent of the open/closed
status of ${\bu}$.
We first show that the scaled path starting from the origin converges to the standard Brownian motion.
\begin{prop}
\label{lemma:Bmotion}
%Assume that ${\bf 0} \in V.$
There exist $\gamma:= \gamma(p)$ and $\sigma:= \sigma(p)$ such that
as $n \to \infty$,
\begin{equation*}
  \pi^{\mathbf 0}_n  \Rightarrow  B^{\mathbf 0} \text{ in } (\Pi, d_{\Pi}).
\end{equation*}
% where $ W^{\bx^1}, \dotsc, W^{\bx^k}$ are standard
% independent coalescing Brownian motions starting from $\bx^1,\dotsc,
% \bx^k$ respectively.
\end{prop}

\noindent\textbf{Proof}: We use the fact that the path in between the regeneration steps
can be broken up into i.i.d. pieces. For a path $\pi$ we define the
modified path $\tilde{\pi}$  which is linear between
successive regeneration points of $\pi$.  Using Proposition
\ref{prop:WidthTail} we have that the displacements  between
successive regeneration times are independent and have exponential moments --
this allows an application of Donsker's
invariance principle to the modified path to prove the convergence to the Brownian motion.
%Again inter-regeneration time differences have exponential
%moments, which is used to prove the

Let $\tau_j$ and $T_j$ denote the $j$ th regeneration step and time respectively for the
path starting from $\mathbf{0}$ (see (\ref{r-tau}) and (\ref{def:TimekthReg})).
Remembering that $ g_j (\mathbf{0} )  $ is the position of the path starting from $ \mathbf{0}  $ after the $j$ th step,
let $ Y_j = Y_{j}^{({\bf 0})}  = \overline{g}_{\tau_j} ({\bf 0})
- \overline{g}_{\tau_{j-1}} ({\bf 0}) =  g_{\tau_j} ({\bf 0})(1) - g_{\tau_{j-1}} ({\bf 0})(1)$ (see (\ref{Ydef})). We define a piecewise linear path
$ \tilde{\pi} $
as follows: for $ T_j \leq t < T_{j+1}, j \geq 0 $,
\begin{align*}
  \tilde{\pi} (t) := &~ g_{\tau_j}  ({\bf 0})(1)  + \frac{ t- T_j}{ T_{j+1} - T_j}
(g_{\tau_{j+1}}  ({\bf 0}) (1)- g_{\tau_j}  ({\bf 0}) (1) )
\end{align*}
and its diffusively scaled version $ \tilde{\pi}_n $ by
\begin{align*}
\tilde{\pi}_n(t) = & ~\tilde{\pi}_n(\gamma,\sigma)(t) :=
\frac{1}{n\sigma} \tilde{\pi}( n^2\gamma t)
\end{align*}
 for $ t  \geq 0$.
% which interpolates
% linearly between $  g_{\tau_n}  ({\bf 0}) = (\sum_{j=1}^n Y_j , T_n)   $ and $
% g_{\tau_{n+1}}  ({\bf 0})  = (\sum_{j=1}^{n+1} Y_j , T_{n+1} ) $ and call this path $
% \tilde{\pi}^{1,(\text{Ind})} $. The , i.e.,
% \begin{equation*}
%
% \end{equation*}
Next we define another stochastic process,  $S$ on $[0,\infty)$  as follows: for $ j \leq t < j+1, j \geq 0$,
\begin{equation*}
S(t) = T_{j} + (t-j) (T_{j+1} - T_{j}) .
\end{equation*}
Clearly, $ S(t) $ is a strictly increasing process. Hence, $ t \mapsto S(t) $
admits an inverse $ S^{-1} (t) $ which is
also strictly increasing.
The process $ S(t) $ denotes the time change required to track
the path $ \tilde{\pi} $.
More precisely, we have,   for $ t \geq 0$,
\begin{equation*}
\tilde{\pi}_n(t) = X_{n} \bigl( S^{-1}(n^2\gamma t) / n^2 \bigr) .
\end{equation*}
where the process
$ X_n = X_n(\gamma,\sigma)$ on $[0,\infty)$ is defined
as follows: $ X_n (0) = 0 $ and  for $ t > 0 $,
\begin{equation*}
X_{n}(t) := \frac{1}{n\sigma} \Bigl[  (n^2 \gamma t - \lfloor
n^2 \gamma t \rfloor ) Y_{ \lfloor n^2 \gamma t \rfloor +1} +  \sum_{i=1}^{ \lfloor n^2 \gamma t \rfloor } Y_{i} \Bigr] .
\end{equation*}
% Fix $1\leq i\leq k$ and as earlier take $\bx_i = \mathbf{0}$.
% For $n\geq 1$ we define the random process $X^{i}_n = X^{i}_n(\gamma,\sigma)$
% as follows
% $$
% X^{i}_{n}(t) := \frac{Z^{(i)}_{[n^2 t]} + (n^2 t -[n^2 t])(Z^{(i)}_{[n^2 t]+1} - Z^{(i)}_{[n^2 t]})}{n\sigma}
% $$
% where $Z^{(i)}_j = \overline{g}^{i}_{\tau^{\text{(Ind)}}_j(k)}(\mathbf{0})$ as defined in (\ref{defn:Z_lind}).
From Remark \ref{rem:FromSinglePoint}, $ Y_i$'s are symmetric and  i.i.d., so that $ \E ( Y_1 ) = 0 $. Thus,
from Donsker's invariance principle,  it follows
that, for 
$$ \sigma = \sigma_0 := \sqrt{\text{Var}(Y_1)}, $$  the process $X_{n}$
converges weakly to the standard Brownian motion starting from $0$.
% \ifShowMarginPar
% \marginpar{$\sigma_0$ should include the time change}
% \fi

%Assuming that we are working on a probability space such that this
%convergence happens almost surely in $(\Pi, d)$, we define
%% where $T_{j}^{(\text{Ind})}(k)$ is as defined in (\ref{r-jtindregtime}).
%Observe that
%
%{\color{red} what is $ \tilde{\pi}^{1}_n(t)  $? This proposition needs correction. Why should
%time for $ X_n$ start from $0$? We have assumed only that $ \bx^1_n (2) / n^2 \gamma \to \bx^1 (2)$. }

Let $ N(t) $ be the number of the renewals for the process $ S(t) $ up to time $t$, i.e.,
$ N(t) = \lfloor S^{-1} (t) \rfloor $ so that, $ N(t) \leq S^{-1}(t) \leq N(t) + 1$.
By the renewal theorem (see
Theorem 4.4.1 of Durrett \cite{D10}),
$ S^{-1}(n^2\gamma t)/ n^2  \to g(t) :=\frac{\gamma t}{\E(T_1)},$ $t \geq 0$
almost surely.
Taking 
$$\gamma = \gamma_0 := \E[T_{1}],$$
we conclude that
\begin{equation*}
\tilde{\pi}_n  \Rightarrow B^{\mathbf{0}}.
\end{equation*}

Finally to conclude the result, it is enough to show that, for any $ s > 0 $ and $ \epsilon > 0 $
\begin{align*}
&~ \P \bigl[  \sup \{ | \tilde{\pi}_n (t) - \pi_n^{\mathbf{0}} (t) | :  t \in [0,s] \} > \epsilon \bigr] \\
= & ~  \P \bigl[  \sup \{ | \tilde{\pi} (t) - \pi^{\mathbf{0}} (t) | :  t \in [0, n^2 \gamma s] \} > n \sigma \epsilon \bigr] \to 0
\end{align*}
% $$
% \P \bigl[  \sup \{ | \tilde{\pi}^{1}_n (t) - \pi^{1}_n (t) | :  t \in [0,s] \}
% > \epsilon \bigr]
% =  \P \bigl[  \sup \{ | \tilde{\pi}^{1} (t) - \pi^{1} (t) | :
%  t \in [0, s n^2] \} > \epsilon n  \bigr]
% \to 0
% $$
as $ n \to \infty $.
From the definition of $ W_j ({\mathbf 0}) $ (see  (\ref{r-W2})), for any $ t \in [T_j, T_{j+1}]  $ we have
$  |  \tilde{\pi}  (t) - \pi^{\mathbf{0}}  (t)  | \leq 2W_{j+1} $ for all $j \geq 0$.
Since $ N (n^2 \gamma s) \leq \lfloor n^2 \gamma s \rfloor$, we have
\begin{align*}
&~  \P \bigl[  \sup \{ | \tilde{\pi} (t) - \pi^{\mathbf{0}} (t) | : t \in [0, n^2 \gamma s] \} > n \sigma \epsilon \bigr]  \\
\leq &~ \P \bigl[ 2 \max \{ W_j  ({\mathbf 0})  : j = 1,\dotsc,  \lfloor n^2 \gamma s \rfloor \}  >  n \sigma \epsilon\bigr] \\
\leq &~ \lfloor n^2 \gamma s \rfloor \P (  2 W_1  ({\mathbf 0})  > n \sigma \epsilon ) \to 0 \text{ as } n \to \infty,
\end{align*}
where the last step follows from Proposition \ref{prop:WidthTail}.
This completes the proof. \qed

Henceforth, we assume that we are working with $\gamma = \gamma_0$ and
$\sigma = \sigma_0$ and for the ease of writing we drop $(\gamma,\sigma)$
from our notation unless required.

By translation invariance of our model, we have
\begin{equation*}
\{ g_m (\bu_n ) : m \geq 0 \} \stackrel{d}{=}  \bu_n + \{  g_m ({\bf 0}) : m \geq 0 \}.
\end{equation*}
Using Proposition \ref{lemma:Bmotion}, we conclude the following corollary:
\begin{cor}
\label{cor:anySinglePath}
For any $ \bu \in \R^2 $ and a sequence  $ \{ \bu_n \in \Z^2 : n \geq 1 \} $ such that
$ (\bu_n (1) / n\sigma, \bu_n(2) / n^2\gamma ) \rightarrow \bu$ as $n\rightarrow \infty$,
we have
\begin{equation*}
  \pi^{\bu_n}_n  \Rightarrow  B^{\bu} \text{ in } (\Pi, d_{\pi}).
\end{equation*}
\end{cor}

We now show that, if two paths start close to each other on the x-axis,
they  converge to the same Brownian motion.
%their distance, under scaling, will converge to $0$ in probability.
\begin{prop}
\label{prop:sandwich:DistConvNoHist}
Let $ \bu_n = (u_n , 0) $, $ \bv_n = (v_n, 0) \in \Z^2$ be such that $ u_n < 0 < v_n $ and $  (v_n - u_n )/n \to 0 $. Then,
\begin{equation}
\label{eqn:Sandwich:DistConvNoHist}
 ( \pi_n^{\bu_n}, \pi_n^{\bv_n}) \Rightarrow ( B^{{\bf 0}}, B^{{\bf 0}}).
\end{equation}
\end{prop}

% To show the convergence of $ d_{\Pi} ( \pi_n^{\bu_n}, \pi_n^{\bv_n}) $ in probability,
% it is enough to prove that distance of the starting points of paths $  \pi_n^{\bu_n} $ and
% $  \pi_n^{\bv_n} $  converges to $ 0$ in probability and
% the supremum distance of $ \pi_n^{\bu_n} $ and $  \pi_n^{\bv_n} $ on any compact set $ [0, t] $, converges to $0$ in
% probability. Clearly, distance of the starting points of paths converges to $0$ (as a real
% sequence) by the choice of $ \bu_n $ and $ \bv_n$.

\noindent{\bf Proof}: By Corollary \ref{cor:anySinglePath},  $  \pi^{ {\bu}_n}_n
\Rightarrow B^{\mathbf{0}} $ and $ \pi^{ {\bv}_n}_n \Rightarrow B^{\mathbf{0}} $. Therefore,
$ \{  \pi^{ {\bu}_n}_n  : n \geq 1 \} $ and $ \{  \pi^{ {\bv}_n}_n  : n \geq 1 \} $ are both tight in
$ (\Pi, d_{\Pi}) $, and hence $ \{ ( \pi^{ {\bu}_n}_n ,  \pi^{ {\bv}_n}_n ) : n \geq 1 \} $ is tight
in the product space.
Now, consider  any convergent subsequence and assume  that $ (B , \tilde{B})$ is the subsequential limit.
Since $ \pi^{ {\bu}_n}_n \Rightarrow B^{\mathbf{0}} $ and $  \pi^{ {\bv}_n}_n \Rightarrow B^{\mathbf{0}}  $,
both $ B $ and $ \tilde{B} $ are standard Brownian motions starting from the origin.
Using Skorohod's representation theorem, we can couple so that the convergence is almost sure.
Furthermore, by the non-crossing property of the path family, $  \pi^{ {\bu}_n}_n  (s) \leq
\pi^{ {\bv}_n}_n (s) $ for all $ s \geq  0 $. Hence,
we have that $ B (s) \leq \tilde{B} (s) $ for all $ s \geq 0$. This implies that $ B(s) = \tilde{B} (s) $
for all $ s \geq 0 $ almost surely.
% This implies that $ \P ( \sup_{ s \in [0 ,t] }  |  \pi^{ {\bu}_n}_n  (s)
%- \pi^{ {\bv}_n}_n  (s) | > \epsilon ) \to 0 $ as $ n \to \infty$.
\qed

For verifying condition ($I_1$), we require an estimate on the displacements of
paths in the presence of some information.
Our next proposition estimates the distance traversed by a path (either
laterally or vertically) in terms of the height of the region whose information
is known.
\begin{prop}
\label{prop:HistPathStayClose}
Let ${\cal R} \subseteq \mathbb{H} ( m^{\beta})$ for some $m \geq 1$ and $ 0 <
\beta < 1/2$. For the path $\pi^{\bf{0}}$ staring from the origin ${\bf {0}} =
(0,0)$,
given any configuration  on $ \cal R $ and any $ \delta > 2 \beta $ we have
\begin{equation*}
 \P \bigl( \sup \{ | \pi^{\bf{0}} (s)  |  : 0 \leq s \leq
m^{\beta} \} \geq m^{\delta} \bigl|  \{  U_{\bw} : \bw \in {\cal R}\}  \bigr)
\leq C_{23} m^{\beta} \exp ( - C_{24} m^{\beta} )
\end{equation*}
where $ C_{23} $ and $ C_{24}$ are positive constants.
\end{prop}

% In our next proposition, we consider the case when we are given some history.
% We do not impose any
% condition on the history, except its height.  For each $ n \geq 1 $,
% let $ \Delta_n  \subseteq \mathbb{H} ( m^{\beta} +r_n) $ $( 0 < \beta < 1/2)$ where $ r_n \in \Z$. % be any finite union of isosceles triangles, each triangle
% % having base parallel to the $x$-axis.
% Given any configuration on $ \Delta_n$, we prove that the path from $ \bu =
% (u,r_n)$ will not traverse too much distance (either laterally or vertically)
% until it crosses the line $ \{ y = r_n +  \lfloor m^{\beta} \rfloor + 1 \} $.
% \begin{prop}
% \label{prop:HistPathStayClose}
% Suppose that $\Delta_n $ is as described above. Let $ \bu = (u, r_n) $.
% Given any configuration  on $ \Delta_n $, for any $ \delta > 2 \beta $, as $ n \to \infty$,
% \begin{equation*}
%  \P \bigl( \sup \{ | \pi^{\bu} (s) - u  |  : r_n \leq s \leq r_n + m^{\beta} \} \geq n^{\delta} \bigl|  \{  U_{\bw} : \bw \in \Delta_n \}  \bigr) \to 0.
% \end{equation*}
% \end{prop}

\noindent{\bf Proof}: Consider the horizontal line $
\{ y = \lfloor m^{\beta } \rfloor + 1 \} $  lying above the region $\cal R$.
For the construction of the path $\pi^{\bf{0}}$ there is no information
regarding the configuration on the lattice points on or above this line.
Consider the following events:
\begin{align*}
 E_m^{(1)} & := \bigcap_{ i = 1}^{ \lfloor m^{\beta } \rfloor +1  } \bigcup_{ u
=   2(i-1) \lfloor m^{\beta } \rfloor+1}^{  (2i-1) \lfloor m^{\beta } \rfloor }
 \bigl\{  (u,  \lfloor m^{\beta } \rfloor +1  ) \in V \bigr\};  \\
E_m^{(2)} & := \bigcap_{ i = 1}^{ \lfloor m^{\beta } \rfloor +1  }  \bigcup_{u
=  - (2i-1) \lfloor m^{\beta} \rfloor}^{ - 2(i-1) \lfloor m^{\beta } \rfloor
-1}
\bigl\{   (u,  \lfloor m^{\beta } \rfloor +1   ) \in V  \bigr\} .
\end{align*}
On the event $  E_m^{(1)} $ ($  E_m^{(2)}$),
there are one or more open vertices in each of
the blocks
 $  [2(i-1)
\lfloor m^{\beta} \rfloor + 1, (2i-1) \lfloor m^{\beta} \rfloor ] \times
\{\lfloor m^{\beta } \rfloor + 1\}$
(respectively from $ [- (2i-1) \lfloor m^{\beta} \rfloor ,   -2( i -1)
\lfloor m^{\beta}
\rfloor - 1]  \times
\{\lfloor m^{\beta } \rfloor + 1\}$) of size $ \lfloor m^{\beta} \rfloor
$ for $ i = 1, 2,
\dotsc,  \lfloor m^{\beta} \rfloor +1 $.  Clearly, $ E_m := E_m^{(1)}  \cap E_m^{(2)}$ depends only on the
uniform random variables
on $ \{ y = \lfloor m^{\beta } \rfloor + 1 \} $ and hence, is independent of
the history. Further, $ \P ( E_m^{(1)} )
= \P ( E_m^{(2)} )  =  \bigl( 1 - (1-p)^{  \lfloor m^{\beta} \rfloor } \bigr)^{  \lfloor m^{\beta} \rfloor +1 }
\geq 1 -  (\lfloor m^{\beta} \rfloor +1) (1-p)^{  \lfloor m^{\beta} \rfloor }$.
Therefore $ \P ( E_m^c) \leq  C_{23} m^{\beta} \exp ( - C_{24} m^{\beta} )  $ for suitable
choice of $ C_{23} $ and $ C_{24}$.

Fix any $ \delta \in ( 2\beta , 1 ) $. 
Let $l = \min\{j : h^j({\bf 0})(2) \geq \lfloor m^{\beta } \rfloor + 1\}$. 
Note that at every step the path moves a distance at least $1$ in the $y$ co-ordinate, 
hence $l \leq \lfloor m^{\beta } \rfloor + 1$.
To complete the proof, it is enough to show that on the set
$ E_m$, we have $\{h^j({\bf 0}): 0 \leq j \leq l\} \subseteq [-m^{\delta}, m^{\delta}] 
\times [0, 2( \lfloor m^{\beta } \rfloor + 1) ] $.

On the event $E_m$, the existence of an open vertex in the first block $ [1,   \lfloor m^{\beta} \rfloor]  \times\{
\lfloor m^{\beta} \rfloor + 1 \}$ ensures that  $ || h^1 ({\bf 0}) ||_1 \leq 2m^{\beta} + 1$.
% Since the path moves vertically by a distance at least $1$, $h^1({\bf 0})$ must be to the left
% of the line joining $ ( 2 \lfloor m^{\beta}\rfloor, 1) $ and $ (  \lfloor m^{\beta} \rfloor,  \lfloor m^{\beta} \rfloor+1 )$
% and to the  right
% of the line joining $ ( -2 \lfloor m^{\beta}\rfloor, 1) $ and $ (  -\lfloor m^{\beta} \rfloor,  \lfloor m^{\beta} \rfloor+1 )$.
The construction of the set $ E_m$ ensures that
this argument can be repeated for each of the steps until the $l$ th step of the path, i.e. the step when the path crosses $\{y = \lfloor m^{\beta } \rfloor + 1\}$.
Since $l \leq  \lfloor m^{\beta } \rfloor +1$, we have that until the $l$ th step the path stays inside the rectangle
$[-( 2 \lfloor m^{\beta } \rfloor +1)  \lfloor m^{\beta } \rfloor, (2 \lfloor m^{\beta } \rfloor +1)  \lfloor m^{\beta } \rfloor] \times
[0,  2(\lfloor m^{\beta } \rfloor +1)]$.
% Finally at the $l$th step, the existence of open vertices on the line
% $[-2  \lfloor m^{\beta } \rfloor, 2  \lfloor m^{\beta } \rfloor]
% \times \{ \lfloor m^{\beta } \rfloor +1\}$ ensures that
% $||h^l(\mathbf{0}) - h^{l-1}(\mathbf{0})||_1 \leq 2m^{\beta}$.  

The proposition follows for any $\delta > 2 \beta$. \qed

Returning to the verification of condition ($I_1$)
we start with a map $o_n: \R^2 \to V$ given by 
\begin{align}
\label{defn:o_n}
 o_n(\bz ) =  \bz_n, 
\end{align}
% where $\bz_n(2) := \lfloor n^2 \gamma \bz(2)\rfloor$ and 
% $\bz_n(1) := \min\{j+\lfloor n \sigma  \bz(1) \rfloor : j \geq 0,  (j+\lfloor n \sigma  \bz(1) 
% \rfloor, $ $\bz_n(2)) 
% \text{ is open}\}$. 
where  
$\bz_n(1) := \min\{j+\lfloor n \sigma  \bz(1) \rfloor : j \geq 0,  (j+\lfloor n \sigma  \bz(1) 
\rfloor, $ $\lfloor n^2 \gamma \bz(2)\rfloor) 
\text{ is open}\}$ and $\bz_n(2) := \lfloor n^2 \gamma \bz(2)\rfloor$.
We now define the path $\zeta_n^{\bz} \in {\cal X}_n$ as follows
\begin{equation}
 \label{defn:zeta_n}
 \mbox{$\zeta_n^{\bz} := \pi_n^{o_n(\bz)}$, for any $\bz \in \R^2$.  
}
\end{equation}
Corollary \ref{cor:anySinglePath}  proves Condition ($I_1$) for $k = 1$.

We proceed to prove it for $k \geq 2$, assuming that it is true for $k-1$.
Fix $\bx^1,\dotsc, \bx^k \in \R^2$, and  without loss of generality
we assume $\bx^k(2) = \min\{\bx^i(2): 1\leq i\leq k\} = 0$.

The strategy we adopt  is to show
that until the time when the $k$ th path comes {\em close to one of the other $(k - 1)$ paths},
it can be approximated by an independent path with the same distribution as itself,
and after that time, it quickly coalesces with the path which is close to it and
both of them converge to the same Brownian motion.

Following the ideas introduced in Ferrari {\it et al.}  \cite{FFW05}, we
 consider the product metric space $(\Pi^k, d^k_{\Pi})$ where
\begin{equation*}
d^k_{\Pi} \bigl( (\pi_1,\ldots,\pi_k),(\theta_1,\ldots,\theta_k) \bigr) := \sum_{i=1}^{k}
d_{\Pi}(\pi_i,\theta_i).
\end{equation*}
We define a subset $A$ of $\Pi^k$ as follows:
\begin{align*}
A = & \Bigl\{  (\pi_1,\dotsc,\pi_k) \in \Pi^k :  \text{ such that } \\
& a)  \quad \pi_k(\sigma_{\pi_j}) \neq \pi_j(\sigma_{\pi_j}) \text{ for all } j \neq k;\\
& b)  \quad t^{k} := \inf\{t > \max\{\sigma_{\pi_i}, \sigma_{\pi_k}\} : \pi_i (t) = \pi_k(t) 
\text{  for some } 1 \leq i\leq k-1 \}< \infty; \\
& c)  \quad \text{ for any }\delta > 0 \text{ there exist }1 \leq i \leq k-1, t^{k} - \delta < t < t^{k} < s < t^{k} + \delta 
% t \in (t^{k} - \delta, t^k) \text{ and } s \in (t^k, t^{k} + \delta)
\\
& \quad \quad \text{ such that } (\pi_k(t)- \pi_i(t))(\pi_k(s)- \pi_i(s))< 0 \text{ where }
\pi_i (t^k) = \pi_k(t^k) \\
& \quad \quad \text{ and } \pi_j (t^k) \neq \pi_k(t^k)\text{ for all }1 \leq j < i \Bigr\}.
\end{align*}
% \begin{align*}
% A = & \Bigl\{  (\pi_1,\dotsc,\pi_k) \in \Pi^k :  \text{ such that } \\
% & a)  \quad \pi_k(\sigma_{\pi_j}) \neq \pi_j(\sigma_{\pi_j}) \text{ for all } j \neq k;\\
% % & b)  \quad t^{k} := \inf\{t > \max\{\sigma_{\pi_i}, \sigma_{\pi_k}\} : \pi_i (t) = \pi_k(t) \text{  for some } 1 \leq i\leq k-1\}< \infty; \\
% & b)  \quad \text{ there exist }t^k \in (\sigma_{\pi_k},\infty) ,1\leq i\leq k-1 
% \text{ such that }\pi_i (t^k) = \pi_k(t^k) \text{ and }\\
% & \quad \quad \pi_j (t^k) \neq \pi_k(t^k)\text{ for all }1 \leq j < i; \\
% & c)\quad \text{ for any }\delta > 0\text{ there exist }t \in (t^{k} - \delta, t^k) \text{ and }
% s \in (t^k, t^{k} + \delta)\text{ such that } \\
% & \quad \quad (\pi_k(t)- \pi_i(t))(\pi_k(s)- \pi_i(s))< 0 \}.
% \end{align*}
Note that $ A$ consists of all $k$-tuples of continuous paths such that
the $k$ th path intersects at least one of the other $k -1 $
paths $\pi_1,\dotsc,\pi_{k-1}$ and it immediately crosses one particular such path.
Let $B^{\bx^k}$ be a standard Brownian motion starting
at $\bx^k$ and independent of $W^{\bx^1}, \dotsc, W^{\bx^{k - 1}}$, the independent coalescing Brownian
motions starting from $\bx^1,\dotsc,\bx^{k-1}$.  From the path property
of independent Brownian motions, we have
\begin{equation}
\label{BrownInCoaSetProb}
 \P\bigl[   (W^{\bx^1}, \dotsc, W^{\bx^{k-1}}, B^{\bx^k}) \in A \bigr] = 1.
\end{equation}
We define a {\em coalescence map} $ f:  \Pi^k \to \Pi^k$ as follows:
\begin{equation*}
f(\pi_1,\dotsc,\pi_{k}) :=
\begin{cases}
(\pi_{1}, \dotsc,\pi_{k - 1},  \overline{\pi}_{k}) & \text{ for } (\pi_1, \dotsc, \pi_{k}) \in A \\
(\pi_1, \dotsc, \pi_{k}) & \text{ otherwise }
\end{cases}
\end{equation*}
with
\begin{align*}
\overline{\pi}_{k}(t) :=
\begin{cases}
\pi_{k}(t) & \text{ for } t \leq t^{k}\\
\pi_{i}(t) & \text{ for } t > t^{k}
\end{cases}
\end{align*}
where $i$ is the index such that $\pi_{i}(t^k)= \pi_{k}(t^k)$  and $\pi_{j}(t^k) \neq \pi_{k}(t^k)$
for all $1 \leq j < i$. It follows that 
\begin{equation}
\label{StBrownToCoaBrown}
 f (W^{\bx^1}, \dotsc,W^{\bx^{k-1}},B^{\bx^k}) \stackrel{d}{=}
 (W^{\bx^1}, \dotsc, W^{\bx^k}) .
\end{equation}

Next, we define a sequence of subsets of $\Pi^k$
where the $k$ th path {\em comes close to one of the $k-1$ paths} and a sequence of
`coalescing functions'. The idea of the subsets of $\Pi^k$ and coalescing functions
is motivated from Ferrari {\it et al.} \cite{FFW05}. Ferrari {\it et al.} considered
a sequence of subsets of $\Pi^k$ where any two of the $k$ paths come close to
each other and defined a sequence of coalescing maps such that before
coalescing those two paths are independent. As we proceed by method of induction, we
consider subsets of $\Pi^k$ where the $k$ th path comes close to one of the $k-1$ paths.
Our coalescing map  ensures that the probability that before coalescence
the $k$ th path is independent of the $k-1$ paths converges to $1$.

We fix $ \alpha \in (0,1/2)$ for the rest of this section.
For $n \geq 1$, define
\begin{align}
 \label{Analpha}
A_n^\alpha = & \Bigl\{  (\pi_1,\dotsc,\pi_k) \in \Pi^k :  \text{ such that } \nonumber \\
% & a)  \quad \pi_1(\sigma_{\pi_1}) \leq \pi_2(\sigma_{\pi_2}) \leq \cdots \leq \pi_{k - 1}(\sigma_{\pi_{k-1}});  \nonumber \\
& \quad t^{k}_n : = \inf\{t \geq \max\{\sigma_{\pi_i}, \sigma_{\pi_k}\}:
|\pi_i(t)-\pi_k(t)|\leq n^{\alpha - 1} \nonumber \\
& \quad \quad \quad\text{  for some }1 \leq i\leq k-1\} < \infty \Bigr\}.
\end{align}
%
%  let $ \subseteq C[0,\infty)^k$ be such that
% \begin{flalign}
%
% (\pi_1,\dotsc,\pi_k)\in A_n^\alpha &\text{ if and only if} &
% \end{flalign}
% s\begin{itemize}
%  \item[a)]  $\pi_1(0) < \pi_2(0) < \cdots < \pi_k(0)$,
%  \item[b)]  $t^{(i,j)}_n : = \inf\{t\geq 0: |\pi_i(t)-\pi_j(t)|\leq 3n^{\alpha - 1}\} < \infty$ for all $1\leq i< j \leq k$,
%  \item[c)]
% $|t^{(i_1,j_1)}_n - t^{(i_2,j_2)}_n| > \frac{1}{n^2}$ for all $1\leq i_1<j_1\leq k$ and $1\leq i_2<j_2\leq k$ with $(i_1,j_1)\neq (i_2,j_2)$.
%  \end{itemize}
% % Let $T^{(i_1,j_1)}_n < \cdots < T^{(i_{k-1},j_{k-1})}_n$ be the ordering of the `$\alpha$-coalescing
% % times' with $j_l  \in \{i_l+1,\dotsc,k\}\setminus\{j_1,\dotsc , j_{l-1}\}\}$ for every $1 \leq l \leq k$.
We now define the `{\em $\alpha$-coalescence map}'
$f_n^{(\alpha)}:  \Pi^k \to \Pi^k$, as follows:
\begin{equation*}
f_n^{(\alpha)} (\pi_1,\dotsc,\pi_{k}) :=
\begin{cases}
(\pi_1,\dotsc,\pi_{k-1},\overline{\pi}_{k}) & \text{for }(\pi_1,\dotsc,\pi_{k}) \in A_n^\alpha \\
(\pi_1,\dotsc,\pi_{k}) & \text{ otherwise }
\end{cases}
\end{equation*}
with
\begin{align*}
\overline{\pi}_{k}(t) :=
\begin{cases}
 \pi_{k}(t) & \text{ for } t \leq t^{k}_n\\
 \pi_{k}(t^{k}_n) +
\frac{( t - t^{k}_n)}{s_n^k - t^{k}_n} \bigl[ \pi_{i}( s_n^k )
-  \pi_{k}( t^{k}_n) \bigr] &  \text{ for } t^{k}_n < t \leq s_n^k \\
 \pi_{i}(t) & \text{ for } t > s_n^k
\end{cases}
\end{align*}
where  $ s_n^k = (\lfloor n^2\gamma t^{k}_n\rfloor + 1)/(n^2\gamma) $ and
$i$ is the index such that  $|\pi_{i}( t^{k}_n) - \pi_{k}( t^{k}_n)| \leq n^{\alpha - 1}$ and
$|\pi_{j}( t^{k}_n) - \pi_{k}( t^{k}_n)| > n^{\alpha - 1}$ for all $ 1 \leq j < i$.
% Since From the definition of $f_n^{(\alpha)}$ it follows that for all $t \geq \sigma_{\pi_i}$
% $$
% \min\{\pi_{k}(t), \pi_{i}(t)\}\leq \overline{\pi}_{k}(t)\leq \max\{\pi_{k}(t), \pi_{i}(t)\}.
% $$

Before proceeding, we state the following deterministic
lemma (which is a slightly stronger version of Lemma 19 of
Coletti {\it et al.}\/ \cite{CFD09}). The proof of this lemma has been
relegated to the Appendix and it
will be used later in the proof of Proposition \ref{prop:fddbweb2}.
\begin{lemma}
\label{lemma:fnconv}
Let $(\pi_1, \dotsc , \pi_k) \in A$ and $ \{ (\pi_{1, n}, \dotsc,  \pi_{k, n}) : n \geq 1 \}
\subseteq \Pi^k $ be such that for all $1 \leq i \leq k$, $ d_{\Pi}(\pi_{i, n},\pi_i) \rightarrow 0 $ as $ n \to \infty $. Then,
for $n$ large enough, we have $(\pi_{1,n}, \dotsc , \pi_{k,n}) \in A^{\alpha}_n$
and
$\lim_{n\rightarrow \infty} t^{k}_n = t^{k}$, where $t^{k}, t^{k}_n$
are as defined above. Further,
\begin{equation}
d_{\Pi}^k \Bigl( f_n^{(\alpha)} (\pi_{1,n}, \dotsc , \pi_{k,n}), f (\pi_{1}, \dotsc ,  \pi_{k} ) \Bigr) \to 0 \text{ as }
n\rightarrow \infty.
\end{equation}
\end{lemma}

We  now describe a construction which will be used to prove the general case. Let  $\{U^{r}_{\bw} : \bw \in \Z^2\}$
and $\{U^{g}_{\bw} : \bw \in \Z^2\}$ be two independent collections of i.i.d. $U(0,1)$ random variables.
Given a set of points $ {\bx^1}, \dotsc,  {\bx^k} \in \R^2$, 
let
$ \{ (\bx^1_n,\dotsc, \bx^k_n)  : n \geq 1 \} $ be such that for all $ i=
1,2, \dotsc, k $,  $ \bx_n^i \in \Z^2 $ for  $ n \geq 1$ with $ \bx^i_n(2) \geq 0 $ and $ (\bx^i_n (1) / n\sigma,
\bx^i_n(2) / n^2\gamma ) \rightarrow \bx^i$
as $n\rightarrow \infty$.
We construct the paths
$ \pi^{1}, \dotsc, \pi^{k-1}$  starting from $  \bx^1_n,\dotsc, \bx^{k-1}_n$ using only the collection 
$ \{U^{r}_{\bw} : \bw \in \Z^2\}$,
while, for  the construction of the path $\tilde{\pi}^{k}$  
starting from $ \bx^k_n$ we use the collection $ \{U^{g}_{\bw} : \bw \in \Z^2\}$.
 From the independence of the collections of uniform
random variables, the scaled paths $(\pi^{1}_n, \dotsc, \pi^{k-1}_n) $ and the scaled path $ \tilde{\pi}^{k}_n$
are independent. Further,

% We  now describe a construction which will be used to prove the general case. Let  $\{U^{r}_{\bw} : \bw \in \Z^2\}$
% and $\{U^{g}_{\bw} : \bw \in \Z^2\}$ be two independent collections of i.i.d. $U[0,1]$ random variables.
% Given  $ {\bx^i} \in \R^2$, $i = 1, \dotsc , k$, we construct the paths
% $ \pi_n^{i} = \zeta_n^{x^i} \in {\cal X}_n$, $i = 1, \dotsc , k-1 $  starting from $ o_n(\bx^i)$ as in (\ref{defn:o_n})
%  using only the collection $ \{U^{r}_{\bw} : \bw \in \Z^2\}$,
% while, for  the construction of the path $\tilde{\pi}^{k}$  starting from $ \bx^k_n$ we use the collection $ \{U^{g}_{\bw} : \bw \in \Z^2\}$.
%  From the independence of the collections of uniform
% random variables, the scaled paths $(\pi^{1}_n, \dotsc, \pi^{k-1}_n) $ and the scaled path $ \tilde{\pi}^{k}_n$
% are independent. Further,
\begin{equation}
\label{eqn:SameDistk-1Andk}
(\pi^{1}_n, \dotsc, \pi^{k-1}_n) \stackrel{d}{=} (\pi^{\bx^1_n}_n, \dotsc, \pi^{\bx^{k-1}_n}_n) \quad \text{ and } \quad
\tilde{\pi}^{k}_n \stackrel{d}{=}  \pi^{\bx^k_n}_n.
\end{equation}

% the marginal distribution
% with  where both $\pi^{\bx^k_n}_n$ and
% $\tilde{\pi}^{\bx^k_n}_n$ marginally have the same distribution but $\tilde{\pi}^{\bx^k_n}_n$ is
% independent of $(\pi^{\bx^1_n}_n, \dotsc, \pi^{\bx^{k-1}_n}_n)$.
% Let  First
% Therefore, $n$ th order diffusively scaled path $\tilde{\pi}^{\bx^k_n}_n$ is  independent of the
% $n$ th order diffusively scaled paths $(\pi^{\bx^1_n}_n, \dotsc, \pi^{\bx^{k-1}_n}_n)$.

Next, we consider the region $ E^{(r)} $ which is explored by the paths   $ \pi^{1},
\dotsc, \pi^{k-1}$, constructed using the collection
$ \{U^{r}_{\bw} : \bw \in \Z^2\}$ only. On the complement set of $ E^{(r)} $, we consider the collection
$\{U^{g}_{\bw} : \bw \in \Z^2\}$ and construct the path, $ \pi^{k} $, starting
from $ \bx^k_n $. More precisely, the set %for $\bw \in \Z^2$ and for $j \geq 1$,  let $h^{j,r}(\bw) := h^{1,r}(h^{j-1,r}(\bw))$
%where $h^{1,r}(\bw)$ denotes $h(\bw)$ based on $\{U^r_{\bw}: \bw\in \Z^2\}$.
\begin{align*}
 E^{(r)}   := \bigcup_{1\leq i \leq k - 1} \bigcup_{ m \geq 0}
S^{+}\Bigl( h^{m}(\bx^i_n),|| h^{m}(\bx^i_n) - h^{m + 1}(\bx^i_n)||_1 \Bigr)
%\setminus\mathbb{H}\bigl(h^{m}(\bx^i_n) (2) \bigr)
\end{align*}
represents the explored region by the paths $ \pi^{1}, \dotsc, \pi^{k-1}$ using the
collection $ \{U^{r}_{\bw} : \bw \in \Z^2\}$ only. We define $\{ U^{\text{mixed}}_{\bw}:
\bw\in \Z^2\}$ by
\begin{align*}
U^{\text{mixed}}_{\bw} :=
\begin{cases}
U^{r}_{\bw} & \text{ if }\bw \in E^{(r)} \\
U^{g}_{\bw} & \text{ otherwise. }
\end{cases}
\end{align*}
Let $\pi^{k}$ be the path starting from $\bx^k_n$ constructed using the collection $ \{ U^{\text{mixed}}_{\bw} :
\bw \in \Z^2\}$. We also observe that the distribution of $ \pi^k$, given the realization of the uniform random variables
in $ E^{(r)} $, is the same as the conditional distribution of $ \pi^{\bx^{k}_n} $ given the paths
$ \pi^{\bx^1_n}, \dotsc, \pi^{\bx^{k-1}_n}. $ Hence from the above observation and definition (\ref{defn:zeta_n}),
\begin{equation}
 \label{eqn:SameDistk}
(\pi^{1}_n, \dotsc, \pi^{k-1}_n, \pi^{k}_n) \stackrel{d}{=} (\pi^{\bx^1_n}_n, \dotsc, 
\pi^{\bx^{k-1}_n}_n, \pi^{\bx^{k}_n}_n)  .
%= (\zeta_n^{\bx^1}, \dotsc , \zeta_n^{\bx^k}).
\end{equation}

%For $1\leq i \leq k-1$, let $\pi^{i,r}_n$ be the $n$ th order diffusively scaled path
%starting from $\bx^1_n$ based on the uniform random variables $\{U^r_{\bw}: \bw \in \Z^2\}$
%and $\pi^{k,g}_n$ be the $n$ th order diffusively scaled path
%starting from $\bx^k_n$ based on the uniform random variables
%$\{U^g_{\bw}: \bw \in \Z^2\}$. Since $\pi^{k,g}_n$ depends on the uniform random variables
%$\{U^g_{\bw}: \bw \in \Z^2\}$ only, $\pi^{k,g}_n$ is
%independent of $(\pi^{1,r}_n, \dotsc, \pi^{k-1,r}_n)$ for all $n$. Next we
%construct the joint paths $(\pi^1_n, \dotsc, \pi^k_n)$
%using both
%$\{U^{r}_{\bw} : \bw \in \Z^2\}$ and $\{U^{g}_{\bw} : \bw \in \Z^2\}$.

% For $\bw \in \Z^2$ and $l > 0$
% let $B(\bw, l)$ denote the upper $||\quad||_1$ ball of
% radius $l$ given by
% $B(\bw , l):= \{\bu : \bu \in \Z^2, \bu(2)>\bw(2), ||\bu-\bw||_1 \leq l\}$.

%based on  $\{U^r_{\bw}: \bw \in \Z^2\}$.

\begin{prop}
\label{prop:fddbweb2}
We have, as $ n \to \infty $,
 \begin{align*}
(a)~ & f^{(\alpha)}_n (\pi^{1}_n, \dotsc, \pi^{k-1}_n, \tilde{\pi}^{k}_n)
\Rightarrow  (W^{\bx^1}, \dotsc, W^{\bx^k}); \\
(b)~ & f^{(\alpha)}_n (\pi^{1}_n, \dotsc, \pi^{k-1}_n, \pi^{k}_n)
\Rightarrow  (W^{\bx^1}, \dotsc, W^{\bx^k}); \\
(c)~ &  (\pi^{1}_n, \dotsc, \pi^{k-1}_n, \pi^{k}_n)
\Rightarrow  (W^{\bx^1}, \dotsc, W^{\bx^k}).
 \end{align*}
\end{prop}
Since $(o_n(\bz)(1)/(n\sigma), o_n(\bz)(2)/(n^2\gamma)) \to \bz$ almost surely,  
where $o_n(\bz)$ is defined in (\ref{defn:o_n}), by (c) of the proposition above, we have
$$
(\zeta_n^{\bx^1}, \dotsc , \zeta_n^{\bx^k}) \Rightarrow  (W^{\bx^1}, \dotsc, W^{\bx^k}),
$$
which verifies ($I_1$).

\noindent \textbf{Proof}: From Corollary \ref{cor:anySinglePath} and (\ref{eqn:SameDistk-1Andk}) it follows that the
scaled path $\tilde{\pi}^{k}_n$ converges in distribution to  $ B^{\bx^k}$, the standard Brownian motion starting at $\bx^k$.
Using the induction hypothesis and (\ref{eqn:SameDistk-1Andk}) we have that $(\pi^{1}_n, \dotsc, \pi^{k-1}_n)$
converges in distribution to $(W^{\bx^1},\dotsc, W^{\bx^{k-1}})$. From the independence of paths
 we have $ (\pi^{1}_n, \dotsc, \pi^{k-1}_n,  \tilde{\pi}^{k}_n)$
converges in distribution to $(W^{\bx^1},\dotsc, W^{\bx^{k-1}}, B^{\bx^k})$
 where $B^{\bx^k}$ is  independent of $(W^{\bx^1},\dotsc, W^{\bx^{k-1}})$. Now  Lemma \ref{lemma:fnconv}
and (\ref{BrownInCoaSetProb}) enable us to use the extended continuous mapping theorem
(see Theorem 4.27 in Kallenberg \cite{K02}) to conclude that
\begin{equation*}
  f^{(\alpha)}_n (\pi^{1}_n, \dotsc, \pi^{k-1}_n, \tilde{\pi}^{k}_n)
\Rightarrow  f  (W^{\bx^1},\dotsc, W^{\bx^{k-1}}, B^{\bx^k}) \stackrel{d}{=}
 (W^{\bx^1}, \dotsc, W^{\bx^k})
\end{equation*}
where the last relation follows from (\ref{StBrownToCoaBrown}). This proves (a).

For  $(b)$, let $ f_n^{(\alpha)} (\pi^{1}_n,\dotsc,\pi^{k-1}_n,\pi^{k }_n)
=  (\pi^{1}_n,\dotsc,\pi^{k-1}_n, {\overline{\pi}}^{k}_n)$ and
$  f_n^{(\alpha)} (\pi^{1}_n,\dotsc,\pi^{k-1}_n, $ $\tilde{\pi}^{k }_n)
=  (\pi^{1}_n,\dotsc,\pi^{k-1}_n, \overline{\overline{\pi}}^{k}_n)$.
Note that $ d^k_{\Pi}  \bigl( (\pi^{1}_n,\dotsc,\pi^{k-1}_n, {\overline{\pi}}^{k}_n),
 (\pi^{1}_n,\dotsc,\pi^{k-1}_n, \overline{\overline{\pi}}^{k}_n) \bigr) = d_{\Pi} (  {\overline{\pi}}^{k}_n,
\overline{\overline{\pi}}^{k}_n) $.  Hence it is enough to show that $ d_{\Pi} (  {\overline{\pi}}^{k}_n,
\overline{\overline{\pi}}^{k}_n)  \stackrel{\P}{\to} 0 $ as $n \to \infty$. Since both the paths start at the same point, it is enough
to prove that,
\begin{equation}
\label{RTP}
\text{ for any $ t > 0$, $ \sup \{ | {\overline{\pi}}^{k}_n (s) - \overline{\overline{\pi}}^{k}_n (s) |
: 0 \leq s \leq t \} \prob 0 $ as $n \to \infty$.}
\end{equation}
Towards this end, we show that on a set whose probability
converges to $1$,  $ {\overline{\pi}}^{k}_n (s) $ and $ \overline{\overline{\pi}}^{k}_n(s)$  agree for  $s \in [0,t]$.

For any $s > 0 $ and $ i = 1, 2, \dotsc, k $, set
\begin{equation*}
 l_i^{(s)} = l_i^{(s)} (n) := \min\{ j \geq 0: T_{j} (\bx^i_n) \geq n^2 \gamma s \}
\end{equation*}
where $T_{j}(\bu)$ is the $j$ th regeneration time for path starting from $\bu$ defined  in  (\ref{def:TimekthReg}).
Here $  l_i^{(s)}  $ stands for the number of regenerations needed for the $i$ th path $ \pi^{i}$ (starting
from $ \bx^i_n$) to cross the line $\{y = n^2 \gamma s\} $.
Now, for $ 0 < \beta < \alpha $ and any  $ s > 0$, define
the event
\begin{equation*}
A_{n}^{(s)} := \{W_{j}(\bx^i_n) < n^{\beta}:\text{ for all }1\leq i\leq k, 1 \leq j \leq l_i^{(s)}\}.
\end{equation*}
with $W_{j}(\bu)$  as defined in (\ref{r-W2}). On the event $ A_{n}^{(s)} $, each of the
regeneration steps till the $l_i^{(s)}  $ th regeneration of all the $k$ paths,
is of length at most $ n^{\beta}$.  Since $  l_i^{(s)} \leq \lfloor n^2 \gamma s \rfloor +1$,  and using the fact
that  the individual regenerations are i.i.d., having bounds on tail probabilities given in Proposition \ref{prop:WidthTail},
we have, as $ n \to \infty$,
\begin{align}
\label{RTP2}
\P\bigl( (A_{n}^{(s)})^c \bigr) & = \P \Bigl( \bigcup_{1 \leq i \leq k}\bigcup_{1 \leq j \leq l_i^{(s)} }
\{ W_{j}(\bx^i_n) \geq n^{\beta} \} \Bigr) \nonumber\\
& \leq k(\lfloor n^2 \gamma s \rfloor +1) \P(W_1(\mathbf{0})
\geq n^\beta) \to 0.
\end{align}

If $ t^{k}_n \geq t $, on the event  $ A_{n}^{(t)} $,  
\begin{itemize}
\item[(a)] from the definition of $ t^{k}_n $, for  $ s \leq n^2 \gamma t$, we have
$\min_{1 \leq i \leq k-1} |\tilde{\pi}^k(s) -  \pi^i(s)| \geq  n^{\alpha}$,
\item[(b)] the restriction on the sizes of the regeneration steps by $ n^{\beta} $ together with the choice $ \beta < \alpha $ ensures that  the path
 $ \tilde{\pi}^k $ has not explored the region $ E^{(r)}$ before it crosses $\{y = n^2 \gamma t\}$.
 \end{itemize}
Thus
the paths $ \tilde{\pi}^k $ and  ${\pi}^k$ agree on $[0,n^2 \gamma t]$ and 
so  $  {\overline{\pi}}_n^{k} (s) = \overline{\overline{\pi}}_n^{k} (s) $
for $ 0 \leq s \leq t $.

If $ t^{k}_n \leq t $, on the event  $ A_{n}^{(t^{k}_n)} $, we have
\begin{itemize}
\item[(a)] from the definition of $ t^{k}_n $, for every $ s \leq n^2 \gamma t^{k}_n $, we have
$\min_{1 \leq i \leq {k-1}} |\tilde{\pi}^k(s) -  \pi^i(s)| \geq  n^{\alpha}$,
\item[(b)] the restriction on the sizes of the regeneration steps by $ n^{\beta} $ 
together with the choice $ \beta < \alpha $ ensures that  the path
 $ \tilde{\pi}^k $ has not explored the region $ E^{(r)}$ before it crosses $\{y = n^2 \gamma t^{k}_n \}$.
 \end{itemize}
Thus
the paths $ \tilde{\pi}^k $ and  ${\pi}^k$ agree on $[0,n^2 \gamma t^{k}_n ]$ and thereby  $  {\overline{\pi}}^{k} (s) = \overline{\overline{\pi}}^{k} (s) $
for $ 0 \leq s \leq  t^{k}_n  $. The rest of the path ${\pi}^k$ (from $ n^2 \gamma t^{k}_n $ onwards)
depends only on the position $ \tilde{\pi}^k ( n^2 \gamma t^{k}_n) $ and the paths $ (\pi^{1}, \dotsc, \pi^{k-1}) $ and
hence, by the definition of the
$\alpha$-coalescing map, we have
\begin{equation*}
  f_n^{(\alpha)} (\pi^{1}_n,\dotsc,\pi^{k-1}_n,\pi^{k }_n) (s) =
  f_n^{(\alpha)} (\pi^{1}_n,\dotsc,\pi^{k-1}_n, \tilde{\pi}^{k }_n) (s)
\end{equation*}
for $ s \in [ t^{k}_n, t] $. Hence, $ {\overline{\pi}}^{k}_n (s) = \overline{\overline{\pi}}^{k}_n (s) $
for $ 0 \leq s \leq t$.

Since $  A_{n}^{(t)} \subseteq A_{n}^{(t^{k}_n)}  $ when $ t^{k}_n \leq t $,   (\ref{RTP2}) along with the two observations above implies (\ref{RTP}). This completes the proof of part (b).

For $(c)$, we show that
$ d^k_{\Pi} \bigl( (\pi^{1}_n,\dotsc,\pi^{k-1}_n, \overline{\pi}^{k}_n) ,
(\pi^{1}_n,\dotsc,\pi^{k-1}_n,\pi^{k }_n) \bigr)  = d_{\Pi} (  \overline{\pi}^{k}_n , \pi^{k }_n) \prob 0$ as $n \to \infty$.
Again it is enough to prove,
\begin{equation}
\label{RTP3}
\text{for any $ t > 0$,  $ \sup \{ |  \pi^{k }_n (s) - \overline{\pi}^{k}_n (s) | : 0 \leq s \leq t \}
\prob 0 $ as $n \to \infty$.}
\end{equation}
 Suppose $i_0 :=  \min \{ j :  |\pi^{k}_n (t^{k}_n)-  \pi^{j}_n(t^{k}_n)| \leq n^{\alpha-1} \}$, i.e. $ \pi^{i_0}_n$ 
 is the path with the minimum index which comes $n^{\alpha-1}$ close to $\pi^{k}_n$. Note that
$ \overline{\pi}^{k}_n(s) $
 is obtained by a linear interpolation between $ \pi^{k}_n (t^{k}_n) $
and $ \pi^{i_0}_n (s_n^k) $ for  $ s \in [t^{k}_n, s_n^k] $ and
$ \overline{\pi}^{k}_n(s) = \pi^{i_0}_n(s)  $ for $s \in [s_n^k, \infty)$
where $ s_n^k = (\lfloor n^2\gamma t^{k}_n\rfloor + 1)/( n^2\gamma)$.
Since paths in ${\cal X}_n$ are noncrossing almost surely, we have $ \overline{\pi}^{k}_n (s) \in [ \pi^{i_0}_n(s) , \pi^{k}_n(s)] $
 for $ [ t^{k}_n, \infty) $. Also, note that  both the paths $ {\pi}^{k}_n  $ and $ \overline{\pi}^{k}_n$
start at the same point and agree till $  t^{k}_n $.   Thus (\ref{RTP3}) holds for  $t_n^k \geq t$.

If $ t_n^k < t $ , we have
\begin{equation*}
 \sup \{ |  \pi^{k }_n (s) - \overline{\pi}^{k}_n (s) | : 0 \leq s \leq t \}
\leq \sup \{ |  \pi^{k }_n (s) - {\pi}^{i_0}_n (s) | : t_n^k \leq s \leq t \}.
\end{equation*}
% \begin{equation*}
% \sup_{s \in [t^{k}_n, t]}| \overline{\pi}^{k}_n(s)- \pi^{k}_n(s)|
% \leq \sup_{s \in [t^{k}_n, t]}| {\pi}^{i_0}_n(s)- \pi^{k}_n(s)| .
% \end{equation*}
%
% , we need to show, as $ n \to \infty $,
% \begin{equation*}
%  \P \bigl[ \sup_{s \in [t^{k}_n, t]}| \overline{\pi}^{k}_n(s)- \pi^{k}_n(s)| \geq \epsilon, t^{k, n} < t \bigr] \to 0.
% \end{equation*}
%
% Again it is enough to prove that the restrictions of
% $(\pi^{1}_n,\dotsc,\pi^{k-1}_n, \overline{\pi}^{k}_n) $ and $(\pi^{1}_n,\dotsc,\pi^{k-1}_n,\pi^{k }_n)$
% on $ [0,t]$ converge to the same limit.  From
% the definition of $  f_n^{(\alpha)} $,  $\overline{\pi}^{k}_n (s) = {\pi}^{k}_n (s) $ for all $ s \leq  t^{k}_n$.
% Hence, we need to consider  $ t^{k}_n < t $ only.
Again, we restrict ourselves to the event $ A_{n}^{(t^{k}_n)} $. Let $ \Delta_n $ be the set of vertices explored by the
paths until all of them have crossed the line $ \{y = \lfloor n^2 \gamma t^{k}_n \rfloor \} $, i.e.,
\begin{equation*}
\Delta_n =  \bigcup_{1\leq i \leq k } \bigcup_{ m = 0}^{ \tau_{ l_i^{( t^{k}_n )} } -1}
S^{+}\Bigl( h^{m}(\bx^i_n),|| h^{m}(\bx^i_n) - h^{m + 1}(\bx^i_n)||_1 \Bigr)
%\setminus\mathbb{H}\bigl(h^{m}(\bx^i_n) (2) \bigr) .
\end{equation*}
We observe that,  on the event $  A_{n}^{(t^{k}_n)} $, the set $ \Delta_n $ is contained
in $ \mathbb{H} ( \lfloor n^2 \gamma t^{k}_n  \rfloor + \lfloor n^{\beta} \rfloor )$. Now we
choose two points $ u_n, v_n \in \Z $ such that $ u_n < \pi^{i_0} ( \lfloor n^2 \gamma t^{k}_n  \rfloor)  ,
\pi^{k} ( \lfloor n^2 \gamma t^{k}_n  \rfloor)  < v_n $ and $ (v_n - u_n)/n \to 0 $. Since
$ |\pi^{i_0} ( \lfloor n^2 \gamma t^{k}_n  \rfloor) - \pi^{k} ( \lfloor n^2 \gamma t^{k}_n  \rfloor) | \leq n^{\alpha - 1}$, such a choice of $ u_n, v_n$ is possible for $n$ large.
Set $ \bu_n :=  ( u_n,   \lfloor n^2 \gamma t^{k}_n  \rfloor) $
 and $ \bv_n :=  ( v_n,   \lfloor n^2 \gamma t^{k}_n  \rfloor) $. We consider scaled paths $ \pi_n^{\bu_n} $
and $ \pi_n^{\bv_n} $ and by noncrossing property of the paths, we see that the paths $ \pi_n^{i_0} $ and
$ \pi_n^{k} $  lie between the paths $ \pi_n^{\bu_n} $ and $ \pi_n^{\bv_n} $ from $ t_n^k $ onwards, so that
\begin{equation*}
 \sup \{ |  \pi^{k }_n (s) - {\pi}^{i_0}_n (s) | : t_n^k \leq s \leq t \} \leq \sup \{ |  \pi_n^{\bv_n} (s) - {\pi}_n^{\bu_n} (s) | : t_n^k \leq s \leq t \}.
\end{equation*}

Fix any $ \delta \in (2 \beta, 1) $ and consider the points $\bu_n^{\prime} =
 (u_n - \lfloor n^{\delta} \rfloor -1,  \lfloor n^2 \gamma t^{k}_n  \rfloor + \lfloor n^{\beta } \rfloor +1) $ and $\bv_n^{\prime} =
 (v_n + \lfloor n^{\delta} \rfloor +1,  \lfloor n^2 \gamma t^{k}_n  \rfloor + \lfloor n^{\beta } \rfloor +1) $.
Let
\begin{align*}
 F_n (\bu_n) &:= \bigl\{ \sup \{ | \pi^{\bu_n} (s) - u_n  |  :  \lfloor n^2 \gamma t^{k}_n  \rfloor \leq s \leq
 \lfloor n^2 \gamma t^{k}_n  \rfloor + n^{\beta} \} \leq n^{\delta} \bigr\},\\
 F_n (\bv_n) &:= \bigl\{ \sup \{ | \pi^{\bv_n} (s) - v_n  |  :  \lfloor n^2 \gamma t^{k}_n  \rfloor \leq s \leq
 \lfloor n^2 \gamma t^{k}_n  \rfloor + n^{\beta} \} \leq n^{\delta} \bigr\} .
\end{align*}
By Proposition \ref{prop:HistPathStayClose} (taking ${\cal R } = \Delta_n$ and translating so that $\bu_n$ is the origin) we have $ \P \bigl( F_n (\bu_n)  \bigr) \to 1 $.
Similarly,  $ \P \bigl( F_n (\bv_n)  \bigr) \to 1 $.

On the set $ F_n (\bu_n) \cap F_n (\bv_n)  $, we have that $ \bu_n^{\prime} (1) \leq \pi^{\bu_n} (\lfloor n^2 \gamma t^{k}_n  \rfloor  + \lfloor n^\beta
\rfloor +1),  \pi^{\bv_n} (\lfloor n^2 \gamma t^{k}_n  \rfloor + \lfloor n^\beta \rfloor +1) \leq \bv_n^{\prime} (1) $ and hence, by non-crossing property of paths,
 $ \pi^{\bu_n^{\prime}}(s) \leq \pi^{\bu_n} (s) \leq  \pi^{\bv_n} (s) \leq
\pi^{\bv_n^{\prime}} (s)$ for all $s >   \lfloor n^2 \gamma t^{k}_n  \rfloor + \lfloor n^\beta \rfloor +1$. Using  the estimate in Proposition \ref{prop:HistPathStayClose}, we have
\begin{align*}
 \lefteqn{ \sup \{ | \pi_n^{\bu_n} (s) -  \pi_n^{\bv_n} (s) | : t_n^k \leq s \leq t \} } \\
   & \leq  ( |u_n -v_n| + 2 n^{\delta}) / (n \sigma) %\qquad \qquad \qquad \qquad
 +  \sup \{ | \pi_n^{\bu_n^{\prime}} (s) -  \pi_n^{\bv_n^{\prime}} (s) | : t_n^k + n^{\beta}/( n^2\gamma) \leq s \leq t \}.
\end{align*}

The paths starting from $ \bu_n^{\prime} $
and $ \bv_n^{\prime} $ depend only on the uniform random variables defined on
$ \{ y >   \lfloor n^2 \gamma t^{k}_n  \rfloor + \lfloor n^{\beta } \rfloor + 1 \} $
and hence independent of the realizations in the set $ \Delta_n$.
Using translation invariance of our model and Proposition \ref{prop:sandwich:DistConvNoHist}
we conclude that
\begin{equation*}
 \sup \{ | \pi_n^{\bu_n^{\prime}} (s) -  \pi_n^{\bv_n^{\prime}} (s) | : t_n^k + n^{\beta}/( n^2\gamma) \leq s \leq t \}
  \prob 0.
\end{equation*}
This proves the proposition.
\qed

\subsection{ Verification of $(B_1)$ and $(E_1^{\prime})$}

The verification of condition $(B_1)$ is standard and follows from the same argument as in
Ferrari {\it et al.}\/ \cite{FFW05}. Recall that for $t_0 \in \R, t > 0$ and  $-\infty < a
< b < \infty$ and for $\Gamma \subseteq \Pi$
\begin{equation*}
\eta_{\Gamma}(t_0,t;a,b) = \#\{\pi(t_0 + t): \pi \in \Gamma^{t_0} \text{ and }\pi(t_0)
\in [a,b]\}.
\end{equation*}
By translation invariance, it suffices to consider $\eta_{\overline{\cal X}_n}(0,t;0,\epsilon)$.
Let  $ \bv^n = \bigl( (\lfloor n \sigma \epsilon \rfloor +1) , 0\bigr) $. As noted in \cite{FFW05}, using $(I_1)$, we have
$ \P \bigl( \eta_{\overline{\cal X}_n}(0,t;0,\epsilon) \geq 2 \bigr) \leq \P \bigl( \pi_n^{\bf 0} (t) \neq   \pi_n^{\bv^n} (t) \bigr)
\to \P \bigl( W^{\bf 0} ( t ) \neq W^{ (\epsilon, 0)} (t) \bigr) =  2\phi(\epsilon/\sqrt{t}) - 1 $.

% \begin{align*}
% \lefteqn{\P(\eta_{\overline{\cal X}_n}(0,t;0,\epsilon)\geq 2)}\\
% & \leq \P\{\text{two independent ${\cal X}_n$ paths, one starting at $(0,0)$ and another starting at $(\epsilon, 0)$,}\\
% & \qquad \text{do not meet in time $t$}\}\\
% & \to \P\{\text{two independent Brownian paths, one starting at $(0,0)$ and another starting at $(\epsilon, 0)$,}\\
% & \qquad \text{do not meet in time $t$}\} \text{ as $n \to \infty$}\\
% & = 2\phi(\frac{\epsilon}{\sqrt{t}}) - 1,
% \end{align*}
%
Hence we have
\begin{align*}
%\label{verifyB1}
\limsup_{n \to \infty}\P(\eta_{\overline{\cal X}_n}(0,t;0,\epsilon)\geq 2)
\leq 2 \Phi(\frac{\epsilon}{\sqrt{t}}) - 1 \to 0
\end{align*}
as $ \epsilon \to 0 $ which verifies $(B_1)$.

In order to verify $(E_1^{\prime})$, recall that
$$\overline{\cal X}_{n}^{t_0}= \{ \pi : \pi \in \overline{\cal X}_{n} \text{ with } \sigma_{\pi} \leq  t_0\}.$$
Since $\overline{\cal X}_{n}^{t_0} \subseteq \overline{\cal X}_{n}$ for any $t_0 \in \R$,
and $\{\overline{\cal X}_{n}: n \geq 1\}$ is tight, we have $\{\overline{\cal X}_{n}^{t_0} : {n \geq 1}\}$
is also tight. Let ${\cal Z}^{t_0}$
be a subsequential limit of $\{\overline{\cal X}_{n}^{t_0}: {n \geq 1}\}$. For ease of notation, we assume 
that $\{\overline{\cal X}_{n}^{t_0}: {n \geq 1}\}$ is itself the subsequence
which converges to ${\cal Z}^{t_0}$.

For $ \Gamma \subseteq \Pi$, let  $ \Gamma(s):=
\{ (\pi(s),s) : \pi \in \Gamma^{s} \} \subseteq \R^2$.
For $ t > 0$ define $\Gamma^{s; (s+t)_T} := \{ \pi : \sigma_{\pi} = s+t
\text{ and there exists }\pi' \in \Gamma^s
\text{ such that } \pi(u) = \pi'(u) \text{ for all }u \geq s+t \}$. 
Note $ \Gamma^{s; (s+t)_T} $ is the restriction of paths
in $ \Gamma^{s} $ on $ [s+t, \infty)$.
Following the argument in Newman {\it et al.}\/ \cite{NRS05}, our strategy to check $(E_1^{\prime})$ is to
first show that the point set ${\cal Z}^{t_0}(t_0 + t)$
is locally finite and then using $(I_1)$ we show that $({\cal Z}^{t_0})^{t_0;(t_0 + t)_T} = 
{\cal Z}^{t_0;(t_0 + t)_T}$
has the same distribution as
coalescing Brownian motions ${\cal W}^{{\cal Z}^{t_0}(t_0 + t)}$
starting from a random point set distributed as ${\cal Z}^{t_0}(t_0 + t)$.

\begin{figure}[htb]
\centering
%\begin{center}\leavevmode

\begin{pspicture}(0,-1.0)(14,1.0)

\pscircle[fillcolor=black](0.5,0){.1}
\pscircle(0,0.5){.05}
\pscircle(0.5,0.5){.05}
\pscircle(1,0.5){.05}
\pscircle[fillcolor=black,fillstyle=solid](0.5,1){.1}
\psline(0.5,0)(0.5,1)
\rput(0.5,-0.4){$(0,0)$}

\pscircle[fillcolor=black](3,0){.1}
\pscircle(2.5,0.5){.05}
\pscircle(3,0.5){.05}
\pscircle[fillcolor=black,fillstyle=solid](3.5,0.5){.1}
\pscircle(3,1){.05}
\pscircle[fillcolor=black,fillstyle=solid](3.5,1){.1}
\psline(3,0)(3.5,0.5)(3.5,1)
\rput(3.0,-0.4){$(m,0)$}

\rput(1.5,-1.0){Case : $ m \geq 3 $}

\pscircle[fillcolor=black](6.5,0){.1}
\pscircle(6,0.5){.05}
\pscircle(6.5,0.5){.05}
\pscircle(7,0.5){.05}
\pscircle[fillcolor=black,fillstyle=solid](6.5,1){.1}
\psline(6.5,0)(6.5,1)
\rput[r](6.5,-0.4){$(0,0)$}

\pscircle[fillcolor=black](7.5,0){.1}
\pscircle(7.5,0.5){.05}
%\pscircle(8,0.5){.05}
\pscircle[fillcolor=black,fillstyle=solid](8,0.5){.1}
\pscircle(7.5,1){.05}
\pscircle[fillcolor=black,fillstyle=solid](8,1){.1}
\psline(7.5,0)(8,0.5)(8,1)
\rput[l](7.5,-0.4){$(m,0)$}

\rput(7,-1.0){Case : $ m =2 $}

\pscircle[fillcolor=black](12.5,0){.1}
\pscircle[fillcolor=black](13,0){.1}
\pscircle(12,0.5){.05}
\pscircle(12.5,0.5){.05}
\pscircle(13,0.5){.05}
\pscircle[fillcolor=black,fillstyle=solid](13.5,0.5){.1}
%\pscircle(0,1){.05}
\pscircle[fillcolor=black,fillstyle=solid](12.5,1){.1}
%\pscircle(1,1){.05}
\psline(12.5,0)(12.5,1)
\pscircle(13,1){.05}
\pscircle[fillcolor=black,fillstyle=solid](13.5,1){.1}
\psline(13,0)(13.5,0.5)(13.5,1)
\rput(12.0,-0.4){$(0,0)$}

\rput(13.5,-0.4){$(m,0)$}

\rput(12.5,-1.0){Case : $ m = 1 $}

\end{pspicture}

\caption{One possible realization of the event $\{ Z_{j+1} = m +1 \mid  Z_{j} = m \} $.
The bold vertices are open and all other vertices depicted are closed.}
\label{UpperBound}
%\end{center}
\end{figure}
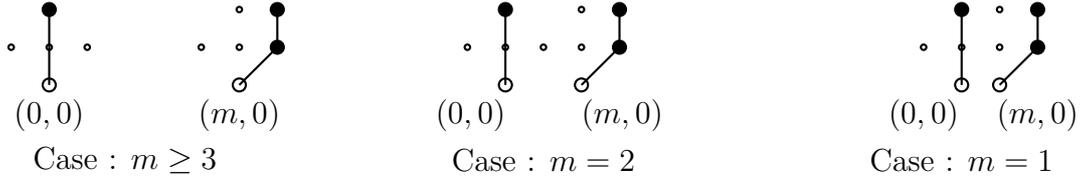

To show ${\cal Z}^{t_0}(t_0 +t)$ is locally finite,
we need to control the tail of the distribution of the coalescing time of two paths starting at the same instant
of time.
\begin{prop}
\label{prop:m_eqn}
For $\bu, \bv \in \Z^2$, $\bu(2) = \bv(2)$ consider the process,
$\{Z_j(\bu,\bv) : j \geq 0\}$. %,  defined in (\ref{def:Z_l}).
We have
\begin{equation*}
% \label{m_eqn}
\sup \{ \P(Z_{j+1}(\bu,\bv)=m |Z_{j}(\bu,\bv)=m) : m \geq 1 \} \leq \theta
\end{equation*}
for some $ \theta \in (0,1)$ which is independent of $j$.
\end{prop}
\noindent{\bf Proof :}
To prove Proposition \ref{prop:m_eqn}, we observe that for $ m \geq 3 $,
$ \P(Z_{j+1}(\bu,\bv)=m+1 |Z_{j}(\bu,\bv)=m) \geq (1-p)^{6}p^{3},
\P(Z_{j+1}(\bu,\bv)=3 |Z_{j}(\bu,\bv)=2) \geq (1-p)^{5}p^{3} $ and
$ \P(Z_{j+1}(\bu,\bv)= 2 |Z_{j}(\bu,\bv)=1) \geq (1-p)^{4}p^{3} $ (see
Figure \ref{UpperBound}).
Therefore, we have
\begin{align*}
\lefteqn{ \P(Z_{j+1}(\bu,\bv)=m|Z_{j}(\bu,\bv)=m) }\\
& \leq  1 - \P(Z_{1}(\bu,\bv)=m+1|Z_{0}(\bu,\bv)=m)\\
& \leq  1 - \min \{ (1-p)^{6}p^{3},  (1-p)^{5}p^{3},  (1-p)^{4}p^{3} \} = 1 - (1-p)^{6}p^{3} .
\tag*{\qed}
\end{align*}
% !TEX root = ./../DSFrf.tex
% !TEX program = XeLaTeX

Now, we prove an estimate on the tail of coalescing time.
We will use the following result (Theorem 4  of   Coletti {\it et al.} \cite{CFD09})%\/ (Theorem 4  of \cite{CFD09}).
\begin{theorem}
\label{thm:CFD09Thm04}
Let $ \{ V_j : j\geq 0\}$ be a Markov chain on the state space $ \{ 0, 1, \dotsc \} $ with
$0$ being the only absorbing state. Further,
assume that $ \{  V_j  : j \geq 0 \} $ is a martingale and $ \sup \{ \P(V_{1} = m | V_{0} = m ) :
m \geq 1 \} \leq \theta $  for some constant $ \theta \in (0,1)$. Let $\tau^{V} := \inf\{j \geq 1 : V_j = 0\}$.
Then, for some constant $ C_{25} $, we have
\begin{equation*}
\P (\tau^{V} \geq  n | V_0 = 1) \leq C_{25}/\sqrt{n} \text{ for all } n \geq 1.
 \end{equation*}
\end{theorem}

\begin{prop}
\label{lemma:CoalescingTimeTail}
Fix $\bu =(1,0),\; \bv =(0,0)\in \Z^2$, let $\nu =
\inf\{l:g_{\tau_l(\bu,\bv)}(\bu) =
g_{\tau_l(\bu,\bv)}(\bv)\}$, where $\tau_l(\bu,\bv)$ is the $l$ th regeneration
step as defined in
(\ref{r-tau}). For the $\nu$ th regeneration time $T_{\nu}(\bu,\bv)$ as
defined in (\ref{def:TimekthReg}), there exist positive constants $C_{26}$ and $C_{27}$, such that,
we have
\begin{equation*}
\P( \nu \geq t) \leq C_{26}/\sqrt{t} \quad \text{ and } \quad
\P(T_{\nu}(\bu,\bv) \geq t) \leq C_{27}/\sqrt{t}.
\end{equation*}
\end{prop}

%\begin{cor}
%\label{cor:CoalescingTimeTail}
%For $\bu =(a,0),\; \bv =(0,0)\in \Z^2$, let $ \nu (a) =
%\inf\{l:g_{\tau_l(\bu,\bv)}(\bu) = g_{\tau_l(\bu,\bv)}(\bv)\}$, where $ a > 0 $ is any positive integer. Then
%\begin{equation*}
%\P(\nu (a)  >t)\leq \frac{C_{21} a }{\sqrt{t}} \quad \text{ and } \quad
%\P(T_{\nu (a) }(\bu,\bv)>t)\leq \frac{C_{22} a}{\sqrt{t}}.
%\end{equation*}
%\end{cor}
%The corollary follows from Proposition \ref{lemma:CoalescingTimeTail},
%using the fact that the the coalescence time  $ \nu (a) $
%is given by the maximum of the coalescence times of the martingales
%starting from the pairs $ (i,0) $ and $ (i-1, 0) $ for $ i = 1, \dotsc, a $.

\noindent{\bf Proof }: The process $\{Z_j(\bu,\bv) : j \geq 0 \}$ satisfies the conditions
of Theorem \ref{thm:CFD09Thm04} and therefore, it follows that  $\P( \nu \geq t)
\leq C_{26}/\sqrt{t}$ where $ C_{26}$ is a constant.
% Consider any non-negative integer valued martingale $\{V_j : j\geq 0\}$
% such that for all $m \geq 1$,
% $\P(V_{j+1} = m | V_{j} = m ) \leq \pi $ for some constant $\pi \in (0,1)$ which is independent
% of $m$.
% n Theorem 4 Coletti {\it et al.}\/ \cite{CFD09}
% showed that there exists a positive constant $c_{25}$ independent of $n$ such that

To achieve the bound on $ T_{\nu}(\bu,\bv) $, we choose $ C_{28} = 1/ \bigl( 2 \E ( W^M ) \bigr) $
where $ W^M $ is as in Proposition \ref{prop:WidthTail}. Note that, it is  also the case that $ T_l (\bu,\bv)
\leq \sum_{i=1}^l W^M (i) $, for any $ l \geq 1 $, where $ \{  W^M (i) : i \geq 1 \} $ is an i.i.d. sequence,
each having the same distribution as that of $  W^M $ (see discussion before (\ref{def:TimekthReg}) and the
statement of  Proposition \ref{prop:WidthTail}). We have, % for all $ t > 2$,  %Using (\ref{eqn:NuBound}),
\begin{align*}
& \P(T_{\nu}(\bu,\bv) \geq t) \leq  \P (T_{\nu}(\bu,\bv) \geq t, \nu < C_{28} t )  + \P ( \nu \geq C_{28} t )  \\
& \leq \P ( T_{ \lfloor C_{28} t \rfloor } (\bu,\bv) \geq t )  + \frac{ C_{26} }{ \sqrt{ C_{28} t } }  \\
& \leq \P \bigl[  \sum_{i = 1}^{ \lfloor C_{28} t \rfloor } \bigl( W^M (i) - \E ( W^M (i) ) \bigr) \geq t  - \E ( W^M  )  \lfloor C_{28} t \rfloor  \bigr]
+ \frac{ C_{26} }{ \sqrt{ C_{28} t } }  \\
& \leq \frac{ \mathbb{V}\text{ar} \bigl(  \sum_{i = 1}^{ \lfloor C_{28} t \rfloor } W^M (i) \bigr) }{
  \bigl[ t  - \E ( W^M  )  \lfloor C_{28} t \rfloor  \bigr]^2}
+ \frac{ C_{26} }{ \sqrt{ C_{28} t } }  \\
&\leq  \frac{  \lfloor C_{28} t \rfloor \mathbb{V}\text{ar}  ( W^M  ) }
{ (t/2 - 1)^2} + \frac{ C_{26} }{ \sqrt{ C_{28} t } }  \\
& \leq  \frac{C_{27}}{\sqrt{t}}
\end{align*}
for a suitable choice of constant $ C_{27}$. This completes the proof. \qed

Before we proceed further, we introduce the following notation: for any $ A \subseteq \Z^2 $,
let $ {\cal X}^A $ and $ {\cal X}^A_n $ be the collection of paths starting at the vertices of $A$
and their scaled versions respectively, i.e., $ {\cal X}^A := \{ \pi^{\bu} : \bu \in A \} $ and 
$ {\cal X}^A_n := \{ \pi^{\bu}_n  : \bu \in A \}$.
% For any $ t $, denote by $ \chi^A (t) = \{ \pi^{\bu} (t) : \bu \in A \} $
% and $ \chi^A_n (t) = \{ \pi^{\bu}_n (t) : \bu \in A \} $.
Now, we prove the following proposition which is an adaptation of Lemma 2.7 of Newman {\it et al.}\/ \cite{NRS05}.
\begin{prop}
\label{prop:DensityEstimate}
For $ a,b \in \R $ and $ t > 0 $, for all $ n \geq 1 $, we have
\begin{equation}
\label{eqn:estDensity}
\E \bigl[  \# \bigl( {\cal X}^{ \Z \times \{0\}}_n (t) \cap ( [a,b] \times \{t\})  \bigr) \bigr] \leq \frac{C_{29}(b-a)}{\sqrt{t}}
\end{equation}
where $ C_{29}  $ is a constant, independent of $t$ and  $ a, b $.
\end{prop}

\noindent{\bf Proof }: Fix $ s > 0 $ and let $ u(s) = \E \bigl[  \# \bigl( {\cal X}^{ \Z \times \{0\}} (s) \cap ([0,1)  \times \{s\})  \bigr) \bigr] $.
Set $ M \geq 1 $. Then by translation invariance, we have $  \E \bigl[  \# \bigl( {\cal X}^{ \Z \times \{0\}} (s) \cap ([0,M)  \times \{s\})  \bigr) \bigr]
= M u(s) $.
% $a_n = \lfloor an \sigma \rfloor $ and
% $ b_n = \lfloor bn \sigma \rfloor + 1$ and set $ l_n = b_n - a_n $. Now, scaling back and using translation invariance
% of our model,
Now, we have
\begin{align*}
  M u(s) 
%=  \E \bigl[  \# \bigl( \chi^{ \Z \times \{0\}} (s) \cap ([0,M)  \times \{s\})   \bigr) \bigr] \\
& \leq \sum_{ i = -\infty}^{\infty} \E \bigl[  \# \bigl( {\cal X}^{ [i M , (i+1)M)  \times \{0\}} (s)
\cap ([0,M)  \times \{s\})  \bigr) \bigr] \\
& = \sum_{ i = -\infty}^{\infty} \E \bigl[  \# \bigl( {\cal X}^{ [0 , M)  \times \{0\}} (s)
\cap\bigl( [ -i M , -(i-1)M)  \bigr)  \times \{s\} \bigr)  \bigr] \\
& =  \E \bigl[  \# \bigl( {\cal X}^{ [0 , M)  \times \{0\}} (s)
  \bigr) \bigr] \\
& \leq M - (M - 1) \bigl[  1 - \P \bigl( \pi^{(0,0)} (s) \neq  \pi^{(1,0)} (s) \bigr) \bigr]
\\ & \leq  M - (M - 1) \bigl[  1 - \frac{C_{27}}{ \sqrt{ s }} \bigr] \qquad (\text{from Proposition \ref{lemma:CoalescingTimeTail}}
)\\
& \leq 1 + \frac{ C_{27} (M-1)}{  \sqrt{ s } }.
\end{align*}
Dividing both sides by $M$ and letting $ M \to \infty$, we have $ u(s) \leq C_{27}/\sqrt{s} $.

For any $ n \geq 1$, let $ l_n = \lfloor n (b-a)/2 \rfloor + 1. $
Now, we have  $ \E \bigl[  \# \bigl( {\cal X}^{ \Z \times \{0\}}_n (t) \cap ([a,b]  \times \{t\})  \bigr) \bigr]
= \E \bigl[  \# \bigl( {\cal X}^{ \Z \times \{0\}}_n (t) \cap ( [-(b-a)/2, (b-a)/2]   \times \{t\})  \bigr) \bigr]
\leq \E \bigl[  \# \bigl( {\cal X}^{ \Z \times \{0\}} (n^2 \gamma t) \cap ( [-l_n, l_n)  \times \{n^2 \gamma t \})  \bigr) \bigr]
= 2 l_n u (n^2 \gamma t ) \leq 2  C_{27}  l_n / \sqrt{ n^2 \gamma t } \leq C_{29}(b-a)/\sqrt{t}$
for a proper choice of $C_{29}$. \qed
%This proves the proposition.

Let $({\cal P}, \rho_{{\cal P}})$ be the space of compact subsets of $(\R^2_c, \rho)$ with the induced
Hausdorff metric. Since $ \overline{\cal X}_{n}^{t_0}  $ converges weakly to $ {\cal Z}^{t_0} $, by the continuous
mapping theorem, we have that $  \overline{\cal X}_{n}^{t_0}  (t_0 + t) $ converges weakly
to $ {\cal Z}^{t_0}(t_0 + t)$ in $({\cal P}, \rho_{{\cal P}})$.  Next we prove that ${\cal Z}^{t_0}(t_0 + t)$
is a.s. locally finite. 

\begin{prop}
\label{prop:LocalFinite}
For any $t > 0$, ${\cal Z}^{t_0}(t_0 + t)$ is a.s. locally finite and
\begin{equation*}
\E \bigl[  \# \bigl( {\cal Z}^{t_0}(t_0 + t) \cap ( (a,b)  \times \{t_0 + t\})  \bigr) \bigr] \leq \frac{C_{29}(b-a)}{\sqrt{t}}
\end{equation*}
for  $ C_{29} $ as in the previous proposition.
\end{prop}

\noindent{\bf Proof }: For the first part it is enough to consider $ t_0 = 0 $ and 
prove that $  \#\bigl(  {\cal Z}^{0}( t) \cap ( (-m, m)  \times \{t\})  \bigr)$ is
finite a.s. for any $ m \geq 1 $.  First, we observe that if $ \bu = ( u(1), u(2) )$ is such that $ u(2) < 0 $ and $ h(\bu) (2) > 0 $,
then it must be the case that $ h ( u(1), 0) = h( \bu)$. Therefore,  $ \#\bigl( \overline{\cal X}_{n}^{0}  (t) \cap ( (-m, m)  \times \{t\})   \bigr)
\leq  \#\bigl(  {\cal X}^{ \Z \times \{0\}}_n (t)  \cap ( (-m, m)  \times \{t\})   \bigr)   $.
From  Proposition \ref{prop:DensityEstimate}, 
\begin{align*}
& \E \bigl[  \liminf_{n \to \infty} \#\bigl( \overline{\cal X}_{n}^{0}  (t)  \cap ( (-m, m)  \times \{t\})    \bigr) \bigr]\\
 & \leq    \E \bigl[  \liminf_{n \to \infty}   \#\bigl(  {\cal X}^{ \Z \times \{0\}}_n (t)  \cap ( (-m, m)  \times \{t\})   
 \bigr) \bigr]\\
& \leq \liminf_{ n \to \infty}  
\E \bigl[     \#\bigl(  {\cal X}^{ \Z \times \{0\}}_n (t)  \cap  ( (-m, m)  \times \{t\})   \bigr) \bigr] \\ 
& \leq 2C_{29}m /\sqrt{t}
.
\end{align*}
Therefore, we conclude that  $   \liminf_{n \to \infty} \#\bigl( \overline{\cal X}_{n}^{0}  (t)  \cap  ( (-m, m)  \times \{t\})  \bigr) < \infty $ almost surely.

Since $  \overline{\cal X}_{n}^{0}  (t) \Rightarrow {\cal Z}^{0}( t)$ in $({\cal P}, \rho_{{\cal P}})$,
using Skorohod's representation theorem, we may couple the processes so that the above convergence is almost sure.
We now claim that, almost surely,
\begin{equation*}
 \#\bigl( {\cal Z}^{0}(t) \cap  ( (-m, m)  \times \{t\})   \bigr) \leq  \liminf_{n \to \infty} \#\bigl( \overline{\cal X}_{n}^{0}  (t)
\cap ( (-m, m)  \times \{t\})  \bigr).
\end{equation*}
Fix $ \omega $ such that $ \liminf_{n \to \infty} \#\bigl( \overline{\cal X}_{n}^{0}  (t)  \cap ( (-m, m)  \times \{t\})   \bigr)
= l (\omega) < \infty $.  For this $\omega$ we may choose a subsequence $ n_k$ along which $  \#\bigl( \overline{\cal X}_{n_k}^{0}  (t)  \cap  ( (-m, m)  \times \{t\})   \bigr)
=  l (\omega) $ for all large $ k $. If $  {\cal Z}^{0}(t) \cap ( (-m, m)  \times \{t\})  $ has at least $ l (\omega) + 1 $ distinct points, say, $
\{ (y_i, t) : i = 1, \dotsc, l(\omega) + 1 \} $, we may choose $ \delta > 0 $ so small that the intervals $ (y_i - \delta, y_i + \delta) \subseteq (-m, m)$,
for $ i = 1, 2,  \dotsc, l(\omega) + 1 $ are mutually disjoint. Since $ \overline{\cal X}_{n_k}^{0}  (t)  \cap  ( [-m, m]  \times \{t\})    $
converges to $  {\cal Z}^{0}(t) \cap ( [-m, m]  \times \{t\}) $ in $ ({\cal P}, \rho_{{\cal P}})$, each of these intervals should contain
at least one point of $ \overline{\cal X}_{n_k}^{0}  (t)  \cap  ( (-m, m)  \times \{t\})   $  for all large $k$. This is a contradiction as
for all large $ k$, $ \overline{\cal X}_{n_k}^{0}  (t)  \cap ( (-m, m)  \times \{t\})   $ has exactly $ l(\omega)$ many points.

The expectation bound now follows from the above inequality and completes the proof. \qed

Since $ \E \bigl[  \# \bigl( {\cal Z}^{t_0}(t_0 + t) \cap( (x - \delta, x+\delta)  \times \{ t_0 + t\} ) \bigr) \bigr] \leq 2C_{29}\delta/\sqrt{t}
\to 0$, as $ \delta \to 0$,  we may conclude
\begin{cor}
\label{cor:ZHavingAGivenPoint}
For any $ x \in \R $, $ \P \bigl( (x, t_0 + t)  \in {\cal Z}^{t_0}(t_0 + t) \bigr) = 0 $.
\end{cor}

We now state the main proposition of this subsection which is
similar to Lemma  6.3 of Newman {\it et al.}\/ \cite{NRS05}.
\begin{prop}
\label{prop:randomBM} For $ t > 0$, we have
\begin{equation*}
{\cal Z}^{t_0; (t_0 + t)_T} \stackrel{d}{=} {\cal W}^{{\cal Z}^{t_0}(t_0 +t)}
\end{equation*}
where ${\cal W}^{{\cal Z}^{t_0}(t_0 +t)}$ is the set of paths given by
coalescing Brownian motions starting
from a random point set distributed as ${\cal Z}^{t_0}(t_0 +t)$.
\end{prop}
We first complete the proof of $(E_1^{\prime})$ assuming the validity of the above proposition.
For  $ 0 < \nu  < t $, we have,
\begin{align*}
\E[\hat{\eta}_{{\cal Z}^{t_0}}(t_0,t;a,b)]  & =
\E\bigl [\hat\eta_{{\cal Z}^{t_0;(t_0+\nu)_T}}(t_0+\nu,t-\nu;a,b)\bigr] \\
 & \leq \E[\hat\eta_{{\cal W}}(t_0+\nu,t-\nu;a,b)] =\frac{b-a}{\sqrt{\pi (t-\nu)}}.
\end{align*}
Letting $ \nu \to 0 $ we obtain  $(E_1^{\prime})$.

Before proving Proposition \ref{prop:randomBM} we observe that the paths 
in ${\cal X}^{0; (t)_T}_n$ carry their own history region along with them.
Hence Proposition \ref{prop:fddbweb2} cannot be applied directly to obtain the finite dimensional 
distributions of ${\cal Z}^{0; (t)_T}$.

\noindent {\bf Proof of Proposition \ref{prop:randomBM}:}
% Before embarking on the formal proof, we explain our method briefly. 
It suffices to prove the result for $ t_0 = 0$.
We restrict our attention to the set of paths
which start in  $ [-m,m] \times\{t\}$ for some $m \geq 1$. Let us
denote by $ {\cal Z}^{0; (t)_T}_m = \{ \pi \in {\cal Z}^{0; (t)_T} :
\pi (t) \in [-m, m] \} $ and $ {\cal W}^{{\cal Z}^{t}_m} = {\cal W}^{{\cal Z}^{0}(t) \cap \bigl( [-m,m] \times \{t\} \bigr)}  $.
Observe that it is enough to show  $  {\cal Z}^{0; (t)_T}_m \stackrel{d}{=} {\cal W}^{{\cal Z}^{t}_m}$
for any $m \geq 1$.
For the rest of the section we fix $ m \geq 1$. 

Consider the mapping $ g : {\cal H} \to {\cal H} $ given by $ g (K) = \{ \pi \in K^{0; (t)_{T}}
: \pi (t) \in [-m, m] \} $. Using Corollary \ref{cor:ZHavingAGivenPoint} we have
$ \P ( {\cal Z}^{0} \in D_g ) = 0 $ where $ D_g $ is the set of discontinuity points of the map $ g $.
Since $  \overline{\cal X}_{n}^{0}
\Rightarrow {\cal Z}^{0}$, from Theorem 5.1 of Billingsley \cite{B68} we  have $ g (\overline{\cal X}_{n}^{0})
\Rightarrow  g({\cal Z}^{0}) = {\cal Z}^{0; (t)_T}_m $.

% In our case, we are going to take $ D = {\cal Z}(t)
% \cap \bigl( [-m, m] \times \{t\} \bigr) $. 
Now we will obtain ${\cal W}^{{\cal Z}^{t}_m}$ 
as a weak limit of $g (\overline{\cal X}_{n}^{0})$ to establish the required equality.
For any $ K \in {\cal H}$, we consider the map $ f : {\cal H} \to
{\cal P} $ given by $ f(K) = \{ (\pi (t) , t) : \sigma_{\pi}\leq t,   \pi ( t) \in [-m,m] \} $.
Again, using Corollary \ref{cor:ZHavingAGivenPoint}, we observe that $ \P ( {\cal Z}^{0} \in D_f ) = 0 $ 
where $ D_f $ is the set of all discontinuity points of the map.
For $ t > 0$, taking $ D := {\cal Z}^{0}(t)
 \cap \bigl( [-m, m] \times \{t\} \bigr) $, we have that, as $n \to \infty$,
\begin{equation}
\label{eqn:WeakConvXnToD}
 \overline{\cal X}_{n}^{0}  (t )
\cap \bigl( [-m, m] \times \{t\} \bigr) \Rightarrow D  \text{ in } ({\cal P}, \rho_{{\cal P}}).
\end{equation}

Fix $ n \geq 1, \beta  < 1/2 $ and $ 2 \beta < \delta < 1$. Define the sets
$$ D_n  := \bigl\{   \bigl(\lfloor n \sigma \bx(1)\rfloor ,  \lfloor n^2 \gamma t \rfloor
 + \lfloor n^{\beta}\rfloor + 1 \bigr) :  (\bx (1), \bx(2) )  \in {\cal X}^0_n(t), \bx (1) \in [-m, m] \bigr\}
 $$
 and $D_n^{(\text{scaled})}  := \{(\by(1)/(n \sigma), \by(2)/(n^2\gamma)): \by \in D_n\}$.
 Note that $D_n$ (and hence $D_n^{(\text{scaled})}$) is a finite set.
For each $ \by = (\by(1), \by(2)) \in D_n^{(\text{scaled})} $, we have $ \bx \in {\cal X}^0_n(t)
\cap \bigl( [-m, m] \times \{t\}  \bigr) $ such that
$ || \bx - \by ||_2 \leq   (n^{\beta} + 2)/n $ and vice versa. Thus, %from the definition of Hausdorff metric, we have
$ \rho_{{\cal P}} ( {\cal X}^0_n(t) \cap \bigl( [-m, m] \times \{t\} \bigr)  , D_n^{(\text{scaled})})
\to 0 $ almost surely. Using (\ref{eqn:WeakConvXnToD}),  we conclude
that $ D_n^{(\text{scaled})} \Rightarrow {\cal Z}^{0}( t) \cap \bigl( [-m, m] \times \{t\} \bigr) = D  $
 in $({\cal P}, \rho_{{\cal P}})$.

% ************************************
% Using $\overline{\cal X}_{n}^{0}  ( t )$, we construct a finite set $ D_n \subseteq \Z^2$
% such that $  D_n^{(\text{scaled})} \Rightarrow D$ in $({\cal P}, \rho_{{\cal P}})$.
% % Note that we cannot apply Lemma \ref{lem:randomBMstart} directly here, because
% % the paths in $\overline{\cal X}^{0}_n$ may carry history when crossing the line $y=t$.
% The construction ensures that on a set whose probability converges to $1$,  the evolution of paths
% from $ D_n$ is independent of history so that an application  of Lemma \ref{lem:randomBMstart} below
% implies that $ \chi_n^{D_n} \Rightarrow {\cal W}^D =  {\cal W}^{{\cal Z}^{t}_m} $ in $ ({\cal H}, d_{{\cal H}}) $.
% ************************************

%
% For any compact set $F =  [a_1, b_1] \times [a_2, b_2] $ ($a_1 < b_1, a_2 < b_2)$,
% the Hausdorff distance between  the compact sets $ D_n^{(\text{scaled})} \cap F$ and
% $ \bigl[ {\cal X}^0_n(t) \cap \bigl( [-m, m] \times \{t\} \bigr) \bigr]  \cap F$, under the usual Euclidean metric, is at most
% $  [1/(n \sigma) + (n^{\beta} + 1)/(n^2 \gamma)]$ which converges to $0$
%as $ n \to \infty $. Therefore, it implies that

We show now that it is unlikely that a path $ \pi \in  {\cal X}^{0}_n $, which crosses the x-axis
far from the origin, will cross the line $ \{ y=t\}$ inside $ [-m,m]$.
 Consider the event
\begin{align*}
E_n := \{ \text{there exists } \pi \in  {\cal X}^{0}_n \text{ with } \pi(0)
\not\in [-n, n] \text{ and }\pi(t) \in [-m, m] \}.
\end{align*}
Scaling back to the original lattice and using the non-crossing property of paths, we observe that the if the paths starting from
$ \bu=  (- \lfloor n^2 \sigma \rfloor ,0) $ and $ \bv = ( \lfloor n^2 \sigma \rfloor , 0) $ do not cross the segment
$ [-n\sigma m, n\sigma m] \times \{n^2\gamma t\}$, then the paths
which cross the x-axis to the
left of $ -n^2 \sigma $ or to the right of $ n^2 \sigma$ will also stay away from that segment.
% Taking $\pi^{(-n^3, 0)}$ and $\pi^{(n^3, 0)}$ as the paths starting from $(-n^3, 0)$
% and $(n^3, 0)$ respectively,
Hence, we have
%  we observe that if $\pi \in  E_n$, the non-crossing nature of the paths
% ensure that either (a) $\pi(s) \leq \pi^{(-n^3, 0)}$ for all $s \geq 0$, or (b) $\pi(s) \geq \pi^{(n^3, 0)}$ for all $s \geq 0$.
% Hence, we have
\begin{align*}
 \P(E_n) &
\leq  \P(\{\pi^{\bu}(n^2 \gamma t) \geq -n\sigma m\}
\cup  \{\pi^{\bv}(n^2 \gamma t) \leq n\sigma m\})\\
& \leq  2\P(\pi^{\bu}(n^2 \gamma t) \geq -n\sigma m)\\
&=  2\P(\pi^{\mathbf{0}}(n^2 \gamma t) \geq \lfloor n^2 \sigma \rfloor  -n\sigma m) \\
& = 2\P \bigl( \pi^{\mathbf{0}}_n (t) \geq (\lfloor n^2 \sigma \rfloor  -n\sigma m) /(n \sigma)  \bigr) \\
&\to  0 \text{  as $n \to \infty $, }
\end{align*}
   $ \pi^{\mathbf{0}} $ being the path starting at $(0,0)$, and the last step follows from Proposition \ref{lemma:Bmotion}.

Let
$
\tilde{\cal X}_{n} := \{ \pi \in \overline{\cal X}_{n}^{0} : \pi (0) \in [-n,n] \}.
$
On the event $ E_n^c $, we observe that $ g ( \tilde{\cal X}_{n} ) =  g ( \overline{\cal X}_{n}^{0} ) $ as
 $ ({\cal H}, d_{{\cal H}}) $-valued random variable. Hence, we have
\begin{equation}
\label{eqn:DistHConv0ForStart}
d_{{\cal H}} \bigl( g ( \tilde{\cal X}_{n} ) , g ( \overline{\cal X}_{n}^{0} ) \bigr) \prob 0.
\end{equation}

% needs at most $O(n^2)$  steps to go above the line
% $y = n^2\gamma_0\epsilon$, while the probability of an increment $||h^j(n^3, 0) - h^{j+1}(n^3, 0)||_1$ being larger than $n$
% decays exponentially in $n$.

% Fix $0<\beta < 1/2$ and choose $2\beta < \delta < 1$. We define the point set
% \begin{align}
% \label{PtsetDnm}
% D^{n,m} & := \{\bx^r_n, \bx^l_n : \bx \in {\cal X}^0_n(\epsilon) \cap (-m,m)\times \{\epsilon\}\}\text{ where}\\
% \bx^r_n & := (\lfloor n \sigma_0\bx(1)\rfloor + \lfloor n^{\delta}\rfloor, \lfloor n^2 \gamma_0\epsilon \rfloor + \lfloor m^{\beta}\rfloor + 1)\text{ and }\nonumber\\
% \bx^l_n & := (\lfloor n \sigma_0\bx(1)\rfloor - \lfloor n^{\delta}\rfloor, \lfloor n^2 \gamma_0\epsilon \rfloor + \lfloor m^{\beta}\rfloor + 1)\nonumber.
% \end{align}

Now, we follow the paths in $  g (\tilde{\cal X}_{n}) $ until they cross the line $ \{ y =  t \} $
 and consider the  history that is created in doing so. Scaling back to the original lattice,
for $\bu \in \Z^2$ with $\bu(2) \leq 0$ let $l^{\bu} := \min\{j: h^j (\bu)(2) \geq n^2\gamma t\}$
denote the number of steps taken by the path starting from $\bu$ to cross the line $\{y = n^2\gamma t\}$.
We define, the set of explored regions for paths in unscaled version of $ g (\tilde{\cal X}_{n})$,  by
\begin{equation*}
\tilde{\Delta}_n := \bigcup
\bigcup^{l^{\bu}-1}_{i = 0} S^{+}(h^i(\bu), ||h^i(\bu) - h^{i+1}(\bu)||_1)
\end{equation*}
where the first union is over $\bu \in V, \bu(2)\leq 0, \pi^{\bu}_n(0) \in [-n, n], \pi^{\bu}_n( t) \in [-m, m]$.
Consider the event $ F_n = \bigl\{ \tilde{\Delta}_n  \nsubseteq  \mathbb{H} (\lfloor n^2 \gamma t \rfloor
+ \lfloor n^{\beta}\rfloor) \bigr\}$.
Assuming $ \P (F_n)  \to 0 $ as $ n \to \infty $ (which will be shown shortly), we observe that   on the event $ (E_n \cup F_n)^c$,
using the fact that the vertices of $ D_n $ lie on the line  $ \{ y = \lfloor n^2 \gamma  t \rfloor + \lfloor n^{\beta}\rfloor + 1 \}$,
the evolution of the paths from the set $ D_n$ is independent of the history $ \tilde{\Delta}_n $. This allows us to 
adapt Lemma 6.5 of \cite{NRS05} for our model  and conclude
\begin{equation}
\label{Kumarjitdekho}
{\cal X}_n^{D_n} \Rightarrow {\cal W}^D =  {\cal W}^{{\cal Z}^{ t}_m} \text{  in } ({\cal H}, d_{{\cal H}})  .
\end{equation}
The details of this argument is presented 
%in Lemma \ref{lem:Kumarjitdekho} 
in the appendix.

% 
% 
% \begin{lemma}
% \label{lem:randomBMstart}
% Let $D_n$ and $D$ be $({\cal P}, \rho_{{\cal P}})$ valued random variables defined as above.
% Suppose the evolution of the paths from the set $ D_n$ is independent of the history $ \tilde{\Delta}_n $, where 
% $ \tilde{\Delta}_n $ is as defined above. Then 
% \begin{equation*}
% \chi_n^{D_n} \Rightarrow {\cal W}^D \text{ in } ({\cal H}, d_{{\cal H}})
% \end{equation*}
% where $ {\cal W}^D $ is the 
% coalescing Brownian motions process starting from a random point set distributed as $D$.
% \end{lemma}

% 
% Also suppose $D_n$ is such that the collection of scaled 
% paths starting from $D_n$ is independent of $\sigma(D_n)$. More formally, let $\{U^a_{\bw}:
% \bw \in \Z^2\}$ be a  of i.i.d. $U[0,1]$ random variables, independent of the collection
% $\{U_{\bw}:\bw \in \Z^2\}$ used to construct $\chi_n$. For  $A \subseteq \Z$, let $\chi_{a,n}^A$
% be the collection of scaled paths starting from $A$ constructed using $\{U^a_{\bw}:
% \bw \in \Z^2\}$ only. Suppose for all $n$ large, $D_n$ is such that 
% $$
% \chi_n^{D_n} | \sigma(D_n) \stackrel{d}{=} \chi_{a,n}^{D_n} \text{ almost surely}.
% $$
% Then, the scaled paths starting from
% $ D_n $ converges weakly to $ {\cal W}^D $, i.e.,
% \begin{equation*}
% \chi_n^{D_n} \Rightarrow {\cal W}^D \text{ in } ({\cal H}, d_{{\cal H}})
% \end{equation*}
% where $ {\cal W}^D $ is the 
% coalescing Brownian motions process starting from a random point set distributed as $D$.
% \end{lemma}

To show $ \P (F_n) \to 0 $, we consider paths starting from each $(j,0)$, $ - n^2 \sigma -1 \leq j \leq  n^2 \sigma +1$ and
consider the region explored by  these paths until they cross the line $ \{ y = n^2 \gamma t \} $, i.e., set
\begin{equation*}
\Delta_n^{\prime} := \bigcup_{j =  - \lfloor  n^2 \sigma \rfloor - 1}^{  \lfloor  n^2 \sigma \rfloor + 1 }
\bigcup^{l^{(j,0)}-1}_{i = 0} S^{+} \bigl(h^i(j,0), ||h^i(j,0) - h^{i+1}(j,0)||_1 \bigr) .
\end{equation*}
We observe that on the event $ E_n^c $, we must have $ \tilde{\Delta}_n \setminus \mathbb{H}(0) \subseteq \Delta_n^{\prime} $.
Therefore, we have
%
% be the union of the explored region of each of these path $\pi^{(j,0)}$ in the first
% $l^{(j,0)}$ steps  (i.e. until the step after which crosses the line $y = n^2\gamma_0 t$).
% With $E_n$ defined as above, we have
\begin{align*}
\lefteqn{ \P \bigl( \tilde{\Delta}_n  \nsubseteq \mathbb{H} (\lfloor n^2 \gamma t
\rfloor + \lfloor  n^{\beta}\rfloor ) \bigr) } \\
 & \leq  \P \bigl(  \bigl\{\Delta_n^{\prime} \nsubseteq \mathbb{H} (\lfloor n^2 \gamma t
\rfloor + \lfloor  n^{\beta}\rfloor) \bigl\} \cap E_n^c \bigr) + \P(E_n) \\
& \leq  \P \Bigl[ \bigcup_{j = -  \lfloor  n^2 \sigma \rfloor - 1}^{ + \lfloor  n^2 \sigma \rfloor + 1}
\bigcup_{i = 1}^{\lfloor n^2\gamma t \rfloor + 1} \bigl\{ W_i(j,0) \geq n^{\beta} \bigr\} \Bigr]  + \P(E_n)\\
&\leq  (2\lfloor n^2\sigma \rfloor + 3)( \lfloor n^2\gamma t \rfloor + 1)\P(W_1(0,0) \geq n^{\beta}) + \P(E_n)\\
& \to 0
\text{ as }n \to \infty,
\end{align*}
where the penultimate inequality follows from the fact that each path $\pi^{(j,0)}$ can have at most
$\lfloor n^2\gamma t\rfloor + 1$ many regenerations until it crosses the line
$\{y = n^2\gamma t\}$  and the last step follows from Proposition \ref{prop:WidthTail}.

Finally, we show that $  d_{{\cal H}} \bigl( {\cal X}_n^{D_n} ,
g (\overline{\cal X}_{n}^{0}) \bigr)  \prob 0 $ to complete the proof of  Proposition \ref{prop:randomBM}.
Consider the event that one of the paths in  $  g (\tilde{\cal X}_{n}) $ moves significantly
far in a short period after crossing the line $ \{ y =  t \} $. Define the event
\begin{align*}
G_n := & \{\text{there exists } \pi  \in  g (\tilde{\cal X}_{n}) \text{ with } |\pi(t) - \pi( s)|
> n^{\delta -1}/\sigma \\
& \qquad  \text{ for some } s \in [ t,  t + (n^{\beta - 2}/\gamma)]\}.
\end{align*}
We have $ \# \bigl( g (\tilde{\cal X}_{n}) \bigr) \leq  \#\{\pi(t) \in [-m,m]: \pi \in {\cal X}_n^{\Z \times \{0\}}\}$ 
and hence, on the event $E_n^c$,  $ \# \bigl( g (\tilde{\cal X}_{n}) \bigr) \leq 2\lfloor n^2\sigma \rfloor + 3$.
From Proposition \ref{prop:HistPathStayClose} it follows that $\P(G_n) \leq \P (E_n ) + \P( F_n) +
(2\lfloor n^2\sigma \rfloor + 3) \P \bigl( \sup \{ | \pi^{\bf{0}} (s)  |  : 0 \leq s \leq
n^{\beta} \} \geq  n^{\delta} \bigl|  \{  U_{\bw} : \bw \in  \Delta_n^{\prime}\} \cap (F_n)^c  \bigr)
\leq \P (E_n ) + \P( F_n)  + C_{23} n^{\beta} \exp ( - C_{24} n^{\beta} ) (2\lfloor n^2\sigma \rfloor + 3)
\to 0 $ as $ n \to \infty $.

The finiteness of $g (\overline{\cal X}_{n}^{0}) $ allows us to enumerate the paths in $g (\overline{\cal X}_{n}^{0}) $
as $\pi^1, \dotsc , \pi^N$ for some random $N$. Let $ \bx_j  :=   
\bigl( \lfloor n \sigma \pi^j ( t) \rfloor ,  \lfloor n^2 \gamma  t \rfloor
 + \lfloor n^{\beta}\rfloor + 1 \bigr) $ for  $j = 1, \dotsc , N$. Note that $\bx_j$'s need not be distinct, 
however $D_n = \{\bx_i : 1 \leq i \leq N\}$ and hence 
$  d_{{\cal H}} \bigl( {\cal X}_n^{D_n} , g (\overline{\cal X}_{n}^{0}) \bigr)  \leq \max \{ d_{\Pi} ( \pi^j , \pi^{\bx_j}_n)
: 1 \leq j \leq N  \} $. In other words, for $s > t$ taking $M_n^j(s) := \sup \{ | \pi^j(l) -  \pi^{\bx_j}_n (l) | : 
l \in  [ t + (\lfloor n^{\beta} \rfloor + 1)/(n^2 \gamma),  s ] \}$ we need to show
(a) $ \max \{ M_n^j (s )  : 1 \leq j \leq N \} \prob 0 $ and (b) $\max \{|| (\pi^j(\sigma_{\pi^j}),\sigma_{\pi_j}) - 
(\pi^{\bx_j}_n(\sigma_{ \pi^{\bx_j}_n}),\sigma_{ \pi^{\bx_j}_n}) ||_2: 1 \leq j \leq N\} \prob 0$ as $n \to \infty$.

Since $|| (\pi^j(\sigma_{\pi^j}),\sigma_{\pi_j}) - (\pi^{\bx_j}_n(\sigma_{ \pi^{\bx_j}_n}),
\sigma_{ \pi^{\bx_j}_n})   ||_2 \leq (n^{\beta} + 2)/n$, (b) follows immediately.

Clearly,   $ \E(N) = \E\bigl[ \# \bigl(g (\overline{\cal X}_{n}^{0} ) \bigr) \bigr] \leq
\E \bigl[ \# \bigl( {\cal X}^{ \Z \times \{0\}}_n ( t)
\cap ( [-m,m] \times \{ t\})  \bigr)  \bigr] \leq 2C_{29}m/ t $ from Proposition \ref{prop:DensityEstimate}. 
So, given $ \eta, \eta^{\prime} > 0$,
 we can choose $ L $ (independent of $n$) so large that $ \P (N
 \geq L ) \leq \eta^{\prime}/2 $.
 On the event
$G_n^c$, we observe that  the paths $  \pi^j $ and $ \pi^{\bx_j}_n  $
will lie between the scaled paths starting from $ \bigl( \lfloor n \sigma \pi^j ( t) \rfloor - \lfloor n^{\delta}
\rfloor,  \lfloor n^2 \gamma  t \rfloor+ \lfloor n^{\beta}\rfloor + 1 \bigr)$ and $ \bigl( \lfloor n \sigma \pi^j ( t)
 \rfloor + \lfloor n^{\delta} \rfloor,  \lfloor n^2 \gamma  t \rfloor+ \lfloor n^{\beta}\rfloor + 1 \bigr)$.
Therefore, by translation invariance,
$ \P (  \{M_j (s )  \geq \eta\} \cap F_n^c )  \leq \P ( \sup \{ | \pi^{\mathbf{0}}_n (l)  - \pi^{( 2\lfloor n^{\delta}
\rfloor, 0)}_n (l) | : 0 \leq l \leq s -  t \} \geq \eta )$.
Hence, we have that 
$\P(\max \{ M_n^j (s )  : 1 \leq j \leq N \} \geq \eta)
 \leq  \P( F_n ) + 
\P( G_n)
+ \P  \bigl[ \# \bigl(g (\overline{\cal X}_{n}^{0} ) \bigr)
 \geq L \bigr]  + L \P ( \sup \{ | \pi^{\mathbf{0}}_n (l) - \pi^{( 2\lfloor n^{\delta}
\rfloor, 0)}_n (l) | : 0 \leq l \leq s -  t \} \geq \eta ) $.
 By Proposition \ref{prop:sandwich:DistConvNoHist},  both the paths $  \pi^{\mathbf{0}}_n $ and
 $ \pi^{( 2\lfloor n^{\delta}
\rfloor, 0)}_n $  converge to the same Brownian motion.
Therefore, for all large $ n $, we have $ \P (\max \{ M_n^j (s )  : 1 \leq j \leq N \} \geq \eta)
  < \eta^{\prime}$. This completes the proof. \qed

%\input{sec6rf.tex}
% !TEX root = ./DSFrf.tex
% !TEX program = XeLaTeX

\section{Appendix}

\noindent{\bf Proof of Lemma \ref{lem:MCAuxRes}:} It suffices to prove that,  for some
$ \alpha > 0 $, we have $ \E ( \exp ( \alpha \tau_M ) ) < \infty $.
Since $ M_{n+1} $ is a function of $M_n $ and an independent sequence of random
variables,
$ \{ M_n : n \geq 0 \}$ is a Markov chain.  Furthermore, it is irreducible and
recurrent. Using Proposition 5.5, Chapter 1
of Asmussen \cite{A03}, it suffices to show that there exist a nonnegative function
$f:\N\cup\{0\} \to \R^{+}$, $n_0 \in \N$ and $r > 1$ such that
 $f(j) > \nu$ for some $\nu > 0$ and $\E[f(M_{1})| M_0 = j] < \infty$ for all $j \leq n_0$,
while for $j > n_0$, 
$\E[f(M_{1})| M_0 = j] \leq f(j)/r$.

Taking $f : \{ 0, 1, 2 , \dotsc\}\to \R $ to be $ f (i) = \exp ( \alpha i ) $,
where $\alpha > 0$ is small enough so that $\E[\exp(\alpha\theta_1)] < \infty$
and $\exp(-\alpha) < 1/r$, we have
\begin{align*}
\lefteqn{ \E[\exp(\alpha(M_{1} - M_{0}))|M_0 = m] }\\ & = \exp(-\alpha) \P(\theta_{1} \leq m) 
 + \exp(-\alpha(m + 1))\E[{\mathbf 1} ( \theta_{1} > m )\exp(\alpha\theta_{1})]\\
 & < (1/r) + \exp(-\alpha(m + 1))\E[{\mathbf 1} ( \theta_{1} > m )\exp(\alpha\theta_{1})]\\
 & \leq (1/r) \text{ for $m$ sufficiently large}.
\end{align*}
Here the last inequality follows because $\E[\exp(\alpha \theta_1)] < \infty$ guarantees 
$\exp(-\alpha(m + 1))\E[{\mathbf 1}
( \theta_{1} > m )\exp(\alpha\theta_{1})] \to 0$ as $m \to \infty$. \qed

\medspace

\noindent{\bf Proof of Lemma \ref{lem:RandomSumHavingExpTail}}:
We have $ \P ( N = n ) \leq
\P( N \geq n ) \leq \E ( \alpha N ) \exp ( - \alpha n ) $.
Let $ \Psi $ be the moment generating function of $ \theta_1 $. Then, for all $ \gamma_0 \leq \beta $,
$ \Psi ( \gamma_0) = \E \bigl( \exp ( \gamma_0 \theta_1) \bigr) < \infty$.
Since the function $ \Psi  ( \gamma_0 ) $ is continuous at $0$ and $ \Psi ( 0 )
= 1 $, we may choose $ \gamma  > 0 $
so that $  \Psi ( 2 \gamma )  \exp ( -  \alpha  ) < 1 $. Now, we have
 % Since $\{J^{1}_n: n \geq 1\}$ are i.i.d. random variables, each
% with an exponentially decaying tail probability,  there exists $ \beta_0 > 0
% $ such that
% the moment generating function of $ \Psi_J
% ( \alpha ) := E(\exp (\alpha  J^{1}_{1}))  < \infty $ for all $
% \alpha < \beta_0 $.

\begin{align*}
\lefteqn{   \E \Bigl[  \exp ( \gamma S ) \Bigr]  =  \E \Bigl[ \sum_{n= 1}^{ \infty } {\mathbf 1} (N = n ) \exp  ( \gamma
\sum_{i=1}^n  \theta_i ) \Bigr]  } \\
 & =   \sum_{n= 1}^{ \infty }  \E \Bigl[  {\mathbf 1} ( N = n ) \exp (
\gamma  \sum_{i=1}^n  \theta_i ) \Bigr]
 \leq \sum_{n= 1}^{ \infty }   \bigl[ \P ( N = n ) \bigr]^{1/2} \bigl[ \E
 \bigl( \exp ( 2 \gamma  \sum_{i=1}^n  \theta_i )  \bigr) \bigr]^{1/2} \\
& \leq \sum_{n= 1}^{ \infty } \sqrt{ \E ( \alpha N ) } \exp ( - n \alpha / 2  )  \bigl[  \Psi (
2 \gamma )  \bigr]^{n/2}
= \sqrt{ \E ( \alpha N )  } \sum_{n= 1}^{ \infty }  \bigl[ \exp ( - \alpha  ) \Psi ( 2 \gamma ) \bigr]^{n/2} < \infty ;
\end{align*}
here the first inequality follows from the Cauchy-Schwartz
inequality. This completes the proof. \qed

\noindent{\bf Proof of Lemma \ref{lem:IndRegen} }: 
Define $ L_n := \max \{  R_n^{(1 )},  R_n^{(2 )}  \} $ and set
$ \tau^L := \inf \{ n \geq 1 : L_n = 0 \} $. Then, we have $ \tau^R = \tau^L $. Again, we define
a new Markov chain which dominates $ L_n $ and satisfies the conditions of Lemma \ref{lem:MCAuxRes},
from which we will conclude the result.

We start with $2$ families of independent copies of the inter-arrival times,
say $  \{ \eta_n^{(1)} : n \geq 1 \} $ and $  \{ \eta_n^{(2)} : n \geq 1 \} $
with $ \eta_1^{(i)} \stackrel{d}{=} \xi_1^{(i)} $
for  $ i =1, 2 $. Now keeping the same notation as in the proof of 
Proposition \ref{prop:RecTimeExpTail}, we set $ W_n^{\text{move}} := 
\{ i : R_n^{(i)} = 0 $ for $ i = 1, 2 \}$ and $ W_n^{\text{stay}} :=
\{ 1, 2\} \setminus W_n^{\text{move}} $. Now, for $ i \in W_n^{\text{move}} $, we have  $ S_{l_i(n)}^{(i)} = n $
for some $ l_i (n) \geq 0 $,  and, for $  i \in  W_n^{\text{stay}} $, we have $ S_l^{(i)} \neq n $ for every $ l \geq 0 $. Define
\begin{equation*}
J_{n+1} := \max \bigl\{ \max \{ \xi_{ l_i(n) + 1 }^{(i)} : i \in W_n^{\text{move}} \}, \max \{  \eta_{n+1}^{(i)} : i \in
W_n^{\text{stay}} \} \bigr\}
\end{equation*}
and
\begin{equation*}
M_0 := 0 \text{ and } M_{n+1} := \max \{ M_n , J_{n+1} \} - 1\text{ for } n \geq 0.
\end{equation*}

We now claim $ M_n \geq L_n $ for all $ n \geq 0 $. 
We have $ M_0 =  L_0 = 0 $. Assume that
the result holds for $ n $
% For $ i \in W_n^{\text{stay}}$ (i.e.,
% $ R_n^{(i)} \geq 1$) we have $ R_{n+1}^{(i)} = R_{n}^{(i)}  - 1 $.
% While, for $ i \in W_n^{\text{move}}$ (i.e.,
% $S_{l_i(n)}^{(i)} = n $ for some $ l_i(n) \geq 0 $)
% we have
% \begin{align*}
%  R_{n+1}^{(i)} & =  \inf \{  S_k^{(i)} : S_k^{(i)}  \geq n+1 \}
% - n - 1 =  S_{l_i(n)+1}^{(i)} - n - 1 \\
% & = S_{l_i(n)}^{(i)}+  \xi_{ l_i(n) + 1 }^{(i)} - n - 1
% =  \xi_{ l_i(n) + 1 }^{(i)} - 1.
% \end{align*}
and we have
\begin{equation*}
L_{n+1} = \max\{L_n , \max\{\xi_{ l_i(n) + 1 }^{(i)} :  i \in W_n^{\text{move}}\}\}-1
\leq \max\{M_n, J_{n+1}\}  -1 = M_{n+1}.
\end{equation*}
%  \begin{align*}
%   L_{n+1} & =  \max\{ R_{n+1}^{(i)}  : i = 1,2, \dotsc, k \}\\
% & = \max \bigl\{ \max \{  R_{n+1}^{(i)} :  i \in W_n^{\text{move}}  \},  \max \{  R_{n+1}^{(i)} :  i \in W_n^{\text{stay}} \} \bigr\}\\
% & =  \max \bigl\{ \max \{  \xi_{ l_i(n) + 1 }^{(i)} - 1 :  i \in W_n^{\text{move}}  \},
% \max \{  R_{n}^{(i)} - 1  :    i \in W_n^{\text{stay}} \} \bigr\}\\
% & \leq
% \max \bigl\{ \max \{  \xi_{ l_i(n) + 1 }^{(i)} :  i \in W_n^{\text{move}} \}, \max \{  R_{n}^{(i)}   :    i = 1,  \dotsc, k  \},\bigr.\\
% &   \qquad \qquad \max \{  \eta_{ n + 1 }^{(i)} :  i \in W_n^{\text{stay}}  \}  \bigr\} - 1\\
% & =  M_{n+1}.
% \end{align*}

% The independence of the families of random variables, $ \{ \xi_n^{(i)} \} $ and
% $ \{ \eta_n^{(i)} \} $ and the fact that $  \xi_1^{(i)} \stackrel{d}{=} \eta_1^{(i)}$, 
% for $ i = 1, 2$,
% implies that we can write $ M_{n+1} = \max \{M_n , \theta_{n+1} \} - 1 $ 
% where $ \{ \theta_n : n \geq 1 \} $
% is a sequence of i.i.d. random variables with $ \theta_1 \stackrel{d}{=} \max \{ \xi_1^{(i)} :
% i = 1, 2 \} $. 
The assumptions imposed on $ \xi_n^{(i)} $ imply that the Markov chain
satisfies the conditions of Lemma \ref{lem:MCAuxRes} and the result follows from that. \qed

\medspace

\noindent{\bf Calculations for (\ref{eqn:FirstMomentInd}),(\ref{eqn:SecondMomentInd}) and (\ref{eqn:ThirdMomentInd}):} 
For all $m \geq 1$ we have
\begin{align*}
& \E[(||(\overline{\bu^0} + \psi^{\bu^0}_1) - (\overline{\bv^0} + \psi^{\bv^0}_1)||^2_2 - ||{\rm x}||^2_2)^m]\\
& = \E[(||({\rm x}(1) + \psi^{\bu^0}_1(1) - \psi^{\bv^0}_1(1),{\rm x}(2) +
\psi^{\bu^0}_1(2) - \psi^{\bv^0}_1(2))||^2_2 - ||{\rm x}||^2_2)^m]\\
& = \E\bigl[\bigl((\psi^{\bu^0}_1(1))^2 + (\psi^{\bu^0}_1(2))^2 + (\psi^{\bv^0}_1(1))^2 +
(\psi^{\bv^0}_1(2))^2  - 2\psi^{\bu^0}_1(1)\psi^{\bv^0}_1(1) \\
& \qquad - 2\psi^{\bu^0}_1(2)\psi^{\bv^0}_1(2) + 2{\rm x}(1)(\psi^{\bu^0}_1(1) - 
\psi^{\bv^0}_1(1)) + 2{\rm x}(2)(\psi^{\bu^0}_1(2) - \psi^{\bv^0}_1(2))\bigr)^m\bigr].
\end{align*}
From Proposition \ref{prop:PropertyIncrementsForIndPaths} we have
$\E[(\psi^{\bu^0}_1(j_1))^{m_1}(\psi^{\bv^0}_1(j_2))^{m_2}] = 0$ at least one of $m_1,m_2$ is
odd. Hence for $m = 1$ we have,
\begin{align*}
& \E[||(\overline{\bu^0} + \psi^{\bu^0}_1) - (\overline{\bv^0} + \psi^{\bv^0}_1)||^2_2 - ||{\rm x}||^2_2]\\
& = \E[(\psi^{\bu^0}_1(1))^2 + (\psi^{\bu^0}_1(2))^2 + (\psi^{\bv^0}_1(1))^2 +
(\psi^{\bv^0}_1(2))^2] = 4\E[(\psi^{\bu^0}_1(1))^2].
\end{align*}
For $m = 2$ using Proposition \ref{prop:PropertyIncrementsForIndPaths} we have
\begin{align*}
& \E[(||(\overline{\bu^0} + \psi^{\bu^0}_1) - (\overline{\bv^0} + \psi^{\bv^0}_1)||^2_2 - ||{\rm x}||^2_2)^2]\\
& = \E[(((\psi^{\bu^0}_1(1))^2 + (\psi^{\bu^0}_1(2))^2 + (\psi^{\bv^0}_1(2))^2 +
(\psi^{\bv^0}_1(2))^2)  - 2\psi^{\bu^0}_1(1)\psi^{\bv^0}_1(1) \\
& \qquad -2 \psi^{\bu^0}_1(2)\psi^{\bv^0}_1(2)
 + 2{\rm x}(1)(\psi^{\bu^0}_1(1) - \psi^{\bv^0}_1(1)) + 2{\rm x}(2)(\psi^{\bu^0}_1(2) - \psi^{\bv^0}_1(2)))^2]\\
& \geq \E[(2{\rm x}(1)(\psi^{\bu^0}_1(1) - \psi^{\bv^0}_1(1)))^2 + (2{\rm x}(2)(\psi^{\bu^0}_1(2) - \psi^{\bv^0}_1(2)))^2]\\
& = 4({\rm x}(1))^2\E[(\psi^{\bu^0}_1(1))^2 + (\psi^{\bv^0}_1(1))^2] +
4({\rm x}(2))^2\E[(\psi^{\bu^0}_1(2))^2 + (\psi^{\bv^0}_1(2))^2]\\
& = 8||{\rm x}||^2_2\E[(\psi^{\bu^0}_1(1))^2].
\end{align*}
The inequality follows from the fact that 
$\E[(\psi^{\bu^0}_1(j_1))^{m_1}(\psi^{\bv^0}_1(j_2))^{m_2}] \neq 0$ for all $1\leq j_1, j_2 \leq d-1$, only if both $m_1$ and $m_2$ are
even and $\E[\psi^{\bu^0}_1(1)\psi^{\bu^0}_1(2)] = \E[\psi^{\bv^0}_1(1)\psi^{\bv^0}_1(2)] = 0$.

By the same logic it also follows that for $m = 3$ we have
\begin{align*}
& \E[(||(\overline{\bu^0} + \psi^{\bu^0}_1) - (\overline{\bv^0} + \psi^{\bv^0}_1)||^2_2 - ||{\rm x}||^2_2)^3]\\
& = \E[(((\psi^{\bu^0}_1(1))^2 + (\psi^{\bu^0}_1(2))^2 + (\psi^{\bv^0}_1(2))^2 +
(\psi^{\bv^0}_1(2))^2)  - 2\psi^{\bu^0}_1(1)\psi^{\bv^0}_1(1) \\
& \qquad -2 \psi^{\bu^0}_1(2)\psi^{\bv^0}_1(2)
 + 2{\rm x}(1)(\psi^{\bu^0}_1(1) - \psi^{\bv^0}_1(1)) + 2{\rm x}(2)(\psi^{\bu^0}_1(2) - \psi^{\bv^0}_1(2)))^3]\\
& = 12\E[((\psi^{\bu^0}_1(1))^2 + (\psi^{\bu^0}_1(2))^2 + (\psi^{\bv^0}_1(1))^2 +(\psi^{\bv^0}_1(2))^2
-2(\psi^{\bu^0}_1(1)\psi^{\bv^0}_1(1) + \\
& \qquad \psi^{\bu^0}_1(2)\psi^{\bv^0}_1(2)))(({\rm x}(1)(\psi^{\bu^0}_1(1) -
\psi^{\bv^0}_1(1)))^2 + ({\rm x}(2)(\psi^{\bu^0}_1(2) - \psi^{\bv^0}_1(2)))^2) + \\
& \qquad \text{ terms free of }{\rm x}]\\
& = 24||{\rm x}||^2_2\E[(\psi^{\bu^0}_1(1))^4 + (\psi^{\bu^0}_1(1)\psi^{\bu^0}_1(2))^2 +
4(\psi^{\bu^0}_1(1)\psi^{\bv^0}_1(1))^2] + \text{ terms free of }{\rm x} \\
& = O ( || {\rm x}||_2^2) \text{ as } || {\rm x}||_2 \to \infty.
\end{align*}
\medspace

\noindent{\bf Proof of Lemma \ref{lemma:fnconv}}:
Let $i \in \{1,\ldots, k-1 \}$ be such that $\pi_i(t^k) = \pi_k(t^k)$ and
$\pi_j(t^k) \neq \pi_k(t^k)$ for all $1 \leq j < i$.
Fix $\epsilon$ such that $0 < \epsilon < t^k - \max\{\sigma_{\pi_i}, \sigma_{\pi_k}\}$.
Given $\eta > 0$ let $P_i, P_k \subseteq \R^2$ be defined as
\begin{align*}
 P_i(\eta) & = \{(x,u): ||(x,u) - (\pi_i(s), s)||_1 \leq \eta \text{ for some } \sigma_{\pi_i} \leq s \leq t^k - \epsilon\}\\
 P_k(\eta) & = \{(x,u): ||(x,u) - (\pi_k(s), s)||_1 \leq \eta \text{ for some } \sigma_{\pi_k} \leq s \leq t^k - \epsilon\},
\end{align*}
i.e. $P_i$ and $ P_k$ are the regions obtained by $\eta$-fattening the paths $\pi_i$ and $\pi_k$ respectively. Since 
$t^k = \inf \{s : \pi_i(s) = \pi_k(s)\} > \max\{\sigma_{\pi_i}, \sigma_{\pi_k}\}$ therefore 
we may first choose $0 < \eta < \epsilon/2 $ such that
$$d(P_i(\eta), P_k(\eta)):= \inf\{||(x,u) - (y,v)||_1 : (x,u) \in P_i(\eta), (y,v) \in P_k(\eta)\} > \eta.$$

Next, since $d_{\Pi}(\pi_{i,n}, \pi_i) \to 0$ and  $d_{\Pi}(\pi_{k,n}, \pi_k) \to 0$ as $n \to \infty$  
we may choose $n_0 \geq 1$ such that $ \eta > n_0^{\alpha - 1} $ and, for all $n \geq n_0$, the following hold:
\begin{itemize}
 \item [(a)] $\sigma_{\pi_{i,n}} \leq t^k - \epsilon$ and $ \sigma_{\pi_{i,n}} \leq t^k - \epsilon$,
 \item [(b)] $\{(\pi_{i,n}(s), s): \sigma_{\pi_{i,n}} \leq s \leq  t^k - \epsilon \}  \subseteq P_i(\eta)$ and 
 $ \{(\pi_{k,n}(s), s): \sigma_{\pi_{k,n}} \leq s \leq  t^k - \epsilon \}  \subseteq P_k(\eta)$.
\end{itemize}
Since $P_i$ and $ P_k$ are separated by a minimum distance $\eta$, we have $ t^k_n \geq t^k - \epsilon
$ for all $n\geq n_0$ and hence $\liminf_{n \to \infty} t^k_n  \geq t^k - \epsilon $.

Now assume that $\pi_k(t^k - \epsilon)
> \pi_i(t^k - \epsilon)$.
For the other case the argument is exactly similar.
Fix $ s \in [t^{k},t^{k} + \epsilon] $, such that
$ \pi_i(s)-\pi_k(s) = \nu > 0 $.  
For $n_0$ as above, choose $n_1 > n_0 $ such that for all $n\geq n_1$ we have
$\sup_{t\in [t^k - \epsilon,t^k + \epsilon]}|\pi_{j,n}(t)-\pi_j(t)| <\nu/4$ for $j=i,k$.
For $n > n_1 $, $ \pi_{i,n} (s) - \pi_{k,n} (s) \geq
 \pi_i(s)-\pi_k(s) -
| \pi_{i,n} (s) - \pi_{i} (s) | -  | \pi_{k,n} (s) - \pi_{k} (s) |  > \nu / 2 > 0 $; and
our choice of $n_1$ ensures that $\pi_{k,n}(t^k - \epsilon) - \pi_{i,n}(t^k - \epsilon)
 >  0 $.
Thus, $\pi_{i,n}$ and $\pi_{k,n}$ cross each other before time
$t^{k} + \epsilon$ and hence $\limsup_{n\to \infty}
t^{k}_n \leq t^{k} + \epsilon$.
This completes the proof of first part of the Lemma.

Since $t^k_n \to t^k$ as $n \to \infty$,  $d_{\Pi}(\pi_i, \pi_{i,n}) \to 0$
and $d_{\Pi}(\pi_k, \pi_{k,n}) \to 0$ as $n \to \infty$, to show 
$d_{\Pi}^k \Bigl( f_n^{(\alpha)} (\pi_{1,n}, \dotsc , \pi_{k,n}), f (\pi_{1}, \dotsc ,  \pi_{k} ) \Bigr) \to 0$ as 
$n\rightarrow \infty$ it suffices to how that $\sup_{t\in [t^k - \epsilon,t^{k}+\epsilon]}
| \overline{\pi}_{k,n}(t)- \overline{\pi}_k(t)| \to 0$ as $n\to \infty$.

For $0 < \beta < \epsilon$, writing 
\begin{align*}
 \lefteqn { \sup_{t\in [t^k - \epsilon,t^k + \epsilon]}| \overline{\pi}_{k,n}(t)- \overline{\pi}_k(t)|   \leq
\sup_{t\in [t^k - \epsilon,t^k - \beta]}|\pi_{k,n}(t)-\pi_k(t)|  } \\
& \qquad + \sup_{t\in [t^k - \beta,t^k + \beta]}| \overline{\pi}_{k,n}(t) - \overline{\pi}_k(t)|
+ \sup_{t\in [t^k + \beta,t^k +\epsilon]}|\pi_{i,n}(t)-\pi_i(t)| .
\end{align*}
and observing that 
\begin{itemize}
\item[(a)] the first and the last terms of the expression above can be made arbitrarily small
as in the first part of this proof,
\item[(b)] the middle term can be made small by choosing $\beta$  such that,  for each of $j_1, j_2 \in \{i,k\}$,
$\sup_{s_1,s_2\in [t^k-\beta,t^k+\beta]}|\pi_{j_1}(s_1)-\pi_{j_2}(s_2)| $
is small 
and  noting that  $ \overline{\pi}_{k,n}$ is defined by a linear interpolation
between  $\pi_{k,n}(t^k_n)$ and
$\pi_{i,n}(\frac{\lfloor n^2 \gamma t^k_n \rfloor +1 }{n^2 \gamma})$.
\end{itemize} \qed

Finally, to show (\ref{Kumarjitdekho}) 
% we have the following lemma.
% \begin{lemma}
% \label{lem:Kumarjitdekho}
% Let $D_n$ and $D$ be as defined in the proof of Proposition \ref{prop:randomBM}. 
% Further let $D_n$ be such that, given the starting point set $D_n$, the collection of scaled 
% paths starting from $D_n$ is independent of $\sigma(D_n)$. Then, the scaled paths starting from
% $ D_n $ converges weakly to $ {\cal W}^D $, i.e.,
% \begin{equation*}
% \chi_n^{D_n} \Rightarrow {\cal W}^D \text{ in } ({\cal H}, d_{{\cal H}})
% \end{equation*}
% where $ {\cal W}^D $ is the set of paths given by
% coalescing Brownian motions starting from a random point set distributed as $D$.
% \end{lemma}
% 
% \noindent\textbf{ Proof of \ref{lem:Kumarjitdekho}:}
we first show that for any deterministic finite sets $B_n$ and $B$ with 
$B_n\subset \Z^2$, $B\subset \R^2$ 
$\rho_{{\cal P}}(B_n^{(\text{scaled})},B)\to 0$ as $n \to \infty$
where $B_n^{(\text{scaled})}  := \{(\by(1)/(n \sigma), \by(2)/(n^2\gamma)): \by \in B_n\}$, 
we have ${\cal X}^{B_n}_n$ converges weakly 
to ${\cal W}^B$, i.e.,  coalescing Brownian motions starting from a random point set distributed as $B$.
Since almost surely ${\cal X}^{\Z^2}$ consists of noncrossing paths
only, $(I_1)$ implies that the family $\{\overline{{\cal X}^{\Z^2}_n} : n \in \N\}$ is
tight, and  ${\cal X}^{B_n}_n \subset \overline{{\cal X}^{\Z^2}_n}$ guarantees that
$\{{\cal X}^{B_n}_n : n \in \N\}$ is also tight.
The sequence $\{\overline{{\cal X}^{\Z^2}_n} : n \in \N \}$ also 
satisfies $(I_1)$ and hence satisfies $(B_1)$. The proof of Theorem 5.3 in \cite{FINR04}
shows that  for any subsequential limit ${\cal Z}$ of $\{\overline{{\cal X}^{\Z^2}_n} : n \in \N \}$
and for any deterministic $\bx \in \R^2$ there is almost surely a unique path starting from $\bx$ in
${\cal Z}$.
A coupling argument then shows that the same is true for any subsequential limit ${\cal Z}_B$
of $\{{\cal X}^{B_n}_n : n \in \N \}$.
% and for any ${\rm x} \in D$, ${\cal Z}_D$ has unique path starting from ${\rm x}$ almost surely.
The sequence $\{(\overline{{\cal X}^{\Z^2}_n}, {\cal X}^{B_n}_n ): n \in \N \}$ 
is jointly tight and let $({\cal Z}, {\cal Z}_B)$ be a subsequential limit of this sequence.
By Skorohod's representation theorem we assume that we are working on a probability space
such that $\{(\overline{{\cal X}^{\Z^2}_{n_k}}, {\cal X}^{B_{n_k}}_{n_k}) : k \in \N \}$ 
converges almost surely to $({\cal Z}, {\cal Z}_B)$. Since 
${\cal X}^{B_{n_k}}_{n_k} \subseteq \overline{{\cal X}^{\Z^2}_{n_k}}$ for all ${n_k}$, 
if for any deterministic $\bx \in B$,  with positive probability
${\cal Z}_B$ has more than one path starting from $\bx$  then
so does ${\cal Z}$. Hence for all $\bx \in B$, ${\cal Z}_B$ has unique path starting from $\bx$
almost surely. Now by $(I_1)$ the finite dimensional distributions of ${\cal Z}_B$ 
are the same as that of a process
of a coalescing Brownian motions. Therefore we have that ${\cal Z}_B$ has the same distribution
as ${\cal W}^B$ starting from the  set $B$.

For the general case, it suffices to show that 
$\E[f({\cal X}^{D_n}_n)] \to \E[f({\cal W}^{D})]$
as $n \to \infty$ for all bounded continuous $f$ on $({\cal H}, d_{\cal H})$.
Let $f_n(D_n) := \E[f({\cal X}^{D_n}_n) | D_n] $ and 
$f(D) := \E[f({\cal W}^{D})|D]$.
By Skorohod's representation theorem we can assume that we are working on a probability
space such that $D_n \to D$ almost surely as $n \to \infty$ in $({\cal P}, \rho_{{\cal P}})$.
Let $\{U^a_{\bw}:\bw \in \Z^2\}$ be a collection of i.i.d. $U[0,1]$ random variables
and independent of the collection $\{U_{\bw}:\bw \in \Z^2\}$ used to build the model.
For any $A \subseteq \Z$, let ${\cal X}_{a,n}^A$
be the collection of scaled paths starting from $A$ constructed using $\{U^a_{\bw}:
\bw \in \Z^2\}$ only. 
Since the evolution of the paths from $D_n$ is independent of $\sigma(D_n)$, we have $$
\chi_n^{D_n} | D_n \stackrel{d}{=} \chi_{a,n}^{D_n} \text{ almost surely}.
$$
From our assumptions on $D_n$ we have 
$f_n(D_n) := \E[f({\cal X}^{D_n}_n) | D_n] = 
\E[f({\cal X}^{D_n}_{a,n})]$ almost surely.
Then, for almost every $\omega$, by the deterministic part of this proof we have that
${\cal X}^{D_n(\omega)}_{a,n}$ converges in distribution to ${\cal W}^{D(\omega)}$.
Hence we have almost surely $f_n(D_n) \to f(D)$ as $n \to \infty$.
By the bounded convergence theorem we have $\E[f_n(D_n)] = \E[f({\cal X}^{D_n}_n)] \to
\E[f(D)] = \E[f({\cal W}^{D})]$ as $n \to \infty$.
\qed

\noindent
{\bf Acknowledgements:} We thank the referee for his comments which led to
a significant improvement of this paper. Kumarjit Saha is grateful to Indian Statistical Institute for
a fellowship to pursue his Ph.D.

\vspace{.5cm}
\noindent Rahul Roy, Kumarjit Saha and Anish Sarkar \\
Theoretical Statistics and Mathematics Unit\\
Indian Statistical Institute\\
7 S. J. S. Sansanwal Marg\\
New Delhi 110016, India.\\
{\em rahul, kumarjit10r, anish@isid.ac.in}

\end{document}